\documentclass{amsart}

\usepackage{enumerate}
\usepackage{amsmath}
\usepackage{amssymb}
\usepackage{amsthm}

\usepackage{graphicx}
\usepackage{type1cm}
\usepackage{eso-pic}
\usepackage{color}
\makeatletter
\usepackage{math}
\usepackage{thm}
\usepackage{ref}

\newcommand{\R}{\ensuremath{\mathbb{R}}}

\newcommand{\D}{\ensuremath{\mathcal{D}}}
\newcommand{\N}{\ensuremath{\mathbb{N}}}

\newcommand{\diam}{\operatorname{diam}}

\newcommand{\eps}{\ensuremath{\varepsilon}}

\newcommand{\Div}{\operatorname{div}}
\newcommand{\DIV}{\operatorname{Div}}

\newcommand{\E}{\mathbb{E}}
\newcommand{\F}{\mathbb{F}}

\newcommand{\pa}{\partial}
\newcommand{\grad}{\operatorname{grad}}
\newcommand{\Def}{\operatorname{Def}}

\newcommand{\id}{\operatorname{id}}

\renewcommand{\span}{\operatorname{span}}
\renewcommand{\phi}{\varphi}
\newcommand{\vertiii}[1]{{\left\vert\kern-0.25ex\left\vert\kern-0.25ex\left\vert #1 
    \right\vert\kern-0.25ex\right\vert\kern-0.25ex\right\vert}}
\newcommand{\esssup}{\operatorname{ess\,sup}}

\title
	[On a Stokes-type system arising in fluid vesicle dynamics]
	{On a Stokes-type system arising in fluid vesicle dynamics}

\author
	[Daniel Lengeler]
	{Daniel Lengeler}

\address
	{Department of Mathematics \newline\indent
	 Universit{\"a}t Regensburg --
	 Universit{\"a}tsstr.~31, 93053 Regensburg, Germany}
	
\email
	{daniel.lengeler@mathematik.uni-regensburg.de}

	
\thanks{I would like to express my gratitude to Harald Garcke, University of Regensburg, and Udo Seifert, University of Stuttgart, for numerous helpful comments and fruitful and pleasant discussions on the physics and mathematics of lipid bilayer membranes. Furthermore, I would like to thank Helmut Abels, University of Regensburg, for very helpful discussions on the theory of function spaces. Finally, I gratefully acknowledge support by DFG SPP 1506 "Transport Processes at Fluidic Interfaces".}		


\date
	{\today}

\begin{document}
% \setlength{\parskip}{0.5\baselineskip}
% \setlength{\parindent}{0pt}
% \renewcommand{\baselinestretch}{1.125}
% \normalsize
\begin{abstract}
This article is the first in a series of papers on the analysis of a basic model for fluid vesicle dynamics. There are two variants of this model, a parabolic one decribing purely relaxational dynamics and a non-parabolic one containing the full dynamics. At the heart of both variants lies a linear elliptic system of Stokes-type. Understanding the mapping properties of this Stokes-type system is crucial for all further analysis. In this article we give a basic exposition of the dynamical model and a thorough $L_2$-analysis of the Stokes-type system that takes into account geometric variations of the fluid vesicle.
\end{abstract}
\maketitle

\section{Introduction}

The basic constituent of most biological membranes is a two-layered sheet of phospholipid molecules, a \emph{lipid bilayer}. Phospholipid molecules possess hydrophilic heads and hydrophobic tails. When exposed to water they arrange themselves into a two-layered sheet with tails pointing inward. Due to the hydrophobic effect, biological membranes tend to avoid open edges and form closed configurations, the \emph{vesicles}. While a membrane is only a few nanometers thick, the diameter of the vesicles formed can be up to $10^4$ times larger. This separation of length scales suggests to describe vesicles as two-dimensional surfaces embedded in three-dimensional space. Since the solubility of the phospholipids is very low, the exchange of material between the membrane and the ambient solution is negligible. Thus, vesicle configurations are not determined by a surface tension but rather by a bending elasticity leading to an amazing variety of equilibrium shapes; in contrast to the typically spherical equilibrium 
shapes of fluid interfaces. A basic model for the elastic energy, the \emph{Canham-Helfrich energy}
\begin{equation}\label{eqn:ch}
 F=\frac\kappa 2\int_\Gamma (H-C_0)^2\,dA+\kappa_G\int_\Gamma K\,dA,
\end{equation}
 was proposed independently in three seminal papers \cite{canham,helfrich,evans}. Here, $\Gamma$ is the two-dimensional surface representing the membrane, $H$ and $K$ denote twice its mean curvature and its Gauss curvature, respectively, and $\kappa,\kappa_G$ are the bending rigidities with the physical dimension of an energy. Furthermore, the constant $C_0$ is the \emph{spontaneous curvature} which is supposed to reflect a chemical asymmetry of the membrane or its environment. The Canham-Helfrich energy (with $C_0=0$) can be seen as the lowest order term in an expansion in curvature. By the Gauss-Bonnet theorem, the Gaussian part of the energy is a topological constant and hence irrelevant for vesicles of fixed topology. Due to the above mentioned low solubility of the phospholipids and a low compressibility of the membrane its area is practically fixed. Furthermore, despite the fact that the membrane is permeable to water, the fluid volume enclosed by the vesicle is also fixed, due to strong osmotic 
effects.
 
Thus, in the simplest case, equilibrium vesicle shapes supposedly minimize the Canham-Helfrich energy under the constraints of fixed area and enclosed volume. The known solutions to this minimization problem can indeed be observed experimentally. However, other (non-axisymmetric) shapes found in the laboratory like the starfish vesicles and shape transformations like the budding transition seem not to be contained in this model.
% \footnote{Furthermore, in fact, the experimental observation of \emph{conformal diffusion} rules out a non-vanishing spontaneous curvature, see, e.g., [].} 
The reason for this deficiency is that the bilayer architecture of the membrane is not correctly incorporated. At a characteristic temperature all lipids  transition from a gel to a liquid phase. The lipid bilayers we are dealing with are assumed to be in the liquid phase, allowing the monolayers to freely flow laterally and to slip over each other while the membrane retains its transverse structure. This leads to a lateral adjustment of the monolayer densities to curvature. Incorporating the bilayer architecture in the form of inhomogeneous monolayer densities yields, to lowest order, the additional elastic energy
\begin{equation}\label{eqn:sl}
 F_{SL}=\frac{\kappa_M} 2\int_\Gamma\Big( (\rho^+-dH)^2+(\rho^-+dH)^2\big)\, dA,
\end{equation}
introduced in \cite{seifertlanger}. Here, $\rho^{\pm}$ are the \emph{reduced density deviations} of the monolayers and $d$ is the membrane thickness. By scale considerations one can show that in equilibrium the mean density $(\rho^+ +\rho^-)/2$ is a fixed constant leading again to a fixed membrane area while the density difference $(\rho^+ -\rho^-)/2$ gives a nonlocal energy involving the total mean curvature. Together with the volume constraint this yields the area-difference elasticity model introduced independently in \cite{seifert92,wiese,bozic}. This model seems to be quantitatively consistent with all experiments performed so far. An intermediate model, the bilayer couple model, introduced in \cite{svetina82,svetina83,svetina89}, is based on the Canham-Helfrich energy, but takes into account three hard constraints, namely area, enclosed volume, and total mean curvature. 

%The typical time scale for micrometer scale shape is of the order of seconds, making them observable under the microscope in real time.
In most experimental situations the dynamics of fluid vesicles are of interest, too. On the one hand, thermal excitation always leads to fluctuations around stable equilibrium shapes. On the other hand, real biological membranes are constantly put into non-equilibrium states by biophysical processes, like molecule insertions and extractions, or by external hydrodynamic flows. When tweaking a parameter (like the enclosed volume, e.g. by changing the osmotic conditions) discontinuous shape transitions may occur, that is, a formerly stable shape becomes unstable and decays towards a new minimum. In any case, it is crucial to take into account the interaction of the moving membrane with the viscous bulk fluid. Apart from the bulk viscosity, two other relevant dissipation 
mechanisms have been identified, namely friction between the two monolayers and shear viscosity within each layer; see for instance \cite{seifertlanger,arroyo10}. In most situations of interest both the bulk and the surface Reynolds numbers are very small, typically of the order $10^{-3}$-$10^{-5}$. In such situations it makes sense to neglect inertial effects in the bulk and on the membrane, leading to purely relaxational dynamics. The significance of the bulk fluid then lies in the induced nonlocal self-interaction of the membrane and in the viscous energy dissipation.
% However, the bulk Reynolds number being of the order of $10^{-6}$ to $10^{-3}$, inertial effects in the bulk can be safely neglected. In fact, the significance of the bulk fluid lies in the induced nonlocal selfinteraction of the membrane and in the viscous energy dissipation.  For the majority of applications it is sufficient to restrict to long-wavelength phenomena. However, due to viscous dissipation, these are governed by overdamped dynamics. Therefore, we can safely neglect membrane inertia, too.\marginpar{So in Ordnung? Nachlesen!}

Several models for the dynamics of fluid vesicles have been proposed; see for instance \cite{schneider,seifertlanger,foltin,miao02,miao05,arroyo09,arroyo12}. Some of these are restricted to the minimal energy \eqref{eqn:ch} while others take into account the bilayer architecture using the additional energy \eqref{eqn:sl}. While \cite{seifertlanger} deals only with an almost planar membrane, \cite{schneider,foltin,miao02,miao05} neglect in-plane shear viscosity. Including the latter in an arbitrary geometry leads to the well-known Boussinesq-Scriven surface fluid \cite{boussinesq,scriven} as basic continuum mechanical modeling shows; cf. \cite{arroyo09,arroyo12}. In this paper, we take into account the in-plane shear viscosity while neglecting the bilayer architecture of the membrane. More precisley, we study a single homogeneous Newtonian surface fluid subject to additional stresses induced by the Canham-Helfrich energy, interacting with a homogeneous Newtonian bulk fluid. The full model (and intermediate 
models) will be studied in future publications. For a comprehensive presentation of the physics of fluid vesicles we refer the reader to \cite{seifert97} and the references therein.

\emph{Phase diagrams} of the different curvature models assign the shape of lowest energy to each choice of external parameters (like the enclosed volume). Most of the insight into these phase diagrams and into the bifurcation and stability of stationary shapes has been generated by numerical computations and formal (linear) analysis; again, see \cite{seifert97} and the references therein. From the point of view of a rigorous mathematical analysis, not much is known concerning these questions. Critical points of the energy \eqref{eqn:ch} subject to the volume and area constraint solve the Helfrich equation
\begin{equation}\label{eqn:helfrich}
 \grad_{L_2}F+\lambda_1+\lambda_2H=0
\end{equation}
with Lagrange parameters $\lambda_1,\lambda_2$. While the round sphere is a solution of this equation for any choice of $C_0$, it is the only known analytical solution of spherical topology for $C_0=0$; see \cite{ouyang2} for further analytical solutions of spherical topology in the case of non-vanishing spontaneous curvature. In \cite{nagasawa03} the authors prove the existence and analyze the stability properties of a one-parameter family of critical points bifurcating from the round sphere. In \cite{schygulla} it is shown that there exists a global minimizer of the Canham-Helfrich energy with vanishing spontaneous curvature in the class of smoothly embedded topological spheres with fixed \emph{isoperimetric ratio} $\sigma\in (0,1]$.\footnote{Since the Canham-Helfrich energy is scaling invariant, instead of fixing both the enclosed volume and the area, it is sufficient to fix only one parameter, for instance the isoperimetric ratio.} Furthermore, it is shown that the minimal energy is a 
continuous and strictly decreasing function of $\sigma$, approaching $8\pi$ from above as $\sigma\searrow 0$ and that, roughly speaking, the minimizers converge to a double sphere in a measure-theoretic sense in this limit. The existence result is generalized in \cite{riviere} to higher genus surfaces. Experiments indicate that all minimizers are axi-symmetric (even for non-vanishing spontaneous curvature); however, proving this mathematically seems to be a hard task. Concerning the more realistic models, nothing is known from the point of view of a rigorous mathematical analysis. Insight into the dynamics of fluid vesicles based on numerical computations and formal linearized analysis can be found, for instance, in \cite{schneider,seifertlanger,foltin,seifert94,seifert99,miao02,arroyo09,arroyo10,arroyo12,garcke,garcke2,garcke3}; it seems that \cite{garcke} is the first article presenting fully three-dimensional numerical computations for the model we are studying in the present article. However, again, not much rigorous analysis has been done on this topic. Concerning the \emph{Canham-Helfrich flow}, that is, the $L_2$-gradient flow of the Canham-Helfrich energy with prescribed enclosed volume and area, a partial local well-posedness result has been shown in \cite{nagasawa12}. There exist further results \cite{nagasawa06,wheeler1,wheeler2} concerning a Helfrich-type flow where the Lagrange parameters instead of volume and area are prescribed and which consequently should not be related directly to fluid vesicles. In \cite{Zhang} local-in-time existence and uniqueness for a homogeneous Newtonian surface fluid subject to Canham-Helfrich stresses is shown. While the bulk fluid is neglected the authors keep the inertial term in the equations for the surface fluid, yielding a kind of dissipative fourth order wave-type equation. In \cite{Shkoller} local-in-time existence and uniqueness of a homogeneous Newtonian bulk fluid with inertial term interacting with a compressible, inviscid surface fluid without inertial term is shown in the $L_2$-scale, the membrane model being rather non-standard. 
%Since the authors work in the $L_2$ scale they have to deal with solutions of higher regularity, making the 
%analysis rather involved. Furthermore, they work in the 
%Lagrangian picture, leading to problems with the tangential degeneracy of the elliptic operator arising from the elastic stresses within the membrane.

The present article is the first in a series of papers on the relaxational and oscillatory (that is, involving both inertial terms) dynamics of a homogeneous surface fluid subject to Canham-Helfrich stresses and interacting with a homogeneous bulk fluid. 
%We will see that the relaxational dynamics can be interpreted as a gradient flow of the Canham-Helfrich energy, however, with respect to a Riemannian metric different from the $L_2$-metric (the latter yielding the classical Willmore and Canham-Helfrich flows, respectively). In fact, the corresponding Riemannian metric is induced by a linear Stokes-type system. 
It contains a basic exposition of the dynamical models and a thorough $L_2$-analysis of a linear Stokes-type system that plays a fundamental role in these models. Understanding the mapping properties of this system is an essential first step for all further work on topics like local well-posedness, stability, and long-time behavior. We will see that our Stokes-type system induces a Riemannian metric on the manifold of surfaces (of fixed area and enclosed volume) that is part of a gradient flow structure for the purely relaxational dynamics; this observation will prove very useful for stability analysis. For most further investigations, in fact, it is important to establish an $L_2$-theory that holds uniformly with respect to suitably controlled membrane deformations. This is particularly important and delicate for the analysis of the formation of singularities. It is a well-known and obvious fact that the Canham-Helfrich energy is invariant under isotropic scaling which is a continuous, non-compact symmetry. 
Although this symmetry is broken by the constraints of fixed area and enclosed volume, one can expect that a blow-up of curvature can only be ruled out by proving or assuming non-concentration of energy; see for instance \cite{Tao06} for the general heuristics on energy-critical systems.\footnote{Here, additionally, we need to assume or prove (which would be an interesting task) that no self-contact of the membrane occurs. Furthermore, the Canham-Helfrich energy, involving only the mean curvature, lacks coercivity. Thus, in order to rule out blow-up of curvature, we will have to control the integral of the second fundamental form squared instead of the mean curvature squared.} A similar result has been proven for the Willmore flow in \cite{kuwertschaetzle02}. The proof of such a fact relies on higher order energy estimates; hence the need for an $L_2$-theory of the underlying Stokes-type system that involves a precise control of the dependence on membrane deformations. 

The present article is organized as follows. In \Secref{stat} we 
shall derive the mathematical models and perform 
some basic computations; in particular, we present the gradient flow structure underlying the purely relaxational dynamics. In \Secref{stokes} we perform the $L_2$-analysis of the Stokes-type system for smooth, fixed membranes and then for variable membranes. In \Appref{app:int} we recall some basic facts concerning the kinematics and dynamics of fluid interfaces, and we derive the stresses resulting from the Canham-Helfrich energy. In \Appref{app:cov} we present some useful formulae concerning covariant differentian and curvature, while in \Appref{app:spaces} we prove some facts on the function spaces we shall use.

Before we proceed, let us fix some notation. Throughout the article (apart from the appendix), we denote by $e$ the Euclidean Riemannian metric in $\R^3$, and for surfaces $\Gamma\subset\R^3$ we denote by $g$ the metric on $\Gamma$ induced by $e$. We also use the notation $u\cdot v$ instead of $\langle u,v\rangle_e$ for $u,v\in\R^3$. We denote by $P_\Gamma$ the field of orthogonal projections onto the tangent spaces of $\Gamma$ and by $T\Gamma$ the tangent bundle of $\Gamma$. Furthermore, $[u]_\Gamma$ denotes the trace of the bulk field $u$ on $\Gamma$; however, when there is no danger of confusion we will sometimes omit the brackets. Moreover, we write $k$, $H$, and $K$ for the second fundamental form, twice the mean curvature, and the Gauss curvature of $\Gamma$ with respect to $e$, respectively. With a slight abuse of notation we use same symbol $k$ also to denote the Weingarten map, that is, in coordinates we write $k_{\alpha\beta}$ and $k_\alpha^\beta$. Furthermore, for any metric $\tilde e$ on an 
arbitrary manifold, we write 
$^{\tilde e}\Gamma_{ij}^k$, $\nabla^{\tilde e}$, $\grad_{\tilde e}$, $k_{\tilde e}$, etc. for the associated Christoffel symbols,
 differential operators, and curvature terms, and we use the abbreviations $\Gamma_{ij}^k:={^e\Gamma}_{ij}^k$, $\nabla:=\nabla^e$, $\grad:=\grad_{e}$, etc. for the corresponding Euclidean objects. When working in coordinates and confusion about the underlying metric can be ruled out, we use the semicolon to separate indices coming from covariant differentiation from the original indices; for instance, for a covector field $\omega$ we write $(\nabla^{\tilde e} \omega)_{ij}=\omega_{i;j}$. Again with a slight abuse of notation, we write $\langle T_1 ,T_2\rangle_{\tilde e}$ for the contraction of tensor fields $T_i$ of the same rank using the metric $\tilde e$ and $|T|_{\tilde e}$ for the corresponding norm. For manifolds $N,M$, a diffeomorphism $\varphi:N\rightarrow M$, and tensor fields $T$ on $N$ and $S$ on $M$ we denote by $\phi_*T$ the pushforward of $T$ and by $\phi^*S$ the pullback of $S$ under $\varphi$. 
% \begin{equation*}
% \begin{aligned}
% (\varphi^*u)^j(x)=(\pa_i(\varphi^{-1})^j\,u^i)(\varphi(x)),\\
% (\varphi_*U)^j(x)=(\pa_i\varphi^j\,U^i)(\varphi^{-1}(x))
% \end{aligned}
% \end{equation*}
We denote by $r(a)$ generic tensor fields that are polynomial or analytic functions of their argument $a$ such that $r(0)=0$. Furthermore, for tensor fields $r_1$ and $r_2$ we write $r_1*r_2$ for any tensor field that depends in a bilinear way on $r_1$ and $r_2$, and we use the abbreviations $r*(r_1,\ldots,r_k)=r*r_1 + \ldots + r*r_k$ and $r^k=r*\ldots*r$ (with $k$ factors on the right hand side). For $s>0$, $p,q\in [1,\infty]$, and tensor fields $T$ we denote by $\|T\|_p$, $|T|_{s,p}$, and $|T|_{s,p,q}$ the corresponding Lebesgue, Sobolev-Slobodetskij, and Besov (semi-)norms, see \Appref{app:spaces}; we omit the spatial 
domain in the notation since it will always be identical to the domain of definition of the argument $T$. We denote by $H^s$ and $H^s_0$ the $L_2$-scale of Sobolev-Slobodetskij spaces, where the lower index $0$ refers to vanishing traces (for $s>1/2$). 
%and $r((\bar\nabla^{\tilde e})^k a)=r(\nabla^{\tilde e} a,(\nabla^{\tilde e})^2a,\ldots,(\nabla^{\tilde e})^ka)$ for any metric $\tilde e$ and tensor field $a$.

\section{Modeling and basic observations}\seclabel{stat}

\subsection{The mathematical models}
Let $\Omega$ be a bounded domain in $\R^3$ containing a homogeneous Newtonian fluid and a closed lipid bilayer $\Gamma_t$ depending on time $t$. We have the usual Navier-Stokes system in the bulk $\Omega\setminus \Gamma_t$:
\begin{equation*}
\begin{aligned}
\rho_b\frac{Du}{Dt}&=\Div S,\\
\Div u&=0.
\end{aligned}
\end{equation*}
Here, $u$ is the fluid velocity, $\frac{D\,\cdot}{Dt}=\partial_t \cdot + (u\cdot\nabla)\,\cdot$ is the material derivative, in particular, $\frac{Du}{Dt}$ is the fluid particle acceleration, $S=2\mu_b  Du - \pi I$ is the (bulk) stress tensor, $Du$ is the symmetric part of the gradient of $u$, $\pi$ is the pressure, $\rho_b$ is the constant bulk fluid density, and $\mu_b$ is the constant dynamic viscosity of the bulk fluid. Obviously, we have $\Div S=\mu_b\Delta u-\grad\pi$. We assume that $u$ vanishes on $\partial\Omega$. As is done in all dynamical models mentioned above, we prescribe the no-slip condition at the membrane, that is, the velocity of the bulk fluid at the membrane equals the membrane velocity. For the normal part of the velocity, this can be justified by assuming an instantaneous osmotic equilibrium of the bulk fluids on either side of the membrane. Mathematically, the no-slip condition implies that the fluid velocity is continuous across the membrane. On $\Gamma_t$ we decompose the membrane 
velocity $u$ 
into its tangential and normal part, $u=v+w\,\nu$, where $\nu$ is the outer unit normal. In the following, we have to consider \emph{hybrid tensor fields} on $\Gamma_t$, that is, tensor fields that involve both tangential and non-tangential vectors and covectors; see for example \cite{scriven,aris}. An instance of such a hybrid tensor field is the velocity field $u$ on $\Gamma_t$ that is non-tangential in general. Another instance is the \emph{surface stress tensor} $T$ which is a hybrid $(1,1)$-tensor, taking a tangential direction and returning a force density that is, in general, not a tangent vector. We adopt the standard coordinate notation for such tensors and write $T^i_\alpha$, the greek index taking the values $1,2$ and the latin index taking the values $1,2,3$. In this paper, the latin index will always refer to Cartesian coordinates in $\R^3$ while the greek index refers to arbitrary coordinates on $\Gamma_t$. Likewise, we write $u^i$ for the Cartesian components of (not neccessarily 
tangential)
 vector fields $u$ on $\Gamma_t$. There exists a canonical covariant differential calculus for hybrid tensors; see \cite{aris}. Here, however, we only need the divergence operators, which take the form
\begin{equation*}
\begin{aligned}
\DIV u &=g^{\alpha\beta}\langle\pa_\alpha u,\pa_\beta\rangle_{e},\\
(\DIV T)^i&=g^{\alpha\beta} T^i_{\alpha;\beta},
\end{aligned}
\end{equation*}
where the semicolon in the second line denotes the usual covariant differentiation of the covectors $(T^i_{\alpha})_{\alpha=1,2}$ (for fixed $i$). Basic continuum mechanical considerations for a general interfacial fluid, see \Appref{app:int} or \cite{scriven,aris}, lead to the following conservation laws for linear momentum and mass on $\Gamma_t$:
% \footnote{It is interesting to note that $D\rho/Dt$ is the Lie derivative of the scalar $\rho$ with respect to $u$, while $Du/Dt$ is not the Lie derivative of the vector $u$ with respect to $u$ (although this is true componentwise). This is closely related to the fact that balance of total linear momentum is not a tensorial postulate, cf. [Marsden/Hughes].}
\begin{equation}\label{eqn:membrane}
\begin{aligned}
\rho\frac{Du}{Dt}&=\DIV T + [\![S]\!]\nu,\\
\frac{D\rho}{Dt}&=-\rho\DIV u.
\end{aligned}
\end{equation}
% \begin{equation}\label{eqn:membrane}
% \begin{aligned}
% 0&=\Div_e T + \langle T,k\rangle_e\, \nu + [S]\nu - \text{grad}_{L_2} W,\\
% \Div_e v&=w\, H.
% \end{aligned}
% \end{equation}
Here, $\rho$ is the surface fluid density and $[\![S]\!]$ is the jump of the bulk stress tensor across the membrane (subtracting the \emph{outer limit} from the \emph{inner limit}). We assume $\rho$ to be constant, reducing \eqref{eqn:membrane}$_2$ to $\DIV u=0$. Furthermore, assuming the interfacial fluid to be isotropic and Newtonian, the (strictly tangential) fluid part $^fT^i_\alpha={^f\tilde T_\alpha^\beta}\pa_\beta^i$ of the surface stress tensor takes the form
\[^f\tilde T_\alpha^\beta=-q\,\delta_\alpha^\beta+\zeta\,\delta_\alpha^\beta g^{\gamma\delta} (\D u)_{\gamma\delta}+\mu\Big(\delta^\delta_\alpha g^{\beta\gamma}+\delta_\alpha^\gamma g^{\beta\delta}-\frac23\,\delta_\alpha^\beta g^{\gamma\delta}\Big)(\D u)_{\gamma\delta};\]
see \cite{scriven,aris}. Here, $q$ is the surface pressure, acting as a Lagrange multiplier with respect to the constraint $\DIV u=0$, $\zeta,\mu$ are the surface dilatational and shear viscosities, respectively, and $\D u$ is the \emph{surface rate-of-strain tensor}. Analogously to the bulk case, $\D u$ is the \emph{Lie derivative} of $g$ of with respect to $u$, $(\D u)_{\alpha\beta}=\frac12(v_{\alpha;\beta}+v_{\beta;\alpha})-w\,k_{\alpha\beta}$; see \Appref{app:int} or \cite{scriven,aris}. We have $g^{\gamma\delta}(\D u)_{\gamma\delta}=\Div_g v-w\,H=\DIV u=0$; for the second identity see below. 
% In lipid bilayers, the dilatational viscosity is very small.\marginpar{Ref?} Hence, we set $\zeta=0$. 
We conclude that the fluid part of the surface stress tensor takes the form
\[^f\tilde T_\alpha^\beta=-q\,\delta_\alpha^\beta+2\mu\,(\D u)_\alpha^\beta=-q\,\delta_\alpha^\beta+\mu\,g^{\beta\gamma}(v_{\alpha;\gamma}+v_{\gamma;\alpha}-2w\,k_{\alpha\gamma}).\]
The Canham-Helfrich energy \eqref{eqn:ch} induces stresses in the membrane which are given by
\[^{h}T^i_\alpha=\kappa\,\big((H-C_0)^2/2\,\pa_\alpha^i-(H-C_0)\,k_\alpha^\beta\pa_\beta^i-(H-C_0)_{,\alpha}\nu^i\big),\]
where the comma denotes the usual partial differentiation; see \Appref{app:int} or \cite{guven}. We have
\begin{equation}\label{eqn:stress}
\DIV {^hT}=-\kappa\big(\Delta_g H + H(H^2/2-2K) +C_0(2K-HC_0/2)\big)\nu=-\grad_{L_2}F\,\nu. 
% \DIV {^hT}=-\kappa\big(\Delta_e H + (H-C_0)(H^2/2-2K+HC_0/2)\big)\nu=-\grad_{L_2}F\,\nu. 
\end{equation}
Finally, our constitutive assumption concerning the interfacial fluid is
\[T={^fT}+{^hT}.\]
The full system of equations now reads
 \begin{equation}\label{eqn:nonrex}
 \begin{aligned}
 \rho_b\frac{Du}{Dt}&=\Div S&&\mbox{ in }\Omega\setminus\Gamma_t,\\
 \Div u&=0&&\mbox{ in }\Omega\setminus\Gamma_t,\\
 \rho\frac{Du}{Dt}&=\DIV T + [\![S]\!]\nu &&\mbox{ on }\Gamma_t,\\
 \DIV u&=0&&\mbox{ on }\Gamma_t,\\
 u&=0&&\mbox{ on }\pa\Omega.
 \end{aligned}
 \end{equation}
This is the non-parabolic model we shall study. It can be interpreted as a vector field on the nonlinear manifold consisting of embedded surfaces $\Gamma\subset\Omega$ of fixed area and enclosed volume and of velocities $u:\Omega\rightarrow\R^3$, in the sense that $(\Gamma,u)$ is mapped to $([u]_\Gamma\cdot\nu,\pa_t u)$, where $\pa_t u$ is computed from \eqref{eqn:nonrex}$_{1,3}$; analogously to the classical Navier-Stokes system one can show that the pressure functions $\pi$ and $q$ at some point in time can be computed from the knowledge of $\Gamma$ and $u$ at that same instant alone.  

Let us put the terms involving the surface divergences in a more explicit form. From the decomposition $u=v+w\,\nu$ we infer
\begin{equation*}%\label{eqn:divu}
\begin{aligned}
\DIV u&=g^{\alpha\beta}\langle\pa_\alpha u,\pa_\beta\rangle_{e}=\pa_\alpha v^\alpha + v^\gamma g^{\alpha\beta} \langle\pa_\alpha\pa_\gamma,\pa_\beta\rangle_{e}+ w\, g^{\alpha\beta}\langle\pa_\alpha\nu,\pa_\beta\rangle_{e}\\
&= \pa_\alpha v^\alpha + v^\beta\,^g\Gamma_{\alpha\beta}^\alpha - w\, k_\alpha^\alpha
=\Div_g v-w\,H.
 \end{aligned}
\end{equation*}
Furthermore, decomposing the surface stress tensor into a tangential and a normal part, $T^i_\alpha={^tT}_\alpha^\beta\pa_\beta^i+{^nT}_\alpha\nu^i$, we obtain
\begin{equation*}
 (\DIV T)^i={^tT}^{\alpha\beta}_{\ \ \, ;\alpha}\pa_\beta^i+{^tT}^{\alpha\beta}\pa_{\beta;\alpha}^i+{^nT}^\alpha_{\ ;\alpha}\nu^i+{^nT}^\alpha\nu^i_{\ ,\alpha}.
\end{equation*}
Note that $\pa_{\alpha}^i$ is the $\alpha$-th component of the covector $X\mapsto \langle X,e^i\rangle_{e}$, where $e^i$ is the $i$-th standard basis vector in $\R^3$. Thus, $\pa_{\beta;\alpha}^i=\pa_{\beta,\alpha}^i-{^g\Gamma}_{\beta\alpha}^\gamma\pa_\gamma^i=k_{\beta\alpha}\nu^i$. Furthermore, by definition of the Weingarten map, we have $\nu^i_{\ ,\alpha}=-k_\alpha^\beta\pa_\beta^i$. Hence, we obtain the general formula
\begin{equation}\label{eqn:generaldiv}
 (\DIV T)^i=(\Div_g{^tT})^{\alpha}\pa_\alpha^i+{^tT}^{\alpha\beta}k_{\alpha\beta}\nu^i+(\Div_g{^nT})\, \nu^i-{^nT}^\alpha k_\alpha^\beta\pa_\beta^i.
\end{equation}
From this, let us compute $\DIV{^fT}$. We have
% Here, $\Div_e$ is the tangential divergence, i.e., $(\Div_e T)^\alpha=T^{\alpha\,\beta}_{\beta;}$, $T=2\mu_0\,d - p\, I$ is the tangential stress tensor, $\mu_0\ge 0$ is the constant surface viscosity, the viscous (rate of) strain tensor $d$ is given by the Lie derivative of $e$ of with respect to $u$, $p$ is the surface pressure, $k$ is the Weingarten map, $H$ is twice the scalar mean curvature, i.e., $H=k_\alpha^\alpha$, $[S]$ is the jump of $S$ across the interface, and, finally, $\text{grad}_{L_2} W$ is the $L_2$-gradient of the Willmore energy of $\Gamma_t$. The viscous strain tensor can be written in the form $d=D^e v-w\, k$, where $D^e v$ is the symmetric part of the covariant derivative of $v$, i.e., $(2D^e v)_\beta^\alpha=e^{\alpha\gamma}(v_{\gamma;\beta}+v_{\beta;\gamma})$. 
\[(\Div_g {^f\tilde T})^\beta=-g^{\alpha\beta}q_{,\alpha}+\mu\,g^{\alpha\delta}g^{\beta\gamma}\big(v_{\alpha;\gamma\delta} + v_{\gamma;\alpha\delta}-2(w\,k_{\alpha\gamma})_{;\delta}\big).\]
By definition of the Riemannian curvature tensor $R$, we have $v_{\alpha;\gamma\delta}=v_{\alpha;\delta\gamma}+R_{\delta\gamma\ \alpha}^{\ \ \, \beta}v_\beta$. Recalling that for a surface $R_{\delta\gamma\ \alpha}^{\ \ \, \beta}=K\,(\delta_\gamma^\beta g_{\alpha\delta}-\delta_\delta^\beta g_{\alpha\gamma})$ we obtain $g^{\alpha\delta}v_{\alpha;\gamma\delta}=(\Div_g v)_{,\gamma}+K v_\gamma$. Since $\Div_g v=w\,H$, this gives
\[\Div_g{^f\tilde T}=-\grad_g q+ \mu\,\big(\Delta_g v + \grad_g(w\,H) + K v - 2\Div_g(w\,k)\big).\]
Furthermore, we have
\begin{equation*}
\begin{aligned}
{^f\tilde T^{\alpha\beta}}k_{\alpha\beta}&=-q\,H+2\mu\,(\D u)^{\alpha\beta}k_{\alpha\beta}\\
 &=-q\,H+2\mu\,\big(v_{\alpha;\beta}k^{\alpha\beta}-w\,(H^2-2K)\big). 
\end{aligned}
\end{equation*}
In view of \eqref{eqn:generaldiv} and \eqref{eqn:stress} we conclude that
\begin{equation}\label{fullt}
\begin{aligned}
 \DIV T = &-\grad_g q -q\,H\nu + \mu\,\big(\Delta_g v + \grad_g(w\,H) + K v -2\Div_g(w\,k)\big)\\ &  +2\mu\big(\langle\nabla^g v,k\rangle_g-w\,(H^2-2K)\big)\nu\\
& -\kappa\big(\Delta_g H + H(H^2/2-2K)+C_0(2K-HC_0/2)\big)\nu.
\end{aligned}
\end{equation}
Note that the highest order term in the normal part of $\DIV T$ is a quasilinear fourth order differential operator acting on the membrane configuration while the highest order term in the tangential part is a second order operator acting linearly on the tangential membrane velocity. Recalling that $Du/Dt$ is the fluid particle \emph{acceleration}, we see that \eqref{eqn:nonrex}$_3$ is a fourth order wave-type equation governing membrane deformation which is strongly coupled to a two-dimensional Navier-Stokes type system governing the lateral flow on the membrane.\footnote{More precisely, the equation governing membrane deformation is not hyperbolic, but it is \emph{dispersive}. The principal part of its linearization has the form $\rho\pa_t^2+\kappa\Delta^2$. This operator essentially describes the linear evolution of a thin elastic plate. Note that it factorizes into two Schr\"odinger operators; in particular we can expect an infinite speed of propagation of disturbances.}

Let us have a look at the order of magnitude of the bulk Reynolds number $\mathcal{R}_b$ and the surface Reynolds number $\mathcal{R}$ of the system  \eqref{eqn:nonrex} in typical experiments; cf. \cite{seifert97}. Typical relevant length scales are of the order $L_{\text{typ}}=1\mu m$, while typical time scales (accessible to video-microscopy) are of the order $T_{\text{typ}}=10^{-3}s$. Furthermore, typical values for the bulk density and viscosity are $\rho_b=10^3 kg/m^3$ and $\mu_b=10^{-3}kg/ms$, respectively. This gives
\[\mathcal{R}_b=\frac{\rho_b L_{\text{typ}}^2}{\mu_b T_{\text{typ}}}=10^{-3}.\]
Assuming that the membrane volume density is comparable to that of water, for its surface density we obtain $\rho=10^{-5}kg/m^2$. Typically, the surface viscosity is of the order $\mu=10^{-9}kg/s$. Thus
\[\mathcal{R}=\frac{\rho\, L_{\text{typ}}^2}{\mu\, T_{\text{typ}}}=10^{-5}.\]
Hence, under such circumstances we can expect that the viscous forces strongly dominate the inertial forces. Neglecting the inertial terms in \eqref{eqn:nonrex}, we obtain a set of equations describing purely relaxational fluid vesicle dynamics:
\begin{equation}\label{eqn:final}
\begin{aligned}
\Div S&=0 &&\mbox{ in }\Omega\setminus \Gamma_t,\\
\Div u&=0&&\mbox{ in }\Omega\setminus \Gamma_t,\\
\DIV T &=-[\![S]\!]\nu&&\mbox{ on }\Gamma_t,\\
\DIV u&=0&&\mbox{ on }\Gamma_t,\\
u&=0&&\mbox{ on }\pa\Omega.
\end{aligned}
\end{equation}
Here, phase space consists of the embedded surfaces $\Gamma\subset\Omega$ of fixed area and enclosed volume alone. Note that, by \eqref{eqn:stress}, equation \eqref{eqn:final}$_3$ can be written in the form
\begin{equation}\label{eqn:spong}
 \DIV {^fT}+[\![S]\!]\nu=\grad_{L_2}F\,\nu.
\end{equation}
Thus, at a fixed instant in time we have to solve for the fluid velocity $u$ a linear Stokes-type system with prescribed Neumann-type boundary values given by the $L_2$-gradient of the Canham-Helfrich energy. Mapping $\Gamma$ to $[u]_\Gamma\cdot\nu$ then defines the dynamics of our system. Compared to the classical Canham-Helfrich flow (or Willmore flow), there is an additional Neumann-to-Dirichlet operator involved here. Since the $L_2$-gradient of $F$ is a fourth order operator, the mapping $\Gamma\mapsto [u]_\Gamma\cdot\nu$ can be considered as a nonlinear, nonlocal pseudo-differential operator of third order.\footnote{The situation is analogous to the relation between volume-conserving mean curvature flow and the Mullins-Sekerka system. There, however, we have a Dirichlet-to-Neumann operator involved. The mean curvature operator being of second order, the Mullins-Sekerka gradient is a third order operator, as in the present case.} The Stokes-type operator defined by 
the left hand sides of \eqref{eqn:final}$_{1,2,4,5}$ and \eqref{eqn:spong} will also play a crucial role in the analysis of \eqref{eqn:nonrex}.

\subsection{Some basic observations}
Let $\Gamma^i_t$, $i=1,\ldots,k$, denote the connected components of $\Gamma_t$, and denote by $\Omega_t^i$ the part of $\Omega$ enclosed by $\Gamma_t^i$. Reynolds' transport theorem then gives
\begin{equation*}%\label{eqn:volume}
\frac{d}{dt}|\Omega^{i}_t|=\frac{d}{dt}\int_{\Omega^{i}_t}dx=\int_{\Gamma_t}w\, dA=\int_{\Omega^{i}_t}\Div u\, dx=0.
\end{equation*}
In particular, the volume enclosed by each connected component of the membrane is conserved. Similarly, we can show that the area of each connected component is conserved:
\begin{equation*}%\label{eqn:area}
\frac{d}{dt}|\Gamma^{i}_t|=\frac{d}{dt}\int_{\Gamma_t^i}dA=-\int_{\Gamma_t^i}w\,H¸\, dA=-\int_{\Gamma_t}\Div_g v\, dA=0.
\end{equation*}

Now, let us present the gradient flow structure underlying the parabolic flow \eqref{eqn:final}. Consider the manifold $N$ of embedded surfaces of fixed area and enclosed volume. For $\Gamma\in N$, the tangent space $T_\Gamma N$ can be identified with the space of scalar fields $w$ on $\Gamma$ such that the linearized constraints
\begin{equation}\label{eqn:tmp}
\int_{\Gamma^i} w\, dA=0\quad\text{ and }\quad\int_{\Gamma^i} w\,H\, dA=0
\end{equation}
hold for all connected components $\Gamma^i$ of $\Gamma$. For $w\in T_\Gamma N$, consider the system
\begin{equation}\label{eqn:esystem}
\begin{aligned}
\Div S&=0 &&\mbox{ in }\Omega\setminus \Gamma,\\
\Div u&=0&&\mbox{ in }\Omega\setminus \Gamma,\\
P_\Gamma(\DIV {^fT} + [\![S]\!]\nu) &= 0&&\mbox{ on }\Gamma,\\
\DIV u&=0&&\mbox{ on }\Gamma,\\
u\cdot\nu&=w&&\mbox{ on }\Gamma,\\
u&=0&&\mbox{ on }\pa\Omega.
\end{aligned}
\end{equation}
Note that the conditions \eqref{eqn:tmp} are necessary for the solvability of these equations, due to the incompressibility constraints. For $w_1,w_2\in T_\Gamma N$, define the Riemannian metric on $N$ associated with fluid vesicle dynamics by
\begin{equation}\label{eqn:rmetric}
 \langle w_1,w_2\rangle_V:=2\mu_b\int_{\Omega\setminus\Gamma}\langle Du_1,Du_2\rangle_e\,dx+2\mu\int_{\Gamma}\langle \D u_1,\D u_2\rangle_g\,dA,
\end{equation}
where $u_1,u_2$ solve the system \eqref{eqn:esystem} with data $w_1,w_2$. Note that the length of a curve in $N$ endowed with this metric is given by the energy dissipated during the corresponding forced deformation of the membrane. The representation of $-dF$ with respect to the metric \eqref{eqn:rmetric} is given by $[u]_\Gamma\cdot\nu$, where $u$ solves \eqref{eqn:final}.
% \begin{equation}\label{eqn:esystem2}
% \begin{aligned}
% \Div S&=0 &&\mbox{ in }\Omega\setminus M,\\
% \Div u&=0&&\mbox{ in }\Omega\setminus M,\\
% \DIV {^fT} - [\![S]\!]\nu &= \grad_{L_2}F\,\nu&&\mbox{ on }M,\\
% \DIV u&=0&&\mbox{ on }M,\\
% u&=0&&\mbox{ on }\pa\Omega.
% \end{aligned}
% \end{equation}
Indeed, for all $w\in T_\Gamma N$ and corresponding solutions $\tilde u$ of \eqref{eqn:esystem} we have
\begin{equation}\label{gradient}
\begin{aligned}
\langle u\cdot\nu,w\rangle_V&=2\mu_b\int_{\Omega\setminus\Gamma}\langle Du,D\tilde u\rangle_e\,dx+2\mu\int_{\Gamma}\langle \D u,\D \tilde u\rangle_g\,dA\\
&=-\int_{\Gamma}\langle[\![S]\!]\nu,\tilde u\rangle_{e}\,dA-\int_{\Gamma}\langle\DIV {^fT},\tilde u\rangle_{e}\,dA=\int_{\Gamma}\langle\DIV {^hT},\tilde u\rangle_{e}\,dA\\
&=-\int_\Gamma \grad_{L_2}F\,w\,dA=-dF(w).
\end{aligned}
\end{equation}
Here, $S$ and $^fT$ denote the stress tensors with respect to $u$, and we used integration by parts for the second identity (see below), \eqref{eqn:final}$_3$ for the third identity, and \eqref{eqn:stress} for fourth identity.
We conclude that, indeed, \eqref{eqn:final} is the gradient flow of the Canham-Helfrich energy on $N$ endowed with the Riemannian metric \eqref{eqn:rmetric}. In particular, the energy $F$ is a strict Lyapunov functional, and, along the flow,
\begin{equation*}
\frac{d}{dt} F=dF(u\cdot\nu)=-\langle u\cdot\nu,u\cdot\nu\rangle_V=-2\mu_b\int_{\Omega\setminus \Gamma_t}|Du|^2_{e}\,dx-2\mu\int_{\Gamma_t}|\D u|^2_g\,dA.
\end{equation*}
It remains to prove the second identity in \eqref{gradient}. In view of \eqref{eqn:final}$_{1,2}$ we have
\begin{equation*}
\begin{aligned}
0=\int_{\Omega\setminus \Gamma}\langle\Div S,\tilde u\rangle_{e}\,dx&= -\int_{\Gamma}\langle[\![S]\!]\nu,\tilde u\rangle_{e}\,dA-\int_{\Omega\setminus \Gamma}\langle S,\nabla \tilde u\rangle_{e}\,dx\\
&=-\int_{\Gamma}\langle[\![S]\!]\nu,\tilde u\rangle_{e}\,dA-2\mu_b\int_{\Omega\setminus \Gamma}\langle Du,D\tilde u\rangle_{e}\,dx
\end{aligned}
\end{equation*}
Furthermore, we compute
\begin{equation*}
\sum_i g^{\alpha\beta}\,{^fT}^i_\alpha\,\tilde u^i,_\beta=-q\,g^{\alpha\beta}\langle\pa_\alpha,\pa_\beta \tilde u\rangle_{e} + 2\mu\,g^{\alpha\beta}(\D u)_{\alpha}^\gamma\,\langle\pa_\gamma,\pa_\beta \tilde u\rangle_{e}.
% \int_{\Gamma_t}\langle\DIV {^fT},u\rangle_{\tilde e}\,dA=
% \int_{\Gamma_t}\langle\Div_e {^fT},v\rangle_{e}+w\,\langle{^fT},k\rangle_{e}\,dA\\
% =-\int_{\Gamma_t}\langle{^fT},\nabla^e v\rangle_{e}-w\,\langle{^fT},k\rangle_{e}\,dA\\=-\int_{\Gamma_t}-q\Div_e v +2\eta\langle D,k,\rangle_e+q\,w\, H-2\eta\langle D,w\,k\rangle_e\,dA=-2\eta\int_{\Gamma_t}|D|^2_e\,dA.
\end{equation*}
Note that $g^{\alpha\beta}\langle\pa_\alpha,\pa_\beta\tilde u\rangle_{e}=\DIV\tilde u=0$ and $\langle\pa_\gamma,\pa_\beta\tilde u\rangle_{e}=\tilde v_{\beta;\alpha}-\tilde w\,k_{\alpha\beta}$. Using integration by parts and the fact that $\D u$ is symmetric, we obtain
\begin{equation*}
\begin{aligned}
 \int_{\Gamma}\langle\DIV {^fT},\tilde u\rangle_{e}\,dA=-2\mu\int_{\Gamma} \langle\D u,\D\tilde u\rangle_g\, dA.
\end{aligned}
\end{equation*}

Furthermore, from the computations above and Reynolds' transport theorem we see that the energy identity for \eqref{eqn:nonrex} reads
\begin{equation*}
\frac{d}{dt} \bigg(\frac{\rho_b}{2}\int_{\Omega\setminus \Gamma_t}|u|_e^2\,dx+\frac{\rho}{2}\int_{\Gamma_t}|u|_e^2\,dA+F\bigg)=-2\mu_b\int_{\Omega\setminus \Gamma_t}|Du|^2_{e}\,dx-2\mu\int_{\Gamma_t}|\D u|^2_g\,dA.
\end{equation*}

Finally, let us analyze the equilibria of \eqref{eqn:nonrex} and \eqref{eqn:final} a little more closely. From the energy identities we see that $Du=0$ in $\Omega\setminus\Gamma$; taking into account the no-slip conditions and Korn's equality, this gives $u=0$ in $\Omega$. Now, \eqref{eqn:nonrex}$_1$ and \eqref{eqn:final}$_1$ reduce to $\grad \pi=0$ in $\Omega\setminus\Gamma$, that is, $\pi$ is constant in each connected bulk component. Furthermore, \eqref{eqn:nonrex}$_3$ and \eqref{eqn:final}$_3$ reduce to
\[\grad_g q + q\,H\nu+\grad_{L_2}F\,\nu=-[\![\pi]\!]\nu.\]
Since the surface pressure gradient is the only tangential component in this equation, $q$ is constant on each connected component of $\Gamma$. Thus, we have
\begin{equation*}
 \grad_{L_2}F+[\![\pi]\!]+q\,H=0.
\end{equation*}
This is the Helfrich equation \eqref{eqn:helfrich} with the pressure jump and the surface pressure acting as Lagrange multipliers with respect to the volume and area constraints.

\section{$L_2$-Theory for Stokes-type systems}\seclabel{stokes}

Let $\Omega\subset\R^3$ be a smooth, bounded domain and $\Gamma\subset\Omega$ a smoothly embedded closed surface. We write $\Gamma^i$, $i=1,\ldots,k$, for the connected components of $\Gamma$, $\Omega^i$ for the open set enclosed by $\Gamma^i$, and we let \[\Omega^0:=\Omega\setminus\big(\bigcup_{i=1}^k \Gamma^i\cup\Omega^i\big).\] 
Let us assume that $\mu_b,\mu>0$. We shall analyze two slightly different Stokes-type systems. The first one is related to the Neumann-to-Dirichlet operator discussed above and reads
\begin{equation}\label{s1}
\begin{aligned}
\Div S&=f_1 &&\mbox{ in }\Omega\setminus \Gamma,\\
\Div u&=f_2&&\mbox{ in }\Omega\setminus \Gamma,\\
\DIV {^fT} + [\![S]\!]\nu&=f_3&&\mbox{ on }\Gamma,\\
\DIV u&=f_4&&\mbox{ on }\Gamma,\\
u&=0&&\mbox{ on }\pa\Omega
\end{aligned}
\end{equation}
for appropriate data $(f_1,\ldots,f_4)$. From \eqref{s1}$_2$ we obtain
\[\int_{\Gamma^i} w\,dA=\int_{\Omega^i}f_2\,dx.\]
Thus, recalling $\DIV u=\Div_g v-w\,H$, by \eqref{s1}$_4$ we have
\begin{equation}\label{comp1}
 \int_{\Gamma^i}f_4/H\,dA=-\int_{\Omega^i}f_2\,dx \quad\text{for each $\Gamma^i$ that is a CMC surface}.
\end{equation}
Recall that the only closed CMC (= constant mean curvature) surfaces embedded in $\R^3$ are (collections of) round spheres. Furthermore, of course, we have
\begin{equation}\label{comp2}
 \int_{\Omega}f_2\,dx=0.
\end{equation}
Hence, for $s\ge 1$, we define the space of data
\begin{equation*}
\begin{aligned}
\F^s(\Gamma):=\big\{(f_1,\ldots,f_4)\,|\,f_1\in H^{s-1}(\Omega\setminus\Gamma;\R^3),\ f_2\in H^s(\Omega\setminus\Gamma),\\ P_\Gamma f_3\in H^{s-1}(\Gamma;T\Gamma),\ f_3\cdot\nu\in H^{s-1/2}(\Gamma),\ f_4\in H^s(\Gamma)\\ \text{such that \eqref{comp1} and \eqref{comp2} hold}\big\}
\end{aligned}
\end{equation*}
endowed with the canonical norm. Furthermore, for $s\ge 0$, we introduce the space of solutions
\begin{equation*}
\begin{aligned}
\E^s(\Gamma):=\big\{(u,\pi,q)\,|\,u\in H^{s+1}(\Omega\setminus\Gamma;\R^3)\cap H^1_0(\Omega;\R^3),\ P_\Gamma [u]_\Gamma\in H^{s+1}(\Gamma;T\Gamma),\\
\pi\in H^{s}(\Omega\setminus\Gamma),\ q\in H^{s}(\Gamma)\ \ \text{such that \eqref{comp1}, \eqref{comp2} hold with}\\ \text{$f_2=\pi$ and $f_4=q$}\big\}
\end{aligned}
\end{equation*}
endowed with the canonical norm. The conditions on $\pi$ and $q$ provide a gauge fixing; as is typical for Stokes-type equations the pressure functions in \eqref{s1} are not uniquely determined. The second Stokes-type system we shall analyze is related to the Riemannian metric discussed above and reads
\begin{equation}\label{s2}
\begin{aligned}
\Div S&=f_1 &&\mbox{ in }\Omega\setminus \Gamma,\\
\Div u&=f_2&&\mbox{ in }\Omega\setminus \Gamma,\\
P_\Gamma(\DIV {^fT} + [\![S]\!]\nu)&=f_3&&\mbox{ on }\Gamma,\\
\DIV u&=f_4&&\mbox{ on }\Gamma,\\
u\cdot\nu&=f_5&&\mbox{ on }\Gamma,\\
u&=0&&\mbox{ on }\pa\Omega
\end{aligned}
\end{equation}
for appropriate data. Of course, again, \eqref{comp1} and \eqref{comp2} must hold. Moreover, we have
\begin{equation}\label{comp3}
\int_{\Omega^i}f_2\, dx=\int_{\Gamma^i}f_5\,dA\quad\text{and}\quad \int_{\Gamma^i}f_4+f_5\,H\,dA=0\quad\text{for all }i=1,\ldots,k.
\end{equation}
Note that \eqref{comp1} is contained in \eqref{comp3}. Hence, for $s\ge 1$, we define the space of data
\begin{equation*}
\begin{aligned}
\F^s_\nu(\Gamma):=\big\{(f_1,\ldots,f_5)\,|\,f_1\in H^{s-1}(\Omega\setminus\Gamma;\R^3),\ f_2\in H^s(\Omega\setminus\Gamma),\ f_3\in H^{s-1}(\Gamma;T\Gamma),\\
f_4\in H^s(\Gamma),\ f_5\in H^{s+1/2}(\Gamma)\ \ \text{such that \eqref{comp2}, \eqref{comp3} hold}\big\}
\end{aligned}
\end{equation*}
endowed with the canonical norm. Furthermore, for $s\ge 0$, we introduce the space of solutions
\begin{equation*}
\begin{aligned}
\E^s_\nu(\Gamma):=\big\{(u,\pi,q)\,|\,u\in H^{s+1}(\Omega\setminus\Gamma;\R^3)\cap H^1_0(\Omega;\R^3),\ P_\Gamma [u]_\Gamma\in H^{s+1}(\Gamma;T\Gamma),\\
\pi\in H^{s}(\Omega\setminus\Gamma),\ q\in H^{s}(\Gamma)\ \ \text{such that \eqref{comp3} holds with}\\ \text{$f_2=\pi$, $f_4=q$, and $f_5=0$}\big\}
\end{aligned}
\end{equation*}
endowed with the canonical norm.

For the applications we have in mind it is important to include the case of $s$ being half-integer valued. For instance, both the local well-posedness of \eqref{eqn:nonrex} and the analysis of the formation of singularities in \eqref{eqn:final} rely on energy estimates on the surface level, that is, yielding an integer class regularity for the normal velocity $[u]_\Gamma\cdot\nu$. But then we must apply the theory developed in the present article with $s=n+1/2$ for some $n\in\N$. Since it doesn't cause any additional effort we prove our results for arbitrary real $s\ge 1$.

Throughout the article (except for a few formulae in the proof of Theorem \ref{s1h} and in \Appref{app:spaces}), we keep track of the physical dimensions in all mathematical expressions. This mainly serves as a rough test for the validity of the estimates we prove. We write $c_{l}$, $l\in\mathbb{Z}$, for generic positive constants with the physical dimension of length$^l$ depending only on $\Omega$, $\Gamma$, the \emph{Saffman-Delbr\"uck} length $\mu/\mu_b$, and $s$; let $c:=c_0$.

\subsection{Analysis for a fixed membrane}

Let us start with the analysis of \eqref{s1}. To begin with, let us assume that $f_2$ and $f_4$ vanish, and define
\begin{equation*}
\begin{aligned}
X:=\big\{u\in H^1_0(\Omega;\R^3)\,|\,\Div u=0 \text{ in }\Omega\setminus \Gamma,\ \DIV u=0 \text{ on } \Gamma,\ P_\Gamma[u]_\Gamma\in H^1(\Gamma;T\Gamma)\big\}.
\end{aligned}
\end{equation*}
Multiplying \eqref{s1}$_1$ by $\phi\in X$, integrating over $\Omega\setminus \Gamma$, and  integrating by parts we obtain
\begin{equation}\label{weak}
\begin{aligned}
0&=\int_{\Omega\setminus \Gamma} \langle\Div S,\phi\rangle_e\,dx-\int_{\Omega\setminus \Gamma} \langle f_1,\phi\rangle_e\,dx\\
&=-\int_{\Omega\setminus \Gamma} \langle S,\nabla\phi\rangle_e\,dx-\int_{\Gamma} \langle [\![S]\!]\nu,\phi\rangle_e\,dA-\int_{\Omega\setminus \Gamma} \langle f_1,\phi\rangle_e\,dx\\
&=-\int_{\Omega\setminus \Gamma} \langle S,D\phi\rangle_e\,dx+\int_{\Gamma} \langle \DIV {^fT},\phi\rangle_e\,dA-\int_\Gamma \langle f_3,\phi\rangle_e\,dA-\int_{\Omega\setminus \Gamma} \langle f_1,\phi\rangle_e\,dx\\
&=-2\mu_b\int_{\Omega\setminus \Gamma} \langle Du,D\phi\rangle_e\,dx-2\mu\int_{\Gamma} \langle \D u,\D \phi\rangle_g\,dA\\
&\quad-\int_\Gamma \langle f_3,\phi\rangle_e\,dA-\int_{\Omega\setminus \Gamma} \langle f_1,\phi\rangle_e\,dx\\
&=:-B(u,\phi)+F(\phi).
\end{aligned}
\end{equation}
Let us have a closer look at the bilinear form $B$.

\begin{lemma}\label{lemma:coercive}
$B$ is a symmetric, continuous, and coercive bilinear form on $X$. For all $u\in X$, we have
\begin{equation}\label{coercive}
\begin{aligned}
\lambda \|u\|^2_{2}\le\|\nabla u\|^2_{2}&\le \frac{1}{\mu_b}B(u,u),\\
\|\nabla^g v\|^2_{2}&\le \frac{c}{\mu}B(u,u),
\end{aligned}
\end{equation}
where $\lambda$ is the smallest eigenvalue of the Dirichlet-Laplacian in $\Omega$.
\end{lemma}
\proof Symmetry and continuity are obvious. Let us prove coercivity. Recall that $2|Du|^2_e=|\nabla u|^2_e + \Div((u\cdot\nabla)u)$. Thus, we have
\begin{equation}\label{du}
 2\int_{\Omega\setminus\Gamma}|Du|^2_e\,dx=2\int_{\Omega}|D u|^2_e\,dx=\int_{\Omega}|\nabla u|^2_e\,dx.
\end{equation}
%\begin{equation}\label{du}
% 2\int_{\Omega\setminus\Gamma}|Du|^2_e\,dx=\int_{\Omega\setminus\Gamma}|\nabla u|^2_e\,dx + \int_\Gamma [\![(u\cdot\nabla)u]\!]\cdot\nu\,dA=\int_{\Omega\setminus\Gamma}|\nabla u|^2_e\,dx.
%\end{equation}
By Poincar\'e's inequality, \eqref{coercive}$_1$ follows. Moreover, we have 
\[|\D u|_g^2=|\Def(v)|_g^2+w^2|k|^2_g-2\langle\Def(v),k\rangle_g\,w,\] where $\Def(v)_{\alpha\beta}=\frac12(v_{\alpha;\beta}+v_{\beta;\alpha})$ is the deformation tensor of $v$, and, similarly to the bulk case, $2|\Def(v)|_g^2=|\nabla v|_g^2+\Div_g(\nabla^g_v v) - v_{\beta;\alpha}^{\quad \, \beta}v^\alpha$. However, by definition of (intrinsic) curvature, we have $v_{\beta;\alpha}^{\quad \, \beta}v^\alpha=(\Div v)_{,\alpha} v^\alpha+ K|v|^2.$ Thus, taking into account $\Div_gv=w\,H$, we obtain
\begin{equation*}
\begin{aligned}
2\int_{\Gamma} |\D u|_g^2\, dA=\int_{\Gamma}\big( |\nabla v|_g^2 +w^2(H^2+2|k|_g^2)-K|v|^2_g-4\langle\Def(v),k\rangle_g\,w\big)\,dA.
\end{aligned}
\end{equation*}
This gives
\begin{equation*}
\begin{aligned}
\int_{\Gamma} |\nabla v|_g^2\,dA\le \frac1\mu B(u,u) + \,\|k\|_{\infty}^2\int_{\Gamma}|v|^2_g\,dA +  4\,\|k\|_{\infty}\int_{\Gamma}|\nabla v|_g|w|\,dA.
\end{aligned}
\end{equation*}
%\begin{equation*}
%\begin{aligned}
%\int_{\Gamma} |\nabla v|_g^2\,dA\le \frac1\mu B(u,u) + 3\,\|k\|_{L^\infty(\Gamma)}^2\int_{\Gamma}|v|^2_g+w^2\,dA\\ + \, 4\,\|k\|_{L^\infty(\Gamma)}\int_{\Gamma}|\nabla v|_g|w|\,dA.
%\end{aligned}
%\end{equation*}
Applying H\"older's and Young's inequality, the continuity of the trace operator, and Poincar\'e's inequality, we obtain
\begin{equation}\label{c2}
\begin{aligned}
\int_{\Gamma} |\nabla v|_g^2\,dA&\le \frac2\mu B(u,u) + 16\,\|k\|_{\infty}^2\int_{\Gamma}|v|^2_g+w^2\,dA\\
&\le \frac2\mu B(u,u) + c_{-1}\int_{\Omega}|\nabla u|^2_e\,dx.
\end{aligned}
\end{equation}
Combining \eqref{coercive}$_1$ and \eqref{c2} we obtain \eqref{coercive}$_2$.
\qed

\medskip
By Riesz' representation theorem in Hilbert space, for each $F\in X'$ there exists a unique $u\in X$ such that $B(u,\phi)=F(\phi)$ for all $\phi\in X$.
% and\marginpar{Evt. nochmal pruefen!}
% \begin{equation}\label{B}
%  B(u,u)=\|F\|_{X'}\|u\|_X.
% \end{equation}
Now, we reconstruct the pressure functions. Consider the space
\[Y:=\big\{u\in H^1_0(\Omega;\R^3)\,|\, P_\Gamma [u]_\Gamma\in H^1(\Gamma;T\Gamma)\big\}\]
endowed with the norm
\begin{equation*}
\|u\|_Y^2:=\|\nabla u\|_2^2+|\Omega|^\frac13\|\nabla^g P_\Gamma[u]_\Gamma\|_2^2
\end{equation*}
and the space
\[Z:=\big\{(f_2,f_4)\in L_2(\Omega)\times L_2(\Gamma)\ \ \text{such that \eqref{comp1}, \eqref{comp2} hold}\big\}\]
endowed with the norm
\begin{equation*}
\|(f_2,f_4)\|_Z^2:=\|f_2\|_2^2+|\Omega|^\frac13\|f_4\|_2^2.
\end{equation*}
Note that $X$ is a closed subspace of $Y$. Each $(\pi,q)\in Z$ defines an element of $Z'$ via
\[(f_2,f_4)\mapsto \int_{\Omega\setminus\Gamma}\pi\,f_2\,dx+\int_\Gamma q\,f_4\,dA.\]
In fact, by Riesz' theorem, this defines an isometry
\begin{equation}\label{ziso}
(Z,\vertiii{\cdot}_Z) \stackrel{\sim}{\longrightarrow} (Z',\|\cdot\|_{Z'}),
\end{equation}
where
\[\vertiii{(\pi,q)}_Z^2:=\|\pi\|_2^2+|\Omega|^{-\frac13}\|q\|_2^2\]
defines an equivalent norm on $Z$.
% Furthermore, we need the space $L_2_p(\Gamma):=(L_2(\Omega)\times L_2(\Gamma))/\sim$ endowed with the quotient norm, where 
% \[\|(\pi,q)\|_{L_2(\Omega)\times L_2(\Gamma)}:=\frac{1}{\mu_b}\|\pi\|_2+\frac{1}{\mu}\|q\|_2.\]

\begin{lemma}\label{lemma:div}
 The operator $\underline{\Div}:Y/X\rightarrow Z$, defined by
 \[\phi\mapsto(\Div\phi,\DIV[\phi]_{\Gamma}),\]
 is an isomorphism. In particular, we have
 \begin{equation}\label{estdiv}
\begin{aligned}
\|\phi\|_{Y/X}&\le c\|\DIV\phi\|_Z.
\end{aligned}
\end{equation}
\end{lemma}
\proof Obviously, $\underline{\Div}$ is well-defined and injective, so let us prove surjectivity. To this end, let $(f_2,f_4)\in Z$. Choose a function $w$ on $\Gamma$ in the following way: On $\Gamma^i$ choose $w\in\span\{1,H\}$ such that
\begin{equation}\label{wf}
 \int_{\Gamma^i}wH\,dA=-\int_{\Gamma^i}f_4\,dA\quad\text{and}\quad\int_{\Gamma^i}w\,dA=\int_{\Omega^i}f_2\,dA,
\end{equation}
which is possible by definition of $Z$. Now, let $v\in H^1(\Gamma;T\Gamma)$ solve
\[\Div_g v=f_4+ wH\]
on $\Gamma^i$; note that the right hand side is mean value free and employ, for instance, $L_2$-theory for the Laplacian on $\Gamma^i$. Let $\phi:=v+w\,\nu\in H^1(\Gamma;\R^3)$ and extend this function to the bulk, giving $\phi_0\in H^1_0(\Omega;\R^3)$. Furthermore, by $L_2$-theory for the divergence operator in $\Omega\setminus\Gamma$ there exists a $\phi_1\in H^1_0(\Omega\setminus\Gamma)$ such that
\begin{equation}\label{divphi1}
\Div\phi_1=f_2-\Div\phi_0
\end{equation}
in $\Omega\setminus\Gamma$; note that the right hand side is mean value free in $\Omega^i$ for all $i=0,\ldots,k$. For $i\not=0$ this is obvious from \eqref{wf} while for $i=0$ we additionally take into account \eqref{comp2} and compute
\[\int_{\Omega^0}\Div\phi_0\,dx=-\int_{\Gamma}w\,dA=-\int_{\cup_i\Omega^i}f_2\,dA=\int_{\Omega^0}f_2\,dA.\]
Finally, for $\phi:=\phi_0+\phi_1\in Y$ we have $\underline{\Div}\,\phi=(f_2,f_4)$.
\qed

\begin{corollary}\label{cor:pressure}
 The operator $\underline{\nabla}:Z\rightarrow X^\perp\subset Y'$, defined by
 \begin{equation*}%\label{pressure}
  \underline{\nabla}(\pi,q)(\phi)=-\int_{\Omega\setminus\Gamma}\pi\,\Div\phi\,dx - \int_{\Gamma}q\,\DIV\phi\,dA
 \end{equation*}
 for $\phi\in Y$, is an isomorphism. We have
 \begin{equation}\label{estpressure}
\begin{aligned}
 \vertiii{(\pi,q)}_Z&\le c\|\underline\nabla (\pi,q)\|_{Y'}
\end{aligned}
\end{equation}
with the same dimensionless constant $c$ as in \eqref{estdiv}.
\end{corollary}
\proof This follow from Lemma \ref{lemma:div}, the fact that
\[\underline{\nabla}=-\underline{\Div}':Z'\rightarrow (Y/X)'\]
is an isomorphism, from the isometry $(Y/X)'\simeq X^\perp$, and by Riesz' isometry \eqref{ziso}. Recall also that adjoining preserves the operator norm.
\qed

\medskip
Note that Corollary \ref{cor:pressure} and Lemma \ref{lemma:div} are equivalent to the inequality
\[\inf_{(\pi,q)\in Z}\sup_{\phi\in Y}\frac{1}{\|(\pi,q)\|_Z\|\phi\|_Y}\Big(\int_{\Omega\setminus\Gamma}\pi\Div\phi\,dx+\int_\Gamma q\DIV\phi\,dA\Big)>0,\]
which is a \emph{Lady\v{z}enskaja-Babu\v{s}ka-Brezzi-type condition}; cf. \cite{garcke}.

Consider the subspace $U\subset L_2(\Omega)\times L_2(\Gamma)$ defined as follows: $(\pi,q)\in U$ iff for all $i=1,\ldots,k$ we have
\begin{itemize}
 \item[(i)] $\pi=\kappa_i$ in $\Omega^i$, $\pi=\kappa_0$ in $\Omega^0$, $q=\kappa^i$ on $\Gamma^i$ with $\kappa_i,\kappa_0,\kappa^i\in\R$
 \item[(ii)] If $\Gamma^i$ is a round sphere with $H$ denoting twice the mean curvature, then $\kappa_i-\kappa_0=\kappa^i H$.
 \item[(iii)] If $\Gamma^i$ is not a round sphere, then $\kappa^i=0$ and $\kappa_i=\kappa_0$.
\end{itemize}
It is not hard to prove that $(L_2(\Omega)\times L_2(\Gamma))/U\simeq Z'\simeq Z$. Hence, the subspace $U$ characterizes the gauge freedom of the pressure functions. Of course, this can also be seen directly. To this end, consider $\pi\in L_2(\Omega)$ and $q\in L_2(\Gamma)$ such that $\underline{\nabla}(\pi,q)(\phi)=0$ for all $\phi\in Y$. Then, obviously, $\pi$ must be constant in $\Omega^i$ for all $i=0,\ldots,k$. Integration by parts yields
 \begin{equation}\label{uibp}
 0=\int_{\Gamma}[\![\pi]\!]\,\phi\cdot\nu\,dx + \int_{\Gamma}q\,\DIV\phi\,dA.
 \end{equation}
For $\phi\in Y$ such that $\phi\cdot\nu=0$ on $\Gamma$ this reduces to
\[0=\int_{\Gamma}q\,\Div_g \phi\,dA,\]
proving that $q$ is constant on $\Gamma^i$ for all $i=1,\ldots,k$, and hence $[\![\pi]\!]=q\,H$. If $\Gamma^i$, $i=1,\ldots,k$, is a round sphere, then $q$ is determined on $\Gamma^i$ by the jump of the bulk pressure. On the other hand, if $\Gamma^i$ is not a round sphere, then $[\![\pi]\!]=0$ and $q=0$.

Multiplying \eqref{s1} with test functions $\phi\in Y$, integrating over $\Omega$ and $\Gamma$, respectively, and integrating by parts we obtain a weak formulation for this system. For functions $f_1:\Omega\rightarrow\R^3$ and $f_3:\Gamma\rightarrow\R^3$ we define
\[\|f_1\|_{-1,2}:=\sup_{\varphi\in H^1_0(\Omega)\atop \|\nabla\varphi\|_2\le 1}\int_\Omega f_1\cdot\varphi\,dx\quad\text{and}\quad\|f_3\|_{-\frac12,2}:=\sup_{\varphi\in H^1_0(\Omega)\atop \|\nabla\varphi\|_2\le 1}\int_\Gamma f_3\cdot\varphi\,dA.\]

\begin{theorem}\label{thm:s1}
 For all $(f_1,\ldots,f_4)\in\F^1(\Gamma)$ there exists a unique weak solution $(u,\pi,q)\in\E^0(\Gamma)$ of \eqref{s1}. We have+++
\begin{equation}\label{est1}
\begin{aligned}
&\mu_b^\frac12 \|\nabla u\|_{2} + \mu^\frac12 \|\nabla^g v\|_{2} + \frac{1}{\mu_b^\frac12} \|\pi\|_{2} + \frac{1}{\mu^\frac12} \|q\|_{2}\\
&\quad\le c\Big(\frac{1}{\mu_b^\frac12}\|f_1\|_{-1,2}+\mu_b^\frac12\|f_2\|_{2}+\frac{1}{\mu_b^\frac12}\|f_3\|_{-\frac12,2}+\mu^\frac12\|f_4\|_{2}\Big).
\end{aligned}
\end{equation}
%  \begin{equation}\label{est1}
%  \begin{aligned}
% &\mu_b^\frac12 \|\nabla u\|_{2} + \mu^\frac12 \|\nabla^g v\|_{2} + \frac{1}{\mu_b^\frac12} \|\pi\|_{2} + \frac{1}{\mu^\frac12} \|q\|_{2}\\
% &\quad\le \frac{c_1}{\mu_b^\frac12}\|f_1\|_{2}+c\,\mu_b^\frac12\|f_2\|_{2}+\frac{c_{\frac12}}{\mu_b^\frac12}\|f_3\|_{2}+c\,\mu^\frac12\|f_4\|_{2}.
% \end{aligned}
% \end{equation}
\end{theorem}
\proof By Lemma \ref{lemma:div}, there exists a $u_0\in Y$ such that $\underline\Div  \,u_0=(f_2,f_4)$ and
 \begin{equation}\label{estdivk}
\begin{aligned}
\|u_0\|_Y&\le c\|(f_2,f_4)\|_Z.
\end{aligned}
\end{equation}
% \begin{equation}\label{esttilde}
% \|\nabla\tilde u\|_{L_2(\Omega\setminus\Gamma)}^2+\frac{1}{c_1}\|\nabla\tilde v\|_{L_2(\Gamma)}^2+\frac{1}{c_3}\|\tilde w\|_{L_2(\Gamma)}^2 \le c\,\|f_2\|_{L_2(\Omega\setminus\Gamma)}^2 + c_1\,\|f_4\|_{L_2(\Gamma)}^2.
% \end{equation}
Hence, it remains to solve \eqref{s1} with $\tilde f_2=0$, $\tilde f_4=0$, $\tilde f_1=f_1-2\mu_b\Div Du_0$, and $\tilde f_3=f_3-2\mu\DIV\D u_0-[\![2\mu_b\,Du_0]\!]\nu$. To this end, we define $F_{u_0}\in X'$ by
\begin{equation}\label{fun}
\begin{aligned}
F_{u_0}(\phi)=\int_{\Omega\setminus\Gamma}\big(f_1\cdot\phi+2\mu_b\langle Du_0,D\phi\rangle_e\big)\,dx+\int_{\Gamma}\big(f_3\cdot\phi+2\mu\langle \D u_0,\D\phi\rangle_g\big)\,dA,
 \end{aligned}
\end{equation}
and we let $u_1\in X$ solve $B(u_1,\phi)=F_{u_0}(\phi)$ for all $\phi\in X$. We have
\begin{equation}\label{bf}
\begin{aligned}
B(u_1,u_1)&=F_{u_0}(u_1)\\
&\le \|\nabla u_1\|_{2}\|f_1\|_{-1,2} + 2\mu_b\|D u_1\|_{2}\|\nabla u_0\|_{2}+\|\nabla u_1\|_{2}\|f_3\|_{-\frac12,2}\\
&\quad+2\mu\|\D u_1\|_{2}\big(\|\nabla^g v_0\|_{2}+\|k\|_{\infty}\|w_0\|_{2}\big)\\
&\le B(u_1,u_1)^\frac12\Big(\frac{1}{\mu_b^{\frac12}}\|f_1\|_{-1,2}+(2\mu_b)^\frac12\|\nabla u_0\|_2+\frac{1}{\mu_b^\frac12}\|f_3\|_{-\frac12,2}\\
&\quad+(2\mu)^\frac12\|\nabla^g v_0\|_{2}+\mu^\frac12c_{-\frac12}\|\nabla u_0\|_2\Big),
\end{aligned}
\end{equation}
where we used the continuity of the trace operator and Poincar\'e's inequality for the second estimate. Combining this with \eqref{estdivk} and \eqref{coercive} we obtain
 \begin{equation*}
 \begin{aligned}
\mu_b^\frac12 \|\nabla u_1\|_{2} + \mu^\frac12 \|\nabla^g v_1\|_{2} &\le c\Big(\frac{1}{\mu_b^\frac12}\|f_1\|_{-1,2}+\big(\mu_b^\frac12+c_{-\frac12}\,\mu^\frac12\big)\|f_2\|_{2}\\
&\quad+\frac{1}{\mu_b^\frac12}\|f_3\|_{-\frac12,2}+\big(c_{\frac12}\,\mu_b^\frac12+\,\mu^\frac12\big)\|f_4\|_{2}\Big).
\end{aligned}
\end{equation*}
% \begin{equation}\label{estbuu}
% B(u,u)\le \frac{c_2}{\mu_b}\|f_1\|_{L_2(\Omega\setminus\Gamma)}^2+c\,\mu_b\,\|f_2\|_{L_2(\Omega\setminus\Gamma)}^2+\frac{c_1}{\mu_b}\|f_3\|_{L_2(\Gamma)}^2+c\,\mu\,\|f_4\|_{L_2(\Gamma)}^2.
% \end{equation}
Let $u:=u_1+u_0$. Then, from the last inequality and \eqref{estdivk} we obtain the part of \eqref{est1} estimating the velocities. By Corollary \ref{cor:pressure}, there exists a unique $(\pi,q)\in Z$ such that $\underline\nabla(\pi,q)=B(u_1,\cdot)-F_{u_0}$. We have
 \begin{equation}\label{bmf}
 \begin{aligned}
|B(u_1,\phi&)-F_{u_0}(\phi)|\\
&=\int_{\Omega\setminus\Gamma}\big(2\mu_b\langle Du,D\phi\rangle_e- f_1\cdot\phi\big)\,dx+\int_{\Gamma}\big(2\mu\langle \D u,\D\phi\rangle_g-f_3\cdot\phi\big)\,dA\\
&\le 2\mu_b\|\nabla u\|_2\|\nabla\phi\|_2 +\|f_1\|_{-1,2}\|\nabla\phi\|_2 +\|f_3\|_{-\frac12,2}\|\nabla\phi\|_2\\
&\quad+ 2\mu\big(\|\nabla^g v\|_2+\|k\|_\infty\|w\|_2\big)(\|\nabla^g P_\Gamma[\phi]_\Gamma\|_2+\|k\|_\infty\|[\phi]_\Gamma\|_2\big)\\
&\le c\Big(\mu_b^\frac12\|\nabla u\|_2+\mu^\frac12\|\nabla^g v\|_2+\frac{1}{\mu_b^\frac12}\|f_1\|_{-1,2}+\frac{1}{\mu_b^\frac12}\|f_3\|_{-\frac12,2}\Big)\\
&\quad \times\Big(\mu_b^\frac12+\mu^\frac12c_{-\frac12}+\frac{\mu}{\mu_b^\frac12}c_{-1}\Big)\|\phi\|_Y,
\end{aligned}
\end{equation}
where, again, we used the continuity of the trace operator and Poincar\'e's inequality for the second estimate. Note that so far the constants in our estimates do not depend on $\mu/\mu_b$. Finally, from the last estimate, the estimate of the velocities, and \eqref{estpressure} we obtain \eqref{est1}.
\qed

\medskip
We saw in the preceding proof that we can easily get estimates with constants independent of the Saffman-Delbr\"uck length $\mu/\mu_b$ and in which, in particular, we can pass to the limit $\mu\searrow 0$. Here, however, we are not interested in this limit; cf. the remark following the proof of Theorem \ref{thm:regs1}. Instead we try to keep the form of our estimates as simple as possible.

In a similar way, we can prove the analogous statement for the system \eqref{s2}.
\begin{theorem}
For all $(f_1,\ldots,f_5)\in\F^1_\nu(\Gamma)$ there exists a unique weak solution $(u,\pi,q)\in\E^0_\nu(\Gamma)$ of \eqref{s2}. We have
 \begin{equation}\label{est1nu}
 \begin{aligned}
&\mu_b^\frac12 \|\nabla u\|_{2} + \mu^\frac12 \|\nabla^g v\|_{2} + \frac{1}{\mu_b^\frac12} \|\pi\|_{2} + \frac{1}{\mu^\frac12} \|q\|_{2}\\
&\quad\le c\Big(\frac{1}{\mu_b^\frac12}\|f_1\|_{-1,2}+\mu_b^\frac12\|f_2\|_{2}+\frac{1}{\mu_b^\frac12}\|f_3\|_{-\frac12,2}+\mu^\frac12\|f_4\|_{2}\\
&\quad\quad + \frac{\mu_b}{\mu^\frac12}\|f_5\|_{2}+\mu_b^\frac12|f_5|_{\frac12,2}\Big).
\end{aligned}
\end{equation}
\end{theorem}
\proof We only give a brief sketch of the proof which proceeds analogously to the proof of Theorem \ref{thm:s1}. To begin with, let assume that $f_2$, $f_4$, and $f_5$ vanish. Instead of $X$, $Y$, and $Z$ we consider the spaces
\begin{equation*}
\begin{aligned}
X_\nu&:=\big\{u\in H^1_0(\Omega;\R^3)\,|\,\Div u=0 \text{ in }\Omega\setminus \Gamma,\ \DIV u=0 \text{ on } \Gamma,\ [u]_\Gamma\in H^1(\Gamma;T\Gamma)\big\},\\
Y_\nu&:=\big\{u\in H^1_0(\Omega;\R^3)\,|\, [u]_\Gamma\in H^1(\Gamma;T\Gamma)\big\},\\
Z_\nu&:=\big\{(f_2,f_4)\in L_2(\Omega)\times L_2(\Gamma)\ \ \text{such that \eqref{comp3} holds with $f_5=0$}\big\}.
\end{aligned}
\end{equation*}
Multiplying \eqref{s2} by $\phi\in X_\nu$, integrating, and integrating by parts we obtain again the weak formulation \eqref{weak}, that is, $B(u,\phi)=F(\phi)$ for all $\phi\in X_\nu$. Since $X_\nu\subset X$, Lemma \ref{lemma:coercive} shows that $B$ is coercive on $X_\nu$, hence giving the unique existence of a weak solution $u\in X_\nu$. Analogous to Lemma \ref{lemma:div} we have to show that the operator
\[\underline{\Div}: Y_\nu\rightarrow Z_\nu,\ \phi\mapsto(\Div\phi,\Div_g[\phi]_{\Gamma})\]
is surjective. This, however, is easy to see using $L_2$-theory for the Laplacian on $\Gamma$ and $L_2$-theory for the divergence operator in the bulk $\Omega\setminus\Gamma$. From this we conclude as before that $\underline\nabla:Z_\nu\rightarrow X_\nu^\perp\subset Y_\nu'$ is an isomorphism.
%Compared to Theorem \ref{thm:s1}, here, uniqueness of the pressure functions holds in a different sense because, for vanishing data, instead of \eqref{uibp} we have
% \begin{equation*}
% 0=\int_{\Gamma}q\,\Div_g\phi\,dA
% \end{equation*}
%for $\phi\in Y_\nu$.
In order to deal with full data, we construct an extension $u'\in H^1_0(\Omega)$ such that $[u']_\Gamma=f_5\,\nu$. Then, we let $u''\in Y_\nu$ solve $\underline\Div\, u''=(f_2-\Div u',f_4-\DIV u')$ and define $u_0:=u'+u''$. Finally, we let $u_1\in X_\nu$ solve $B(u_1,\phi)=F_{u_0}(\phi)$ for all $\phi\in X_\nu$ for the same $F_{u_0}$ as in \eqref{fun}, reconstruct the pressure functions using the $\underline\nabla$-operator, and define $u:=u_1+u_0$.
\qed

\medskip
Consider the subspace $U_\nu\subset L_2(\Omega)\times L_2(\Gamma)$ defined as follows: $(\pi,q)\in U_\nu$ iff $\pi=\kappa_i$ in $\Omega^i$, $\pi=\kappa_0$ in $\Omega^0$, $q=\kappa^i$ on $\Gamma^i$ with $\kappa_i,\kappa_0,\kappa^i\in\R$ for $i=1,\ldots,k$. Note that $Z_\nu\simeq Z'_\nu\simeq (L_2(\Omega)\times L_2(\Gamma))/U_\nu$. Again, the fact that the gauge freedom of the pressure functions is characterized by $U_\nu$ can be seen directly. Indeed, $\underline\nabla(\pi,q)(\phi)=0$ for all $Y_\nu$ implies that $\pi$ is constant in each connected component of $\Omega\setminus\Gamma$. Furthermore, instead of \eqref{uibp}, here we obtain
\begin{equation*}
 0=\int_{\Gamma}q\,\Div_g\phi\,dA,
 \end{equation*}
proving that $q$ is constant in each connected component of $\Gamma$.

Now, let us proceed with higher order estimates. 
%Separating the tangential and the normal part of \eqref{eqn:s1}$_3$, we obtain
%\begin{equation}\label{eqn:divt}
%\begin{aligned}
%-\grad_g q + \mu\,\big(\Delta_g v + \grad_g(w\,H) + K v -2\Div_g(w\,k)\big)&=-2\mu_b[\![Du]\!]\nu,\\
%-q\,H\nu + 2\mu\big(\langle\nabla^g v,k\rangle_g-w\,(H^2-2K)\big)&\\
%-\kappa\big(\Delta_g H + H(H^2/2-2K)+C_0(2K-HC_0/2)\big)&=[\![\pi]\!].
%\end{aligned}
%\end{equation}
%Note that $[\![Du]\!]\nu$ is tangential due to the incompressibility constraint.\marginpar{Genauer?}
Separating the tangential and the normal part of \eqref{s1}$_3$ by using the equation \eqref{fullt}, we see that the system \eqref{s1} can be written in the form
\begin{equation}\label{eqn:explicitsystem}
\begin{aligned}
\mu_b\Delta u-\grad\pi&=f_1&&\mbox{ in }\Omega\setminus\Gamma,\\
\Div u&=f_2&&\mbox{ in }\Omega\setminus\Gamma,\\
\mu\big(\Delta_{g} v + \grad_{g}(w\,H) + K\,v -2\Div_{g}(w\,k)\big)&\\
-\grad_{g} q  
+2\mu_b[\![D u]\!]\nu&=P_\Gamma f_3&&\mbox{ on }\Gamma,\\
2\mu\big(\langle\nabla^{g} v,k\rangle_{g}
-w\,(H^2-2K)\big)-q\,H-[\![\pi]\!]&=f_3\cdot\nu&&\mbox{ on }\Gamma,\\
\Div_{g} v-w\, H&=f_4&&\mbox{ on }\Gamma,\\
u&=0&&\mbox{ on }\pa\Omega.
\end{aligned}
\end{equation}
Note that $[\![Du]\!]\nu$ is tangential due to the incompressibility constraint. Indeed, for any vector $X$ on $\Gamma$ we have $[\![(X\cdot\nabla)u]\!]\cdot\nu=0$. If $X$ is tangential, we even have $[\![(X\cdot\nabla)u]\!]=0$ since $u$ is continuous across $\Gamma$. But then, choosing an orthonormal basis $\nu,e_1,e_2$ at some arbitrary point on $\Gamma$, from $\Div u=0$ we deduce that
\[[\![(\nu\cdot\nabla)u]\!]\cdot\nu=-[\![(e_1\cdot\nabla)u]\!]\cdot e_1-[\![(e_2\cdot\nabla)u]\!]\cdot e_2=0.\]

\begin{theorem}\label{thm:regs1}
Let $s\ge 1$. For all $(f_1,\ldots,f_4)\in\F^s(\Gamma)$ there exists a unique solution $(u,\pi,q)\in\E^s(\Gamma)$ of \eqref{s1}. We have
\begin{equation}\label{est:regs1}
\begin{aligned}
&\mu_b^\frac12|u|_{s+1,2} + \mu^\frac12 |v|_{s+1,2} + \frac{1}{\mu_b^\frac12} |\pi|_{s,2} + \frac{1}{\mu^\frac12} |q|_{s,2}\\
&\quad\le c\Big(\frac{1}{\mu_b^\frac12}|f_1|_{s-1,2}+ \mu_b^\frac12|f_2|_{s,2}+\frac{1}{\mu^\frac12}|P_\Gamma f_3|_{s-1,2} + \frac{1}{\mu_b^\frac12} |f_3\cdot\nu|_{s-\frac12,2}\\
&\quad\quad+ \mu^\frac12|f_4|_{s,2}\Big) + c_{-s}\Big(\mu_b^\frac12\|\nabla u\|_{2} + \mu^\frac12\|\nabla^g v\|_{2} + \frac{1}{\mu_b^\frac12}\|\pi\|_{2}+\frac{1}{\mu^\frac12}\|q\|_{2}\Big).
\end{aligned}
\end{equation}
\end{theorem}
\proof Theorem \ref{thm:s1} gives the existence of a unique weak solution $(u,\pi,q)\in\E^0(\Gamma)$. For two reasons we will only give a brief sketch of the proof of regularity. First, the general procedure of localization, transformation, and taking difference quotients is rather classical. And second, when generalizing the present result to the case of a variable membrane, we will have to repeat most of the arguments anyway. 

For integer $s$ we can prove by the techniques just mentioned that our solution lies in $\E^s(\Gamma)$ and
\begin{equation}\label{sest}
\|(u,\pi,q)\|_{\E^s(\Gamma)}\le c\,\|(f_1,\ldots,f_4)\|_{\F^s(\Gamma)}. 
\end{equation}
By interpolation, the same is true for arbitrary $s\ge 1$. In order to prove the estimate \eqref{est:regs1}, we only have to replace the lower order parts of norm of the data in \eqref{sest} by the terms in the second bracket on the right hand side of \eqref{est:regs1}. This, however, follows easily by using the equations \eqref{eqn:explicitsystem} and by interpolation and absorption. We present the idea by considering the term $f_1$. Let $s>1$. For arbitrary integer $0\le l<s-1$, we have 
\begin{equation*}%\label{est:regs1}
\begin{aligned}
|f_1|_{l}\le c\big(\mu_b\|\nabla^{2+l} u\|_{2}+\|\nabla^{1+l} \pi\|_{2}\big)\le c\big(\mu_b|u|_{s+1}^{\theta_0}\|\nabla u\|_{2}^{1-\theta_0}+|\pi|_{s}^{\theta_1}\|\pi\|_{2}^{1-\theta_1}\big).
\end{aligned}
\end{equation*}
Using Young's inequality and absorbtion, the claim follows. The other terms $f_2,\ldots,f_4$ can be handled analogously.

In the proof of \eqref{sest} we consider only the case $s=1$. The case $s\in\N$ then follows by differentiating the equations in the prototype geometries (whole space, half space, and double half space). Furthermore, for the sake of a short presentation, we assume that $(u,\pi,q)\in\E^1(\Gamma)$, that is, we only prove a-priori estimates. By a standard localization argument, we may assume that the solution is supported in an open cube $Q_R$ of side length $R>0$. We only consider the case that $Q_R$ is centered at some point $x_0\in\Gamma$; the other cases, namely $Q_R$ being centered at some point $x_0\in\pa\Omega$ and $\overline Q_R$ being contained in $\Omega\setminus\Gamma$, only involve standard analysis of the classical Stokes system. Rotating and translating the Cartesian coordinate system and choosing $R$ smaller if necessary, we may assume that $x_0=0$ and that $\Gamma\cap Q_R$ is the graph of a 
smooth function $h: Q_R^2:=Q_R\cap (\R^2\times\{0\})\rightarrow (-R/2,R/2)$ such that $h(0)=0$ and $\nabla h(0)=0$. Consider the smooth diffeomorphism 
\[\Phi^{-1}_h: Q_R\rightarrow \tilde  Q_R:=\Phi^{-1}_h(Q_R), (x',x_3)\mapsto (x',x_3-h(x')).\]
This diffeomorphism induces the metric $\tilde e:=\Phi_h^* e$ on $\tilde Q_R$. We denote the restriction of $\tilde e$ to $Q_R^2$ by $\tilde g$. Note that $\Phi_h:(\tilde Q_R,\tilde e)\rightarrow (Q_R,e)$ and $\Phi_h|_{Q_R^2}:(Q_R^2,\tilde g)\rightarrow (\Gamma\cap Q_R,g)$ are isometries. Let us denote the pullbacks of the involved fields by $\tilde u:=\Phi_h^*u$, $\tilde\pi:=\Phi_h^*\pi$, $\tilde v:=\Phi_h^*v$, $\tilde w:=\Phi_h^*w$, $\tilde q:=\Phi_h^*q$, $\tilde f_3^\top:=\Phi_h^*(P_\Gamma f_3)$, $\tilde f_3^\perp:=\Phi_h^*(f_3\cdot\nu)$, and $\tilde f_i=\Phi_h^*f_i$ for $i=1,2,4$. By exploiting naturality of covariant differentiation under isometries, from \eqref{eqn:explicitsystem} we obtain
\begin{equation}\label{eqn:tildesystem}
\begin{aligned}
\mu_b\Delta_{\tilde e}\tilde u-\grad_{\tilde e}\tilde\pi&=\tilde f_1&&\mbox{ in }\tilde Q_R\setminus Q_R^2,\\
\Div_{\tilde e} \tilde u&=\tilde f_2&&\mbox{ in }\tilde Q_R\setminus Q_R^2,\\
\mu\big(\Delta_{\tilde g} \tilde v + \grad_{\tilde g}(\tilde w\,H_{\tilde e}) + K_{\tilde g}\,\tilde v -2\Div_{\tilde g}(\tilde w\,k_{\tilde e})\big)&\\
-\grad_{\tilde g} \tilde q  
+2\mu_b[\![D^{\tilde e}\tilde u]\!]\nu_{\tilde e}&=\tilde f_3^\top&&\mbox{ on } Q_R^2,\\
2\mu\big(\langle\nabla^{\tilde g} \tilde v,k_{\tilde e}\rangle_{\tilde g}
-\tilde w\,(H^2_{\tilde e}-2K_{\tilde g})\big)-\tilde q\,H_{\tilde e}-[\![\tilde\pi]\!]&=\tilde f_3^\perp&&\mbox{ on }Q_R^2,\\
\Div_{\tilde g} \tilde v-\tilde w\, H_{\tilde e}&=\tilde f_4&&\mbox{ on }Q_R^2,\\
\tilde u-\tilde v-\tilde w\,\nu_{\tilde e}&=0&&\mbox{ on }Q_R^2.
\end{aligned}
\end{equation}
Here, we included the transformation of the identity $u=v+w\,\nu$ on $\Gamma$; the reason for this will become clear in a moment. The results from \Appref{app:cov} show that \eqref{eqn:tildesystem} can be written in the form
\begin{equation}\label{eqn:trans}
\begin{aligned}
\mu_b\Delta\tilde u-\grad\tilde\pi&=\hat f_1&&\mbox{ in }\tilde Q_R\setminus Q_R^2,\\
\Div \tilde u&=\hat f_2&&\mbox{ in }\tilde Q_R\setminus Q_R^2,\\
\mu\Delta \tilde v -\grad \tilde q  
+2\mu_b[\![D\tilde u]\!]e_3&=\hat f_3^\top&&\mbox{ on }Q_R^2,\\
-[\![\tilde\pi]\!]&=\hat f_3^\perp&&\mbox{ on } Q_R^2,\\
\Div \tilde v&=\hat f_4&&\mbox{ on } Q_R^2,\\
\tilde u -\tilde v - \tilde w\,e_3 &= \hat f_5&&\mbox{ on }Q_R^2,
\end{aligned}
\end{equation}
where $e_3:=\Phi_h^*\nu$ is the normal to $\R^2\times\{0\}$ and
\begin{equation*}
\begin{aligned}
\hat f_1&=\tilde f_1 + (\tilde e-e)*(\mu_b\nabla^2 \tilde u,\grad\tilde\pi) + \mu_b\, r(\tilde e)*\big((\nabla^2\tilde e,(\nabla\tilde e)^2)*\tilde u+\nabla\tilde e*\nabla \tilde u\big),\\
\hat f_2&=\tilde f_2+r(\tilde e)*\nabla\tilde e*\tilde u,\\
\hat f_3^\top
% &=\tilde f_3^\top+(\tilde g-g)*(\mu(\nabla^g)^2 \tilde v,\grad_g \tilde q) + \mu\,r(\tilde g,\hat g)*\big((\nabla^{g})^2\tilde g*\tilde v,(\nabla^{g}\tilde g)^2*\tilde v,\nabla^g\tilde g*\nabla^g \tilde v\big)\\
% &\quad + (w,v)*r(\tilde e,\hat e)*(k^2,\nabla k,k*\nabla\tilde e,\nabla^2\tilde e,(\nabla\tilde e)^2) + \nabla w*r(\tilde e,\hat e)*(k,\nabla\tilde e)\big)\\
% &\quad + \mu_b\,r(\tilde e,\hat e)*([\nabla \tilde u],\nabla\tilde e*[\tilde u])\\
&=\tilde f_3^\top+(\tilde e-e)*(\mu(\nabla^g)^2 \tilde v,\grad_g \tilde q) + \mu_b\,r(\tilde e)*\big([\nabla \tilde u]+\nabla\tilde e*[\tilde u]\big)\\
&\quad+ \mu\,r(\tilde e)*\big((\nabla^2\tilde e,(\nabla\tilde e)^2)*[\tilde u] + \nabla\tilde e*[\nabla\tilde u]\big),\\
\hat f_3^\perp&=\tilde f_3^\perp + \mu\,r(\tilde e)*\big(\nabla\tilde e*\nabla^g\tilde v + (\nabla\tilde e)^2*[\tilde u]\big)
% &=\tilde f_3^\perp + \mu\,r(\tilde e,\hat e)*\big((k,\nabla\tilde e)*\nabla^g\tilde v + \nabla \tilde e*(k,\nabla\tilde e)*\tilde v + w\,(k^2,k*\nabla\tilde e,(\nabla\tilde e)^2)\big)\\
+ r(\tilde e)*\nabla\tilde e\,\tilde q,\\
\hat f_4&=\tilde f_4 + r(\tilde e)*\nabla\tilde e*[\tilde u],\\
\hat f_5&=(\tilde e-e)*r(\tilde e)\,\tilde w.
\end{aligned}
\end{equation*}
In fact, the precise form of the functions $\hat f_1,\ldots,\hat f_5$ given above is not so relevant here (except for the terms involving highest order derivatives); it will however be crucial when we generalize the present result to the case of a variable membrane. It is not hard to see that
\begin{equation}\label{eme}
\tilde e (x',x_3)-e=r(\nabla h(x')).
\end{equation}
% \begin{equation*}
% \begin{aligned}
% \tilde f_1&=\hat f_1 + (\tilde e-e)*(\nabla^2 \tilde u,\grad\tilde\pi) + p(\tilde e)*\big(\nabla^2\tilde e*\tilde u,(\nabla\tilde e)^2*\tilde u,\nabla\tilde e*\nabla \tilde u\big),\\
% \tilde f_2&=\tilde f_2+\tilde e*\nabla \tilde e*\tilde u,\\
% \tilde f_3^\top&=\hat f_3^\top+(\tilde g-g)*((\nabla^g)^2 \tilde v,\grad_g q) + p(\tilde g)*\big((\nabla^g)^2\tilde g*\tilde v,(\nabla^g\tilde g)^2*\tilde v,\nabla^g\tilde g*\nabla^g \tilde v\big),\\
% \tilde f_3^\perp&=\hat f_3 + r(1,h,\bar \nabla^2 h)*(\tilde q,[\tilde u]_{\Gamma},[\nabla\tilde u]_{\Gamma}),\\
% \tilde f_4&=\hat f_4 + r(h,\bar\nabla^2 h)*[\tilde u]_{\Gamma},\\
% \tilde f_5&=\tilde w\,r(h,\nabla h).
% \end{aligned}
% \end{equation*}
% Here and in the following, we denote by $r(a)$ generic tensor fields that are analytic functions of their argument $a$ such that $r(0)=0$, for tensor fields $r_1$ and $r_2$ we write $r_1*r_2$ for any tensor field that depends in a bilinear way on $r_1$ and $r_2$, and we use the abbreviations $\bar\nabla^k h:=(\nabla h,(\nabla^g)^2h,\ldots,(\nabla^g)^kh)$ and $r*(r_1,r_2,\dots):=r*r_1 + r*r_2 +\ldots$ Finally, we note that
% \begin{equation}
%  \tilde e - e=r(\nabla\Phi_h)=r(h,\nabla h), \quad \tilde g - g=r(h,\nabla h).
% \end{equation}
We may assume without restriction that $\hat f_2$, $\hat f_4$, and $\hat f_5$ vanish; cf. the proof of Lemma \ref{lemma:div}. We multiply equation \eqref{eqn:trans}$_1$ by $\pa_i^2 \tilde u$, where $i=1,2$, integrate over $\tilde Q_R$, and integrate by parts to obtain
\begin{equation*}
\begin{aligned}
\mu_b\int_{\tilde Q_R\setminus Q^2_R} |\nabla \pa_i \tilde u|^2\,dx + \mu\int_{Q_R^2} &|\nabla \pa_i \tilde v|^2\,dA \le \int_{\tilde Q_R} \hat f_1\cdot \pa_i^2\tilde u\,dx + \int_{Q_R^2} \hat f_3\cdot\pa_i^2\tilde u\,dA\\
&\le \|\hat f_1\|_{2}\|\pa_i^2\tilde u\|_{2} + \|\hat f_3^\top\|_{2}\|\pa_i^2\tilde v\|_{2}
+ c|\hat f_3^\perp|_{\frac12}|\pa_i\tilde w|_{\frac12}\\
&\le c\big(\|\hat f_1\|_{2} + |\hat f_3^\perp|_{\frac12}\big)\|\nabla\pa_i\tilde u\|_{2} + \|\hat f_3^\top\|_{2}\|\pa_i^2\tilde v\|_{2}.
 \end{aligned}
\end{equation*}
For the second inequality we used the Fourier transform and Plancherel's identity, while for the last inequality we used the mapping properties of the trace operator and Poincar\'e's inequality; recall that all involved functions are supported in $\tilde Q_R$. Using Young's inequality, absorption, and a standard trick from the theory of the classical Stokes system to obtain estimates of \emph{all} derivatives, see for instance \cite{Galdi:Navier-Stokes-1}, we infer that
\begin{equation}\label{est:regs1med}
\begin{aligned}
&\mu_b^\frac12|\tilde u|_{2,2} + \mu^\frac12 |\tilde v|_{2,2} + \frac{1}{\mu_b^\frac12} \|\tilde \pi\|_{2,2} + \frac{1}{\mu^\frac12} \|\tilde q\|_{2,2}\\
&\quad\le c\Big(\frac{1}{\mu_b^\frac12}\|\hat f_1\|_{2}+\frac{1}{\mu^\frac12}\|\hat f_3^\top\|_{2} + \frac{1}{\mu_b^\frac12} |\hat f_3^\perp|_{\frac12,2}\Big).
\end{aligned}
\end{equation}
% Using Young's inequality, absorption, \eqref{eqn:trans}$_{1,2,3}$, and again Poincar\'e we obtain
% \begin{equation}\label{est:h2}
% \begin{aligned}
% \|u\|_{H^2(\tilde Q_R\setminus Q_R^2)}^2& + \|v\|_{H^2(Q^2_R)}^2 + \|\grad_e\pi\|_{L_2(Q_R)}^2 +  \|q\|_{H^1(Q_R^2)}^2\\
% &\le c\,\big(\|\tilde f_1\|_{L_2(\tilde Q_R\setminus Q^2_R)}^2 + \|\tilde f_2\|_{H^1(\tilde Q_R\setminus Q^2_R)}^2 + \|\tilde f_3^\top\|_{L_2(Q_R^2)}^2\\
% &\quad + \|\tilde f_3^\perp\|_{H^{1/2}(Q_R^2)}^2 + \|\tilde f_4\|_{H^1(Q_R^2)}^2+ \|\tilde f_5\|_{H^{3/2}(Q_R^2)}^2\big).
% \end{aligned}
% \end{equation}
Due to \eqref{eme}, for $R\searrow 0$, the factor $\tilde e-e$ gets arbitrarily small in $L^\infty$. Thus, we can absorb the highest order terms in $\hat f_1$ and $\hat f_3^\top$ (and $\hat f_5$ for full data) on the left hand side of \eqref{est:regs1med}. 
% Then, by interpolation and Poincar\'e's inequality we can also get rid of the lower order terms, giving \eqref{est:h2} with $\hat f_i$ instead of $\tilde f_i$ and $\tilde f_5=0$. 
Transforming \eqref{est:regs1med} back to $Q_R$, summing the finite number of local estimates, and eliminating the lower order terms on the right hand side by using interpolation and \eqref{est1} (or a standard contradiction argument) we obtain \eqref{sest} for $s=1$.
\qed

\medskip
In the preceding proof we took the geometric pullback of the unkowns when transforming the system to the prototype geometry. Apart from being very convenient for computations, here, this procedure is indispensible since the normal and the tangential part of the velocity field on $\Gamma$ have different regularities. Without this geometric transformation we would not be able to control the perturbations in $\hat f_3^\top$ and $\hat f_3^\perp$. 

The constants in \eqref{est:regs1} depend on $\mu/\mu_b$; in particular, the estimate doesn't hold uniformly for $\mu\searrow 0$. The reason for this is not that our proof is too rough; for $\mu=0$ we have to assume $P_\Gamma f_3\in H^{s-\frac12}(\Gamma)$ in order to obtain $u\in H^{s+1}(\Omega\setminus\Gamma)$. Of course, it is possible to establish estimates that hold uniformly for $\mu\searrow 0$. However, since then \eqref{eqn:explicitsystem}$_3$ has to hold in $H^{s-\frac12}(\Gamma)$ instead of $H^{s-1}(\Gamma)$ such estimates need a higher regularity for the membrane than \eqref{est:regs1}; cf. the next subsection. For this reason we refrain from giving uniform estimates.

In a very similar way we can prove the following theorem.
\begin{theorem}
Let $s\ge 1$. For all $(f_1,\ldots,f_5)\in\F^s_\nu(\Gamma)$ there exists a unique solution $(u,\pi,q)\in\E^s_\nu(\Gamma)$ of \eqref{s2}. We have
\begin{equation}\label{est:regs2}
\begin{aligned}
&\mu_b^\frac12|u|_{s+1,2} + \mu^\frac12 |v|_{s+1,2} + \frac{1}{\mu_b^\frac12} |\pi|_{s,2} + \frac{1}{\mu^\frac12} |q|_{s,2}\\
&\quad\le c\Big(\frac{1}{\mu_b^\frac12}|f_1|_{s-1,2}+ \mu_b^\frac12|f_2|_{s,2}+\frac{1}{\mu^\frac12}|f_3|_{s-1,2} + \mu^\frac12|f_4|_{s,2}+\mu_b^\frac12|f_5|_{s+\frac12,2}\Big)\\
&\quad\quad + c_{-s}\Big(\mu_b^\frac12\|\nabla u\|_{2} + \mu^\frac12\|\nabla^g v\|_{2} + \frac{1}{\mu_b^\frac12}\|\pi\|_{2}+\frac{1}{\mu^\frac12}\|q\|_{2}\Big).
\end{aligned}
\end{equation}
\end{theorem}
\medskip

\subsection{Analysis for variable membranes}\subseclabel{stokesvar}

Now, we transfer the results from the preceding subsection to the case of variable membranes. More precisely, we will analyze the dependence of the constants in the estimates \eqref{est1}, \eqref{est1nu}, \eqref{est:regs1}, and \eqref{est:regs2} on the membrane. We will represent variations of the membrane in the following way. Let us denote by $S_{\alpha}$, $\alpha>0$, the open set of points in
$\Omega$ whose distance from $\Gamma$ is less than $\alpha$. It's a
well-known fact from elementary differential geometry that there
exists a maximal $\kappa>0$ such that the mapping
\begin{equation*}
 \begin{aligned}
\Lambda: \Gamma\times (-\kappa,\kappa)\rightarrow
S_{\kappa},\
(x,d)\mapsto x + d\,\nu(x)  
 \end{aligned}
\end{equation*}
is a diffeomorphism. For functions $h:\Gamma\rightarrow (-\kappa,\kappa)$ we define $\Gamma_h:=\{\Lambda(x,h(x))\,|\,x\in\Gamma\}$, and we write $x\mapsto(\tau(x),d(x))$ for the inverse mapping $\Lambda^{-1}$. Analogously to the notation used before, we denote by $\Gamma_h^i$ and $\Omega_h^i$, $i=1,\ldots,k$, the connected components of $\Gamma_h$ and the open sets enclosed by $\Gamma_h^i$, respectively, and by $\Omega_h^0$ the exterior part of $\Omega\setminus\Gamma_h$.

The analysis will be based on the \emph{Hanzawa transform} $\Phi_h$ which maps $\overline\Omega$ diffeomorphically to itself such that $\Phi_h(\Gamma)=\Gamma_h$. We choose a real-valued function $\beta\in C^\infty(\R)$ that is $0$ in neighborhoods of $-1$ and $1$, and $1$ in a neighborhood of $0$ such that $|\beta'|< \kappa/\|h\|_{L^\infty(\Gamma)}$ on $\Gamma$.  While, in $\overline\Omega\setminus S_\kappa$, we let $\Phi_h$ be the identity, we define $\Phi_h$ in $S_\kappa$ by
\begin{equation*}
 \begin{aligned}
x&\mapsto x + \nu(\tau(x))\,h(\tau(x))\,\beta(d(x)/\kappa).
\end{aligned}
\end{equation*}
It is not hard to prove that $\Phi_h:\overline\Omega\rightarrow\overline\Omega$ and $\varphi_h:=\Phi_h|_{\Gamma}:\Gamma\rightarrow\Gamma_h$ are diffeomorphisms; see for instance \cite{lengelerphd}. Throughout the whole subsection we denote by $r_h$ analytic functions of $\|h/\kappa\|_{\infty}$, $\|\nabla h\|_{\infty}$, and, in the case of Theorem \ref{s1h} and Theorem \ref{s2h}, of $|h|_{2,2}$ that, additionally, may depend on $\Omega$, $\Gamma$, $\mu/\mu_b$ and $s$. 
%Note that the term $\|h/\kappa\|_\infty$ controls the quantity $\|hk\|_{\infty}$, where $k$ is the second fundamental form of $\Gamma$, since $\kappa$ is bounded by the inverse of the maximal principal curvature of $\Gamma$; see for instance \cite{lengeler14}. 
Furthermore, we write $k_h$, $H_h$, etc. for the geometric quantities corresponding to $\Gamma_h$.

The analysis of the constants in \eqref{est:regs1} and \eqref{est:regs2} will proceed by a perturbation argument. However, this argument inevitably produces lower order terms on the right hand sides of the estimates; for this reason  the analysis of \eqref{est1} and \eqref{est1nu} has to be done in a different, more direct way.\footnote{We could also argue by contradiction to prove that \eqref{est1} and \eqref{est1nu} hold uniformly as long as $h$ remains bounded in a suitable norm. Such an \emph{abstract} result, however, is insufficient, for instance, in the analysis of singularities in the dynamical problems. Furthermore, a Neumann series argument would need the assumption of $h$ being small in some high regularity class which, again, is insufficient for our purposes.} We will need the following lemmas.
\begin{lemma}[Trace operator]\label{lemma:trace}
 Let $h\in W^{1,\infty}(\Gamma)$ such that $\|h\|_{L^\infty(\Gamma)}<\kappa/2$. For $u\in H^1_0(\Omega)$ we have
 \begin{equation*}
  \|[u]_{\Gamma_h}\|_4+|[u]_{\Gamma_h}|_{\frac12,2}\le r_h \|\nabla u\|_2;
 \end{equation*}
 in particular
 \begin{equation}\label{tr22}
  \|[u]_{\Gamma_h}\|_2\le r_h\, |\Gamma|^{\frac14} \|\nabla u\|_2.
 \end{equation}
\end{lemma}
\proof To begin with, let us assume that $h\equiv 0$. Then, it is well-know that
\begin{equation}\label{tr42}
 \|[u]_{\Gamma}\|_4\le c\,\|\nabla u\|_2;
\end{equation}
see for instance \cite{kufner13}. For non-vanishing $h$ we note that
\begin{equation*}
\begin{aligned}
\int_{\Gamma_h} |u|^4\,dA_h &= \int_\Gamma |u\circ\phi_h|^4 |\det{d\phi_h}|\,dA\le r_h \int_\Gamma |u\circ\phi_h|^4\,dA\\
&\le r_h\Big(\int_{\Omega} |\nabla(u\circ\Phi_h)|^2\,dx\Big)^2\le r_h\Big(\int_{\Omega} |(\nabla u)\circ\Phi_h)|^2\,dx\Big)^2\\
&=r_h\Big(\int_{\Omega} |\nabla u|^2|\det{d(\Phi_h^{-1})}|\,dx\Big)^2\le r_h\Big(\int_{\Omega} |\nabla u|^2\,dx\Big)^2.
\end{aligned}
\end{equation*}
Here, the determinant $\det{d\phi_h}$ has to be taken with respect to orthonormal bases in the respective tangent spaces. For the inequalities we used \eqref{tr42} and the obvious facts that $d\phi_h$ and $d(\Phi_h^{-1})$ can be bounded in $L^\infty$ by $r_h$. In particular, we have $|\Gamma_h|\le r_h|\Gamma|$, proving \eqref{tr22}. The estimate of the trace in $H^{1/2}(\Gamma_h)$ proceeds very similarly.
\qed

\medskip
\begin{lemma}[Divergence equation]\label{lemma:divh}
 Let $h\in W^{1,\infty}(\Gamma)$ such that $\|h\|_{\infty}<\kappa/2$, and let $f\in L_2(\Omega)$ be mean value free in each connected component of $\Omega\setminus\Gamma_h$. Then there exists a function $u\in H^1_0(\Omega\setminus\Gamma_h)$ such that $\Div u=f$ in $\Omega\setminus\Gamma_h$ and
\begin{equation}\label{divi}
\|\nabla u\|_2\le r_h\|f\|_2.
\end{equation}
\end{lemma}
\proof Theorem III.3.1 in \cite{Galdi:Navier-Stokes-1} gives the existence of a solution. This theorem also provides an upper bound for the constant in \eqref{divi}.\footnote{See \cite{duran12} for an improvement of this bound.} In order to prove \eqref{divi} we have to cover $\Omega\setminus\Gamma_h$ suitably by sets each of which is star-shaped with respect to a ball. However, we will only briefly sketch this procedure since the article already is quite long.
%The explicit constuction of the solution of the divergence equation given in \cite{duran06} in fact shows that the constant of continuity of the solution operator can be bounded by an analytic function of the John constants which itself can easily be bounded by $r_h$.
In fact, we solve the divergence equation in each connected component of $\Omega\setminus\Gamma_h$ separately. Hence, let us fix some $i\in\{0,\ldots,k\}$ and consider the set $\Omega_h^i$. We cover $\Omega_h^i\setminus S_{\frac34 \kappa}$ by finitely many balls whose closures are contained in $S_{\frac{\kappa}{2}}$; these balls do not depend on $h$. Now, let us fix a point $x\in \Gamma^i$, let $Z$ denote the cylinder of length $2\kappa$ and radius $r$ that is centered at $x$ and perpendicular to $T_x\Gamma^i$, and let $B\subset\Omega^i$ be the ball of radius $r/2$ centered in $Z$ at a distance $3\kappa/4$ from $\Gamma^i$. We claim that 
$Z\cap\Omega^i_h$ is star-shaped with respect to the ball $B$ provided that $r=1/c_1\,(1+\|\nabla h\|_{L_\infty(\Gamma)})^{-1}$ for a sufficiently large constant $c_1>0$. Indeed, if $\Gamma^i$ is flat in a neighbourhood of $x$ we may choose $r=\kappa/4\,(1+\|\nabla h\|_{L_\infty(\Gamma)})^{-1}$, since then, as is not hard to see, $B$ is contained in each cone $K\subset\Omega^i_h$ with its vertex on $\Gamma^i_h$ and its axis perpendicular to $T_x\Gamma^i$. On the other hand, if $\Gamma^i$ is not flat in a neighbourhood of $x$, then we reduce the situation to the flat case with the help of a flattening submanifold chart. Even though the cylinders, cones, and balls that we construct in this way are not Euclidean, they do contain corresponding Euclidean objects with smaller radii and open angles, respectively. This results in multiplying $r$ by an additional factor that depends on $\Gamma^i$ but is independent of $h$. Now, we cover $\Omega^i_h\cap S_\kappa$ by such cylinders whose inverse radii are bounded by $c_1r_h$. Finally, in order to control the constants $C_k$ in Lemma III.3.2 in \cite{Galdi:Navier-Stokes-1} we have to enumerate our covering sets suitably. Here, it is convenient to start the enumeration with the cylinders and to end it with the balls that cover $\Omega_h^i\setminus S_{\frac34 \kappa}$ and that do not depend on $h$; then it is easily possible to control the quantities $|F_k|$ in Lemma III.3.2 from below.
\qed

\medskip
\begin{lemma}[Poincar\'e's inequality]\label{lemma:poincare}
 Let $h\in W^{1,\infty}(\Gamma)$ such that $\|h\|_{\infty}<\kappa/2$. For $u\in H^1(\Gamma_h)$ we have
 \[\|u-\bar u\|_2\le r_h\, c_1\|\nabla u\|_2,\]
where $\bar u:=|\Gamma_h^i|^{-1}\int_{\Gamma_h^i}u\,dA_h$ on $\Gamma_h^i$ for all $i=1,\ldots,k$.
\end{lemma}
\proof This follows from Cheeger's inequality, see \cite{cheeger69}, stating that the smallest positive eigenvalue of the Laplacian on $\Gamma_h$ is bounded from below by $i(\Gamma_h)^2/4$, where 
\[i(\Gamma_h)=\inf\frac{|\gamma|}{\min(|A_1|,|A_2|)}.\]
Here, the infimum is taken over all simple closed curves $\gamma$ dividing $\Gamma_h$ into two subsets $A_1$ and $A_2$ with common boundary $\partial A_1=\partial A_2=\gamma$, and $|\gamma|, |A_i|$ denote the length of $\gamma$ and the area of $A_i$, respectively. From \eqref{krh} we see that $i(\Gamma_h)^{-1}$ is bounded by $r_h\, i(\Gamma)^{-1}$.
\qed

\medskip
We begin with the analysis of the constant in \eqref{est1}.
\begin{theorem}\label{thm:fh}
Let $h\in W^{2,\infty}(\Gamma)\cap H^{3}(\Gamma)$ such that $\|h\|_{L^\infty(\Gamma)}<\kappa/2$. For all $(f_1,\ldots,f_4)\in\F^1(\Gamma_h)$ there exists a unique weak solution $(u,\pi,q)\in\E^0(\Gamma_h)$ of \eqref{s1}. The estimate \eqref{est1} holds with the constant $c$ being bounded by
\begin{equation}\label{rhgamma}
\begin{aligned}
r_h\big(1+|h|_{2,\infty}^{15}+|h|_{3,2}^{15}\big)\Big(1+\sum_{i}\big(|\Gamma^i_h|(H_h,H_h)_{\Gamma_h^i}-(1,H_h)^2_{\Gamma^i_h}\big)^{-2}\Big)\,c_{19}.
\end{aligned}
\end{equation}
Here, the sum is taken over all $i$ such that $\Gamma^i_h$ is not a round sphere; note that, by Jensen's inequality, the denominators in the second bracket are positive.
\end{theorem}
\proof The proof consists of four steps. We define $X_h$, $Y_h$, $Z_h$, and $B_h$ just like $X$, $Y$, $Z$, and $B$, respectively, with $\Gamma$ replaced by $\Gamma_h$.

\medskip
\noindent (a) First, we need to analyze the constant $c$ in \eqref{coercive}$_2$. From the proof of Lemma \ref{lemma:coercive} and using Lemma \ref{lemma:trace}, for $u\in X_h$ we obtain
\begin{equation*}
\begin{aligned}
\|\nabla v\|_2^2&\le \frac2\mu B_h(u,u) + 16\|k_{h}\|_{\infty}^2 \|[u]_{\Gamma_h}\|_2^2\\
&\le \frac2\mu B_h(u,u) + r_h\|k_h\|_{\infty}^2 c_1\|\nabla u\|_2^2.
\end{aligned}
\end{equation*}
Thus, in view of \eqref{krh}, the definition of $B_h$, and \eqref{du}, we infer that
\begin{equation*}
\begin{aligned}
\|\nabla v\|_2^2&\le  \frac2\mu B(u,u) + r_h\big(1+|h|_{2,\infty}^2\big)\,c_1\|\nabla u\|_2^2\\
&\le r_h\big(1+|h|_{2,\infty}^2\big)\frac{c_2}{\mu} B_h(u,u).
\end{aligned}
\end{equation*}

\medskip
\noindent (b) Next, we analyze the constant $c$ occuring in \eqref{estdiv} and \eqref{estpressure}. Let $(f_2,f_4)\in Z_h$. For each $i=1,\ldots,k$ we choose $w=\lambda_1^i+\lambda_2^i |\Gamma_h^i|^{\frac12} H_h$ on $\Gamma_h^i$ such that
\begin{equation*}
\begin{aligned}
\lambda_1^i(1,H_h)_{\Gamma_h^i}+\lambda_2^i|\Gamma_h^i|^{\frac12}(H_h,H_h)_{\Gamma_h^i}&=-(1,f_4)_{\Gamma_h^i},\\
\lambda_1^i(1,1)_{\Gamma_h^i}+\lambda_2^i|\Gamma_h^i|^{\frac12}(1,H_h)_{\Gamma_h^i}&=(1,f_2)_{\Omega^i},
\end{aligned}
\end{equation*}
where $(\cdot,\cdot)_{\Gamma_h^i}$ and $(\cdot,\cdot)_{\Omega^i_h}$ denote the $L_2$-scalar products on $\Gamma_h^i$ and $\Omega^i_h$, respectively. If $\Gamma_h^i$ is a round sphere, we choose 
\[\lambda_1^i|\Gamma_h^i|:=(1,f_2)_{\Omega^i}\quad\text{and}\quad\lambda_2^i:=0.\]
If $\Gamma_h^i$ is a not a round sphere, then the linear system is regular, and we obtain
\begin{equation*}
\begin{aligned}
\lambda_1^i&=\frac{1}{|\Gamma^i_h|(H_h,H_h)_{\Gamma_h^i}-(1,H_h)^2_{\Gamma^i_h}}\big((1,f_4)_{\Gamma_h^i}(1,H_h)_{\Gamma_h^i}+(1,f_2)_{\Omega^i_h}(H_h,H_h)_{\Gamma_h^i}\big),\\
\lambda_2^i&=\frac{-|\Gamma^i_h|^{-\frac12}}{|\Gamma^i_h|(H_h,H_h)_{\Gamma_h^i}-(1,H_h)^2_{\Gamma_h^i}}\big((1,f_4)_{\Gamma_h^i}(1,1)_{\Gamma_h^i}+(1,f_2)_{\Omega^i_h}(1,H_h)_{\Gamma_h^i}\big).
\end{aligned}
\end{equation*}
%
%In the latter case we have
%\begin{equation*}
%\begin{aligned}
%|\lambda_1^i|&\le\frac{r_h}{(H^2)_{\Gamma_h^i}-(H)^2_{\Gamma_h^i}}\big((1+|h|_{2,\infty})|\Gamma^i|^{-\frac12}\|f_4\|_2+(1+|h|_{2,\infty}^2)|\Omega|^{\frac12}|\Gamma^i|^{-1}\|f_2\|_2\big),\\
%|\lambda_2^i|&\le\frac{r_h}{(H^2)_{\Gamma_h^i}-(H)^2_{\Gamma_h^i}}\big(|\Gamma^i|^{-\frac12}\|f_4\|_2+(1+|h|_{2,\infty})|\Omega|^{\frac12}|\Gamma^i|^{-1}\|f_2\|_2\big).
%\end{aligned}
%\end{equation*}
%We conclude that
%\begin{equation*}
%\begin{aligned}
%\|w\,H\|_2&\le \bar r_h\big((1+|h|_{2,\infty}^2)\|f_4\|_2+(1+|h|_{2,\infty}^3)(|\Omega|/|\Gamma|)^{\frac12}\|f_2\|_2\big).
%\end{aligned}
%\end{equation*}
Now, let $\psi\in H^2(\Gamma_h)$ be a function with vanishing mean value on $\Gamma^i_h$ solving $\Delta_g\psi=f_4+wH_h$ on $\Gamma_h^i$, for all $i=1,\ldots,k$. Note that $\psi_{;\alpha\beta}^{\quad\, \beta}=\psi_{;\beta\ \alpha}^{\ \ \beta}+\psi_{;\alpha}K$. Hence, for $v:=\grad_g\psi$, integration by parts and \eqref{krh} give
\begin{equation*}
\begin{aligned}
\|\nabla^g v\|_2^2&=\|\Delta_g\psi\|_2^2-\int_{\Gamma_h}K_h|v|^2\, dA_h\\
&\le c\big(\|f_4\|_2^2 + \|w H_h\|_2^2\big) + r_h\big(1+|h|_{2,\infty}^2\big)\|v\|_2^2.
\end{aligned}
\end{equation*}
However, by Lemma \ref{lemma:poincare} we have
\begin{equation}\label{vau2}
\begin{aligned}
\int_{\Gamma_h}|v|^2\,dA_h=-\int_{\Gamma_h}(f_4+wH_h)\psi\,dA_h&\le \|f_4+wH_h\|_2\|\psi\|_2\\
&\le \big(\|f_4\|_2+\|wH_h\|_2\big)\,r_h\,c_1\|v\|_2.
\end{aligned}
\end{equation}
Combining the last two inequalities and making use of \eqref{krh} we obtain
\begin{equation}\label{nabvau}
\begin{aligned}
\|\nabla^g v\|_2&\le c_1\, r_h\big(1+|h|_{2,\infty}\big)\big(\|f_4\|_2+\|wH_h\|_2\big)\\
&\le c_3\,r_h\big(1+|h|_{2,\infty}^3\big)\big(\|f_4\|_2 + \sum_i |\lambda_1^i|+|\lambda_2^i|\big).
\end{aligned}
\end{equation}
%\begin{equation*}
%\begin{aligned}
%\|\nabla^g v\|_2\le \bar r_h\big((1+|h|_{2,\infty}^3)|\Gamma|^{1/4}\|f_4\|_2 + (1+|h|_{2,\infty}^4)|\Gamma|^{1/2}\|f_0\|_2\big).
%\end{aligned}
%\end{equation*}
Now, we extend the function $v+w\,\nu_h$ from $\Gamma_h$ to $\Omega$. Let $\chi:\R\rightarrow [0,\infty)$ be a cutoff function being $1$ in $(-\kappa/2,\kappa/2)$ and $0$ outside of $(-3\kappa/4,3\kappa/4)$. Furthermore, let $\phi_0$ be $0$ in $\Omega\setminus S_\kappa$, and for $x\in S_{\kappa}$ we define 
\[\phi_0(x):=(v+w\,\nu_h)(\phi_h(\tau(x)))\,\chi(s(x)).\]
Then, from \eqref{nabvau}, \eqref{vau2}, and \eqref{krh} we infer that $\varphi_0\in H^1_0(\Omega)$ and
\begin{equation*}
\begin{aligned}
\|\nabla\phi_0\|_2&\le r_h\big(\kappa^{\frac12}(\|\nabla^g v\|_2 + \|\nabla w\|_2 +\|w\, k_h\|_2) + \kappa^{-\frac12}(\|v\|_2+\|w\|_2)\big)\\
&\le c_{\frac72} r_h\big(1+|h|_{2,\infty}^3+|h|_{3,2}^3\big)\big(\|f_4\|_2+\sum_i|\lambda_1^i|+|\lambda_2^i|\big).
\end{aligned}
\end{equation*}
Let $\phi_1\in H^1_0(\Omega\setminus\Gamma_h)$ solve \eqref{divphi1} in $\Omega\setminus\Gamma_h$ and define $\phi:=\phi_0+\phi_1$. From the last estimate and Lemma \ref{lemma:divh} we finally obtain
\begin{equation}\label{nabphi}
\begin{aligned}
\|\nabla\phi\|_2\le r_h\|f_2\|_2 + c_{\frac72}r_h\big(1+|h|_{2,\infty}^3+|h|_{3,2}^3\big)\big(\|f_4\|_2+\sum_i|\lambda_1^i|+|\lambda_2^i|\big).
\end{aligned}
\end{equation}
The fact that this estimate and \eqref{nabvau} involve the numbers $\lambda^i_j$ is very natural; situations in which one of the $\Gamma^i_h$ tends to a sphere without $f_2$, $f_4$ fulfilling (asymptotically) the compatibility condition \eqref{comp1} have to be penalized. Unfortunately, taking the dual of these inequalities doesn't lead to a suitable estimate for $\underline\nabla$, at least not without further work. For this reason, we simplify the above estimates using the fact that
\begin{equation*}
\begin{aligned}
&\sum_i|\lambda_1^i|+|\lambda_2^i|\\
&\le r_h\big(1+|h|_{2,\infty}^2\big)\Big(1+\sum_{i}\big(|\Gamma^i_h|(H_h,H_h)_{\Gamma_h^i}-(1,H_h)^2_{\Gamma^i_h}\big)^{-1}\Big)\big(c_{\frac72}\|f_2\|_2+c_4\|f_4\|_2\big),
\end{aligned}
\end{equation*}
where the sum on the right hand side is taken over all $i$ such that $\Gamma^i_h$ is not a round sphere. Combining this estimate with \eqref{nabvau} and \eqref{nabphi} we obtain
\begin{equation*}
\begin{aligned}
&\|\phi\|_{Y_h}\\
&\le r_h\big(1+|h|_{2,\infty}^5+|h|_{3,2}^5\big)\Big(1+\sum_{i}\big(|\Gamma^i_h|(H_h,H_h)_{\Gamma_h^i}-(1,H_h)^2_{\Gamma^i_h}\big)^{-1}\Big)c_7\|(f_2,f_4)\|_{Z_h}.
\end{aligned}
\end{equation*}

\medskip
\noindent (c) Next, we have to take a closer look at the constants occuring in the proof of Theorem \eqref{thm:s1}. Using Lemma \ref{lemma:trace} and \eqref{krh} we see that the constant $c_{-\frac12}$ in \eqref{bf} can be bounded by
\[r_h(1+|h|_{2,\infty})\,c_\frac12.\]
From this fact and steps (a) and (b) we obtain
 \begin{equation}\label{u1v1}
 \begin{aligned}
&\mu_b^\frac12 \|\nabla u_1\|_{2} + \mu^\frac12 \|\nabla^g v_1\|_{2} \\
&\quad\le r_h\big(1+|h|_{2,\infty}^8+|h|_{3,2}^8\big)\Big(1+\sum_{i}\big(|\Gamma^i_h|(H_h,H_h)_{\Gamma_h^i}-(1,H_h)^2_{\Gamma^i_h}\big)^{-1}\Big)\,c_{10}\\
&\quad\quad\times\Big(\frac{1}{\mu_b^\frac12}\|f_1\|_{-1,2}+\mu_b^\frac12\|f_2\|_{2}+\frac{1}{\mu_b^\frac12}\|f_3\|_{-\frac12,2}+\mu^\frac12\|f_4\|_{2}\Big).
\end{aligned}
\end{equation}
Using Lemma \ref{lemma:trace} and \eqref{krh} we see that the constant $c_{-\frac12}$ in \eqref{bmf} can be bounded by
\[r_h(1+|h|_{2,\infty})\,c_\frac12\]
while the constant $c_{-1}$ can be bounded by
\[r_h(1+|h|_{2,\infty}^2)\,c_1\]
and the constant $c$ is independent of $h$. From these facts, \eqref{u1v1}, and steps (a) and (b) we obtain
\begin{equation*}
 \begin{aligned}
&|B(u_1,\phi)-F_{u_0}(\phi)|\\
&\quad\le r_h\big(1+|h|_{2,\infty}^{10}+|h|_{3,2}^{10}\big)\Big(1+\sum_{i}\big(|\Gamma^i_h|(H_h,H_h)_{\Gamma_h^i}-(1,H_h)^2_{\Gamma^i_h}\big)^{-1}\Big)\,c_{12}\\
&\quad\quad\times\Big(\frac{1}{\mu_b^\frac12}\|f_1\|_{-1,2}+\mu_b^\frac12\|f_2\|_{2}+\frac{1}{\mu_b^\frac12}\|f_3\|_{-\frac12,2}+\mu^\frac12\|f_4\|_{2}\Big).
\end{aligned}
\end{equation*}

\medskip
\noindent (d) Finally, we combine the above estimates to obtain the stated inequality for $u=u_1-u_0$ and $(\pi,q)=\underline\nabla^{-1}\big(B(u_1,\cdot)-F\big)$.
\qed

\medskip
As already pointed out in the preceding proof, the occurence of the possibly small denominators in \eqref{rhgamma} is related to the fact that the whole expression contains no control of the compatibility condition \eqref{comp1} that $f_2$, $f_4$ must fulfill asymptotically as one of the $\Gamma^i_h$ tends to a round sphere. If one wants to do analysis in the neighbourhood of a round sphere it might become necessary to prove a refined version of \eqref{rhgamma} that resolves the singular structure in a more precise way. Furthermore, note that \eqref{rhgamma} is not optimal with respect to the regularity requirements on $h$ and with respect to the exponents of the semi-norms of $h$.

Let us proceed with the analysis of the constant in \eqref{est1nu}.
\begin{theorem}
Let $h\in W^{2,\infty}(\Gamma)$ such that $\|h\|_{L^\infty(\Gamma)}<\kappa/2$. For all $(f_1,\ldots,f_4)\in\F^1_\nu(\Gamma_h)$ there exists a unique weak solution $(u,\pi,q)\in\E^0_\nu(\Gamma_h)$ of \eqref{s2}. The estimate \eqref{est1nu} holds with the constant $c$ being bounded by
\begin{equation*}
\begin{aligned}
r_h\big(1+|h|_{2,\infty}^{8}\big)\,c_{8}.
\end{aligned}
\end{equation*}
\end{theorem}
\proof We can proceed very similarly to the proof of Theorem \ref{thm:fh}. Here, however, there is no singular behavior in the vicinity of round spheres since we don't need to prescribe additional compatibility conditions on these. Furthermore, the construction of $u_0$ is much less involved than in the proof of Theorem \ref{thm:fh}, resulting in lower regularity requirements on $h$ and in a lower exponent.
\qed

\medskip
Next, we analyze the constants in \eqref{est:regs1}.

\begin{theorem}\label{s1h}
 Let $s\ge 1$ and $h\in H^{s+2}(\Gamma)\cap W^{2,\infty}(\Gamma)$ such that $\|h\|_{L^\infty(\Gamma)}<\kappa/2$. For all $(f_1,\ldots,f_4)\in\F^s(\Gamma_h)$ there exists a unique solution $(u,\pi,q)\in\E^s(\Gamma_h)$ of \eqref{s1}. The estimate \eqref{est:regs1} holds with the constant $c$ being bounded by $r_h$ and the constant $c_{-s}$ being bounded by
 \[r_h\big(1+c_{\frac12}|h|_{2,\infty}^s|h|_{s+2,2}^{\frac{1}{2s}}+c_{s+\frac32}|h|_{s+2,2}^{2+\frac{3}{2s}}\big);\]
 note that for $s>1$ we can get rid of the term involving $|h|_{2,\infty}^s$ by interpolation.
 \end{theorem}
\proof It suffices to prove the estimates for smooth $h$ and smooth $(u,\pi,q)$. Indeed, then, approximating $h$ by smooth  functions, applying Theorem \ref{thm:regs1} and Theorem \ref{thm:fh}, and taking the limit, the claim follows. 

As in the proof of Theorem \ref{thm:regs1} we write $\tilde u:=\Phi_h^*u$, $\tilde\pi:=\Phi_h^*\pi$, $\tilde v:=\Phi_h^*v$, $\tilde w:=\Phi_h^*w$, $\tilde q:=\Phi_h^*q$, $\tilde f_i=\Phi_h^*f_i$ for $i=1,2,4$, and $\tilde f_3:=\tilde f_3^\top+\tilde f_3^\perp\, e_3:=\Phi_h^*(P_\Gamma f_3)+\Phi_h^*(f_3\cdot\nu)\,e_3$. As before, by naturality of covariant differentiation under isometries, we obtain
\begin{equation}\label{eqn:tildesystemh}
\begin{aligned}
\mu_b\Delta_{\tilde e}\tilde u-\grad_{\tilde e}\tilde\pi&=\tilde f_1&&\mbox{ in }\Omega\setminus\Gamma,\\
\Div_{\tilde e} \tilde u&=\tilde f_2&&\mbox{ in }\Omega\setminus\Gamma,\\
\mu\big(\Delta_{\tilde g} \tilde v + \grad_{\tilde g}(\tilde w\,H_{\tilde e}) + K_{\tilde g}\,\tilde v -2\Div_{\tilde g}(\tilde w\,k_{\tilde e})\big)&\\
-\grad_{\tilde g} \tilde q  
+2\mu_b[\![D^{\tilde e}\tilde u]\!]\nu_{\tilde e}&=\tilde f_3^\top&&\mbox{ on } \Gamma,\\
2\mu\big(\langle\nabla^{\tilde g} \tilde v,k_{\tilde e}\rangle_{\tilde g}
-\tilde w\,(H^2_{\tilde e}-2K_{\tilde g})\big)-\tilde q\,H_{\tilde e}-[\![\tilde\pi]\!]&=\tilde f_3^\perp&&\mbox{ on }\Gamma,\\
\Div_{\tilde g} \tilde v-\tilde w\, H_{\tilde e}&=\tilde f_4&&\mbox{ on }\Gamma,\\
\tilde u-\tilde v-\tilde w\,\nu_{\tilde e}&=\tilde f_5&&\mbox{ on }\Gamma
\end{aligned}
\end{equation}
with $\tilde f_5=0$.

\medskip
\noindent(a) To begin with, let us assume that $\tilde e$ is constant in $\Omega$ and that $\tilde f_5\in H^{s+1/2}$ is non-vanishing. In this case we can proceed by localization and transformation as in the proof of Theorem \ref{thm:regs1} to prove the estimates for integer $s$ and then interpolate for the general case. The only difference is that one has to consider a non-Euclidean, but constant Riemannian metric in the prototype geometries. In Cartesian coordinates the Laplace operators take the form
\[\Delta_{\tilde e}\tilde u=\tilde e^{ij}\, \pa_{ij}^2 \tilde u\quad\text{in the bulk,}\quad
\Delta_{\tilde g}\tilde v=\tilde g^{\alpha\beta}\, \pa_{\alpha\beta}^2 \tilde v\quad\text{on $Q^2_R$.}\]
In this case, the constants in \eqref{est:regs1} are polynomial functions of the ratio of the largest and the smallest eigenvalue of $\tilde e$.

\medskip
\noindent(b) In the next step, let $\hat e^{ij}:=\tilde e^{ij}(x_0)$ for some fixed $x_0\in\Omega$ and assume $\|\tilde e^{ij}-\hat e^{ij}\|_{\infty}$ to be sufficiently small. We need to write the right hand sides in \eqref{eqn:tildesystemh} in a more precise form. The results from \Appref{app:cov} show that \eqref{eqn:tildesystemh} can be written as
\begin{equation}\label{eqn:tildesystem2}
\begin{aligned}
\mu_b\Delta_{\hat e}\tilde u-\grad_{\hat e}\tilde\pi&=\hat f_1&&\mbox{ in }\Omega\setminus\Gamma,\\
\Div_{\hat e} \tilde u&=\hat f_2&&\mbox{ in }\Omega\setminus\Gamma,\\
\mu\big(\Delta_{\hat g} \tilde v + \grad_{\hat g}(\tilde w\,H_{\hat e}) + K_{\hat g}\,\tilde v -2\Div_{\hat g}(\tilde w\,k_{\hat e})\big)&\\
-\grad_{\hat g} \tilde q  
+2\mu_b[\![D^{\hat e}\tilde u]\!]\nu_{\hat e}&=\hat f_3^\top&&\mbox{ on } \Gamma,\\
2\mu\big(\langle\nabla^{\hat g} \tilde v,k_{\hat e}\rangle_{\hat g}
-\tilde w\,(H^2_{\hat e}-2K_{\hat g})\big)-\tilde q\,H_{\hat e}-[\![\tilde\pi]\!]&=\hat f_3^\perp&&\mbox{ on }\Gamma,\\
\Div_{\hat g} \tilde v-\tilde w\, H_{\hat e}&=\hat f_4&&\mbox{ on }\Gamma,\\
\tilde u-\tilde v-\tilde w\,\nu_{\hat e}&=\hat f_5&&\mbox{ on }\Gamma
\end{aligned}
\end{equation}
with
\begin{equation}\label{eqn:hutfunkt}
\begin{aligned}
\hat f_1&=\tilde f_1 + (\tilde e-\hat e)*r(\hat e)*(\mu_b\nabla^2 \tilde u,\grad\tilde\pi)\\
&\quad + \mu_b\, r(\tilde e,\hat e)*\big((\nabla^2\tilde e,(\nabla\tilde e)^2)*\tilde u+\nabla\tilde e*\nabla \tilde u\big),\\
\hat f_2&=\tilde f_2+r(\tilde e,\hat e)*\nabla\tilde e*\tilde u,\\
\hat f_3^\top
% &=\tilde f_3^\top+(\tilde g-g)*(\mu(\nabla^g)^2 \tilde v,\grad_g \tilde q) + \mu\,r(\tilde g,\hat g)*\big((\nabla^{g})^2\tilde g*\tilde v,(\nabla^{g}\tilde g)^2*\tilde v,\nabla^g\tilde g*\nabla^g \tilde v\big)\\
% &\quad + (w,v)*r(\tilde e,\hat e)*(k^2,\nabla k,k*\nabla\tilde e,\nabla^2\tilde e,(\nabla\tilde e)^2) + \nabla w*r(\tilde e,\hat e)*(k,\nabla\tilde e)\big)\\
% &\quad + \mu_b\,r(\tilde e,\hat e)*([\nabla \tilde u],\nabla\tilde e*[\tilde u])\\
&=\tilde f_3^\top+(\tilde e-\hat e)*r(\hat e)*(\mu(\nabla^g)^2 \tilde v,\grad_g \tilde q) + \mu_b\,r(\tilde e,\hat e)*\big([\nabla \tilde u]+\nabla\tilde e*[\tilde u]\big)\\
&\quad+ \mu\,r(\tilde e,\hat e)*\big((k^2,\nabla k,k*\nabla\tilde e,\nabla^2\tilde e,(\nabla\tilde e)^2)*[\tilde u] + (k,\nabla\tilde e)*[\nabla\tilde u]\big),\\
\hat f_3^\perp&=\tilde f_3^\perp + \mu\,r(\tilde e,\hat e)*\big((k,\nabla\tilde e)*\nabla^g\tilde v + (k^2,k*\nabla\tilde e,(\nabla\tilde e)^2)*[\tilde u]\big)\\
% &=\tilde f_3^\perp + \mu\,r(\tilde e,\hat e)*\big((k,\nabla\tilde e)*\nabla^g\tilde v + \nabla \tilde e*(k,\nabla\tilde e)*\tilde v + w\,(k^2,k*\nabla\tilde e,(\nabla\tilde e)^2)\big)\\
&\quad+ r(\tilde e,\hat e)*(k,\nabla\tilde e)\,\tilde q,\\
\hat f_4&=\tilde f_4 + r(\tilde e,\hat e)*(k,\nabla\tilde e)*[\tilde u],\\
\hat f_5&=(\tilde e-\hat e)*r(\tilde e,\hat e)\,\tilde w.
\end{aligned}
\end{equation}
Here, we replaced $\nabla^{\hat e}$ and $\nabla^{\hat g}$ by $\nabla=\nabla^{e}$ and $\nabla^{g}$ by absorbing the additional coefficients into $r(\hat e)$ and $r(\tilde e,\hat e)$, respectively. It is not hard to see that in $S_\kappa$
\[\tilde e-e=r(h/\kappa,hk,\nabla h).\]
Here and in the following lines, for simplicity of notation, we write $h$ and $k$ instead of $h\circ\tau$ and $k\circ\tau$. From this identity, by induction we infer that for $l\in\N$ in $S_\kappa$ we have
\begin{equation*}
\begin{aligned}
 \nabla^l\tilde e=\sum_{p,p_i,p_i'\ge 0}(\nabla^{p_1} r)(h/\kappa,hk,\nabla h)\,h^{p_0}*\prod_{i=2}^{l+1}((\nabla^g)^i h)^{p_i}*\prod_{i=0}^l((\nabla^g)^ik)^{p_i'}\ \kappa^{-p},
\end{aligned}
\end{equation*}
where $\nabla^{p_1}r$ stands for derivatives of $r$ with respect to the spatial variable and
\begin{equation*}
-p_0+p_1+\sum_{i=2}^{l+1}(i-1)p_i+\sum_{i=0}^{l}(i+1)p_i'+p=l, \quad p_0\le \sum_{i=0}^{l}(i+1)p_i'.
\end{equation*}
In particular, we have $\sum_{i=2}^{l+1}(i-1)p_i\le l$. From this we easily deduce that in $S_\kappa$
\begin{equation}\label{nablare}
\begin{aligned}
 \nabla^lr(\tilde e,\hat e)&=\sum_{q_j\ge 0}r(\tilde e,\hat e)*\prod_{j=1}^{l}(\nabla^j \tilde e)^{q_j}\\
 &=\sum_{p_i\ge 0}r(h/\kappa,hk,\nabla h,\hat e)*\prod_{i=2}^{l+1}((\nabla^g)^i h)^{p_i},
\end{aligned}
\end{equation}
where $\sum_{j=1}^ljq_j=l$ and $\sum_{i=2}^{l+1}(i-1)p_i\le l$. Now, fix some $\delta\in[0,1)$. For $\delta>0$ and $a\in\{1+\delta,2,2+\delta,\ldots\}$ let
\[I_a:=\{1+\delta,2,2+\delta,\ldots,a\},\]
while for $\delta=0$ and $a\in\N_{\ge 2}$ let
\[I_a:=\{2,3,\ldots,a\}.\]
Finally, for $a$ as above and $b,c\in\{0,\delta,1,1+\delta,\ldots\}$ let
\begin{equation*}
\begin{aligned}
J_a^{b,c}:=\Big\{&(p_i,l)_{i\in I_a}\,|\,l\in \{0,\delta,1,1+\delta,\ldots,c\}\ \text{and}\\
&p_i\in\N_{\ge 0}\ \text{such that}\ \sum_{i\in I_a}(i-1)p_i+l\le b\Big\}.
\end{aligned}
\end{equation*}

\medskip
\noindent(b.1) We start with the analysis of the term $\hat f_3^\perp$ which is, with regard to the exponent, the limiting one. We write $s+\frac12=n+\delta$ for $n\in\N_{\ge 1}$ and $\delta\in [0,1)$. Starting from the expression for $\hat f_3^\perp$ in \eqref{eqn:hutfunkt} and using \eqref{nablare}, Lemma \ref{lemma:prod}, Corollary \ref{cor:analytic}, and \eqref{qeinbettung} we infer that
\begin{equation*}
\begin{aligned}
\frac{1}{\mu_b^{\frac12}}|\hat f_3^\perp|_{s-\frac12,2}&\le \frac{1}{\mu_b^{1/2}}|\tilde f_3^\perp|_{s-\frac12,2} + \frac{\mu}{\mu_b^{\frac12}}\sum_{(p_i,l)\in J_{s+3/2}^{s+3/2,s+1/2}}r_h \prod_{i\in I_{s+3/2}}|h|_{i,p_ir_i,2}^{p_i}|\tilde v|_{l,r,2}\\
&\quad+ \frac{\mu}{\mu_b^{\frac12}}\sum_{(p_i,l)\in J_{s+3/2}^{s+3/2,s-1/2}}r_h \prod_{i\in I_{s+3/2}}|h|_{i,p_ir_i,2}^{p_i}|[\tilde u]|_{l,r,2}\\
&\quad+ \frac{1}{\mu_b^{\frac12}}\sum_{(p_i,l)\in J_{s+3/2}^{s+1/2,s-1/2}}r_h \prod_{i\in I_{s+3/2}}|h|_{i,p_ir_i,2}^{p_i}|\tilde q|_{l,r,2},
\end{aligned}
\end{equation*}
where $2\le r,r_i\le\infty$ are chosen such that 
\[\frac{1}{r}+\sum_{i:\,p_i\not=0}\frac{1}{p_ir_i}=\frac12.\]
Here and in the following, with a slight abuse of notation we define $|\cdot|_{0,r,2}:=\|\cdot\|_{r}$ for ease of notation.

To begin with, let us consider the term involving $|\tilde v|_{s+1/2,r,2}$ and $p_2=1$; in this case all other $p_i$ vanish since $\sum (i-1)p_i\le 1$. Assuming $2\le r<4$, we use \eqref{inth} with $t=2$, $p=r_2=2r/(r-2)$, and $\alpha=2/(sr)$ as well as \eqref{intv} with $t=s+1/2$, $p=r$, and $\beta=(s+1-2/r)/(s+1/2)$; for the first estimate we have to exclude the case $r=2$ if $s=1$. Hence, we have
\[\frac{\mu}{\mu_b^{1/2}}|h|_{2,r_2,2}|\tilde v|_{s+\frac12,r,2}\le \frac{r_h}{\eps^{\frac{\beta}{1-\beta}}}\Big(1+\Big(\frac{\mu}{\mu_b}\Big)^{1-\frac{\beta}{2}}|h|_{s+2,2}^\alpha\Big)^{\frac{1}{1-\beta}}\mu_b^{\frac12}\|\nabla\tilde u\|_{2}+ \eps\mu^{\frac12}|\tilde v|_{s+1,2}.\]
The exponent
\[\frac{\alpha}{1-\beta}=\Big(2+\frac1s\Big)\frac{1}{2-r/2}\]
of $|h|_{s+2,2}$ is minimized for $r=2$. Furthermore, $\eps$ has to be chosen so small that the $|\tilde v|_{s+1,2}$-term can be absorbed on the left hand side of the estimate proved in step (a). However, the constants of continuity in the latter estimate are of the form $r_h$; hence, $\eps$ can be chosen to have the form $r_h^{-1}$. Thus, for $s>1$ we have
\begin{equation}\label{prototype}
 \frac{\mu}{\mu_b^{1/2}}|h|_{2,\infty}|\tilde v|_{s+\frac12,2}\le r_h\Big(1+\Big(\frac{\mu}{\mu_b}\Big)^{s+1}|h|_{s+2,2}^{2+\frac1s}\Big)\mu_b^{\frac12}\|\nabla\tilde u\|_{2}+ \eps\mu^{\frac12}|\tilde v|_{s+1,2}.
\end{equation}
For $s=1$ we use the interpolation estimate only for $\tilde v$. This gives
\[\frac{\mu}{\mu_b^{1/2}}|h|_{2,\infty}|\tilde v|_{\frac32,2}\le r_h\Big(1+\Big(\frac{\mu}{\mu_b}\Big)^2|h|_{2,\infty}^3\Big)\mu_b^{\frac12}\|\nabla\tilde u\|_{2}+ \eps\mu^{\frac12}|\tilde v|_{s+1,2}.\]

Next, let us consider the remaining terms involving $|\tilde v|_{s+1/2,r,2}$. Note that in this case all $p_i$ with $i\ge 2$ vanish. Again, we choose $r=2$ and $r_i=\infty$ for all $i$ and consider \eqref{inth} with $t=i$, $p=\infty$, and $\alpha=\theta_i:=(i-1)/s$; note that $\theta_i\in (0,1)$. From this estimate and \eqref{intv} with $t=s+1/2$, $p=2$, and $\beta=s/(s+1/2)$ we obtain
\begin{equation*}
\begin{aligned}
\frac{\mu}{\mu_b^{1/2}}\prod_{i\in I_{s+3/2}}|h|_{i,\infty,*}^{p_i}|\tilde v|_{s+\frac12,2}&\le r_h\Big(1+\Big(\frac{\mu}{\mu_b}\Big)^{1-\frac{\beta}{2}}|h|_{s+2,2}^{\sum\theta_ip_i}\Big)^{\frac{1}{1-\beta}}\mu_b^{\frac12}\|\nabla\tilde u\|_{2}\\
&\quad+ \eps\mu^{\frac12}|\tilde v|_{s+1,2}.
\end{aligned}
\end{equation*}
Since $\sum\theta_ip_i\le 1$ this gives
\begin{equation*}
\begin{aligned}
\frac{\mu}{\mu_b^{1/2}}\prod_{i\in I_{s+3/2}}|h|_{i,\infty,2}^{p_i}|\tilde v|_{s+\frac12,2}&\le r_h\Big(1+\Big(\frac{\mu}{\mu_b}\Big)^{s+1}|h|_{s+2,2}^{2+\frac1s}\Big)\mu_b^{\frac12}\|\nabla\tilde u\|_{2}+ \eps\mu^{\frac12}|\tilde v|_{s+1,2}.
\end{aligned}
\end{equation*}

Now, we assume that $\delta>0$ and consider the terms involving $|\tilde v|_{n,r,2}$. In this case we have $\sum (i-1)p_i\le 1+\delta$ and hence all $p_i$ with $i\ge 3$ vanish. For $s>3/2$, we can choose $r=2$ and $r_i=\infty$ for all $i$ since then we have $\theta_i\in (0,1)$ even for $i=2+\delta$. Using \eqref{inth} with $t=i$ and $\alpha=\theta_i$ as well as \eqref{intv} with $t=n$ and $\beta=(n-1/2)/(s+1/2)$ we obtain
\begin{equation}\label{vn}
\begin{aligned}
\frac{\mu}{\mu_b^{1/2}}\prod_{i\in I_{s+3/2}}|h|_{i,\infty,2}^{p_i}|\tilde v|_{n,2}&\le r_h\Big(1+\Big(\frac{\mu}{\mu_b}\Big)^{1-\frac{\beta}{2}}|h|_{s+2,2}^{\sum\theta_ip_i}\Big)^{\frac{1}{1-\beta}}\mu_b^{\frac12}\|\nabla\tilde u\|_{2}\\
&\quad+ \eps\mu^{\frac12}|\tilde v|_{s+1,2}.
\end{aligned}
\end{equation}
Since $\sum\theta_ip_i\le (1+\delta)/s$ and 
\[\frac{1+\delta}{s}\frac{1}{1-\beta}=\Big(2+\frac1s\Big)\frac{1+\delta}{1+2\delta}< 2+\frac1s\]
we can dominate the right hand side of \eqref{vn} by the right hand side of \eqref{prototype}. However, in the case $s \le 3/2$ we need to proceed differently. In this case we have $n=1$ and $\delta\in [1/2,1)$. We choose $r=4$, $r_j=4$, and $r_i=\infty$ for $i\not=j$; here and in the following we denote by $j$ the largest $i\in I_{s+3/2}$ such that $p_i\not=0$. We want to use \eqref{inth} with $t=j$, $p=4p_j$, and $\alpha=\theta_j':=(j-1-1/(2p_j))/s$. Note that 
\[0\le\frac{\delta-1/(2p_j)}{s}\le\theta_j'\le \frac{1/2+\delta}{s}=1;\]
in the case $\delta=1/2$ we may assume without restriction that $p_j>1$ so that $\theta_j'>0$. Using this estimate as well as \eqref{inth} with $t=i$, $p=\infty$, and $\alpha=\theta_i$ and \eqref{intv} with $t=n=1$, $p=4$, and $\beta=1/(s+1/2)$ we obtain
\begin{equation}\label{vn4}
\begin{aligned}
&\frac{\mu}{\mu_b^{1/2}}\prod_{i\in I_{s+3/2}\atop i\not=j}|h|_{i,\infty,2}^{p_i}|h|_{i,4p_j,2}^{p_j}|\tilde v|_{1,4}\\
&\quad\le r_h\Big(1+\Big(\frac{\mu}{\mu_b}\Big)^{1-\frac{\beta}{2}}|h|_{s+2,2}^{\sum\theta_ip_i+\theta_j'p_j}\Big)^{\frac{1}{1-\beta}}\mu_b^{\frac12}\|\nabla\tilde u\|_{2} + \eps\mu^{\frac12}|\tilde v|_{s+1,2}.
\end{aligned}
\end{equation}
Since $\sum\theta_ip_i+\theta_j'p_j\le (1/2+\delta)/s$ and 
\[\frac{1/2+\delta}{s}\frac{1}{1-\beta}=1+\frac{1}{s-1/2}\le 2+\frac1s\]
for $s\ge 1$ we can dominate the right hand side of \eqref{vn4} by the right hand side of \eqref{prototype}.

Now, we consider the terms involving $|[\tilde u]|_{s-\frac12,r,2}$. In this case we have $\sum (i-1)p_i\le 2$; however, since $\hat f_3^\perp$ only contains up to first derivatives (squared) of $\tilde e$ and no second derivatives, in fact all $p_i$ with $i>2$ vanish. We choose $r=2$ and $r_i=\infty$ for all $i$ and use \eqref{inth} with $t=i$, $p=\infty$, and $\alpha=\theta_i$ as well as \eqref{intu} with $t=s-1/2$, $p=2$, and $\gamma=(s-1)/s$; again, we have to exclude the case $s=1$ so that $\theta_2=1/s\in (0,1)$. This gives
\begin{equation*}
\begin{aligned}
\frac{\mu}{\mu_b^{\frac12}}r_h \prod_{i\in I_{s+3/2}}|h|_{i,\infty,2}^{p_i}|[\tilde u]|_{s-\frac12,2,2}&\le r_h\Big(1+\frac{\mu}{\mu_b}|h|_{s+2,2}^{\sum\theta_ip_i}\Big)^{\frac{1}{1-\gamma}}\mu_b^{\frac12}\|\nabla\tilde u\|_{2}\\
&\quad+ \eps\mu_b^{\frac12}|\tilde u|_{s+1,2}.
\end{aligned}
\end{equation*}
Since $\sum\theta_ip_i\le 2/s$ and $1/(1-\gamma)=s$ we can infer that 
\begin{equation}\label{prototypeu}
\begin{aligned}
\frac{\mu}{\mu_b^{\frac12}}r_h \prod_{i\in I_{s+3/2}}|h|_{i,\infty,2}^{p_i}|[\tilde u]|_{s-\frac12,2,2}&\le r_h\Big(1+\Big(\frac{\mu}{\mu_b}\Big)^s|h|_{s+2,2}^{2}\Big)\mu_b^{\frac12}\|\nabla\tilde u\|_{2}\\
&\quad+ \eps\mu_b^{\frac12}|\tilde u|_{s+1,2}.
\end{aligned}
\end{equation}
For $s=1$ we simply don't use the estimate \eqref{inth} for $t=2$ and obtain
\begin{equation*}
\begin{aligned}
\frac{\mu}{\mu_b^{\frac12}}r_h \prod_{i\in I_{5/2}}|h|_{i,\infty,2}^{p_i}|[\tilde u]|_{\frac12,2,2}&\le r_h\Big(1+\frac{\mu}{\mu_b}|h|_{3,2}^{\sum_{i\not=2}\theta_ip_i}|h|_{2,\infty}^{p_2}\Big)\mu_b^{\frac12}\|\nabla\tilde u\|_{2}\\
&\quad+ \eps\mu_b^{\frac12}|\tilde u|_{2,2}.
\end{aligned}
\end{equation*}
An application of Young's inequality yields
\begin{equation}\label{prototypeus1}
\begin{aligned}
\frac{\mu}{\mu_b^{\frac12}}r_h \prod_{i\in I_{5/2}}|h|_{i,\infty,2}^{p_i}|[\tilde u]|_{\frac12,2,2}&\le r_h\Big(1+\frac{\mu}{\mu_b}\big(|h|_{3,2}^2+|h|_{2,\infty}^{2}\big)\Big)\mu_b^{\frac12}\|\nabla\tilde u\|_{2}\\
&\quad+ \eps\mu_b^{\frac12}|\tilde u|_{2,2}.
\end{aligned}
\end{equation}

Next, we assume that $\delta>0$ and consider the terms involving $|[\tilde u]|_{n-1,r,2}$. In this case we have $\sum (i-1)p_i\le 2+\delta$; however, all $p_i$ with $i\ge 3$ vanish. As before, we have to distinguish the cases $s>3/2$ and $s\le 3/2$. In the case $s>3/2$ we choose $r=2$ and $r_i=\infty$; note that $\theta_i\in (0,1)$ even for $i=2+\delta$. We apply \eqref{inth} with $t=i$ and $\alpha=\theta_i$ and \eqref{intu} with $t=n-1$ and $\gamma=(n-3/2)/s$ and proceed as in the case $l=s-1/2$. Since $\sum\theta_ip_i\le (2+\delta)/s$ and
\[\frac{2+\delta}{s}\frac{1}{1-\beta}=\frac{2+\delta}{1+\delta}\le 2\]
we can dominate our terms by the right hand sides of \eqref{prototypeu} and \eqref{prototypeus1}, respectively. In the case $s\le 3/2$ we note that $\delta\ge 1/2$ and choose $r=4$, $r_j=4$, and $r_i=\infty$ for all $i\not=j$. Note that $i\le 2$ and hence $\theta_i\in (0,1)$, except for $\theta_2$ in the case $s=1$; this case has to be treated separately, exactly as above. Furthermore, note that
\[\frac{\delta-1/(2p_j)}{s}\le \theta_j'\le \frac{1/2+\delta}{s}\le 1;\]
without loss we may assume $p_j$ to be so large that $\theta_j'>0$. Now, we apply \eqref{inth} with $t=i$ and $\alpha=\theta_i$, \eqref{inth} with $t=j$ and $\alpha=\theta_j'$, and \eqref{intu} with $t=n-1$ and $\gamma=(n-1)/s$ and proceed as in the case $l=s-1/2$. Since $\sum\theta_ip_i+\theta_j'p_j\le (3/2+\delta)/s$ and
\[\frac{3/2+\delta}{s}\frac{1}{1-\beta}=\frac{3/2+\delta}{1/2+\delta}\le 2\]
we can dominate our terms by the right hand sides of \eqref{prototypeu} and \eqref{prototypeus1}, respectively.

Next, we consider the terms involving $|[\tilde u]|_{l,r,2}$ with $0\le l\le s-3/2$. In this case we have $\sum (i-1)p_i\le s+3/2-l$; however, all $p_i$ with $i> s+3/2$ vanish. We choose $r=4$, $r_j=4$, and $r_i=\infty$ for all $i\not=j$. Note that for all $i\not=j$ we have $i\le n+1$ in the case $\delta>0$ and $i\le n$ in the case $\delta=0$; therefore, we have $\theta_i\in (0,1)$. Furthermore, we note that $p_{s+3/2}\le 1$ so that
\[\theta_{s+3/2}'\le \frac{s+1/2-1/(2p_j)}{s}\le 1;\]
similarly we see that $\theta_{n+1}'\le 1$, and $\theta_j'\le 1$ holds trivially for $j\le s+1/2$.  We apply \eqref{inth} with $t=i$ and $\alpha=\theta_i$, \eqref{inth} with $t=j$ and $\alpha=\theta_j'$, and $\eqref{intu}$ with $t=l$ and $\gamma=l/s$. Since $\sum\theta_ip_i+\theta_j'p_j\le (s-l+1)/s$ and
\[\frac{s-l+1}{s}\frac{1}{1-\beta}=\frac{s-l+1}{s-l}\le \frac53\]
we can dominate our terms by the right hand sides of \eqref{prototypeu} and \eqref{prototypeus1}, respectively. 

The terms involving the surface pressure $\tilde q$ can be treated almost exactly like the terms containing the trace $[\tilde u]$ of the bulk velocity; the main difference is that one has to interpolate between $|\tilde q|_{s,2}$ and $\|\tilde q\|_{2}$ instead of $|\tilde u|_{s+1,2}$ and $\|\nabla \tilde u\|_{2}$. We only present the case $l=s-1/2$ where $\sum (i-1)p_i\le 1$. We choose $r=2$ and $r_i=\infty$ for all $i$; again we have to treat the case $s=1$ separately. Applying \eqref{q2} with $t=s-1/2$ and \eqref{inth} with $t=i$ and $\alpha=\theta_i$ we obtain
\begin{equation*}
\begin{aligned}
\frac{1}{\mu_b^{\frac12}}r_h \prod_{i\in I_{s+3/2}}|h|_{i,\infty,2}^{p_i}|\tilde q|_{s-1/2,2}&\le r_h\Big(1+\Big(\frac{\mu}{\mu_b}\Big)^{s}|h|_{s+2,2}^{2}\Big)\frac{1}{\mu^{\frac12}}\|\tilde q\|_{2}\\
&\quad+ \eps\frac{1}{\mu^{\frac12}}|\tilde u|_{2,2}.
\end{aligned}
\end{equation*}

\medskip
\noindent(b.2) Next, we analyze the term $\hat f_5$. Let $s+\frac12=n+\delta$ for $n\in\N_{\ge 1}$ and $\delta\in [0,1)$. We have
\begin{equation*}
\begin{aligned}
&\mu_b^{\frac12}|\hat f_5|_{s+\frac12,2}\le \|\tilde e-\hat e\|_\infty r_h \mu_b^{\frac12}|\tilde u|_{s+1,2}\\
&\quad+ \mu_b^{\frac12}r_h\bigg(\sum_{(p_i,l)\in J_{s+3/2}^{s+1/2,s-\frac12}}\prod_{i\in I_{s+3/2}}|h|_{i,p_ir_i,2}^{p_i}|[\tilde u]|_{l,r,2}+b(\delta)|h|_{1+\delta,r_{1+\delta},2}|[\tilde u]|_{n,r,2}\bigg),
\end{aligned}
\end{equation*}
where $b(\delta)=1$ if $\delta>0$ and $b(\delta)=0$ if $\delta=0$, and where $2\le r,r_i\le\infty$ are chosen as above. 

First, let us assume that $\delta>0$ and consider the terms involving $|[\tilde u]|_{n,r}$. In this case we have $\sum(i-1)p_i\le \delta$; hence, all $p_i$ with $i\ge 2$ vanish. Choosing $r=2$ and $r_i=\infty$ for all $i$ and applying \eqref{inth} with $t=i$ and $\alpha=\theta_i$ and \eqref{intu} with $t=n$ and $\gamma=(n-1/2)/s$ we obtain
\begin{equation}\label{f5}
\begin{aligned}
\mu_b^{\frac12}r_h \prod_{i\in I_{s+3/2}}|h|_{i,p_ir_i,2}^{p_i}|[\tilde u]|_{n,2}&\le r_h\Big(1+|h|_{s+2,2}^{\sum\theta_ip_i}\Big)^{\frac{1}{1-\gamma}}\mu_b^{\frac12}\|\nabla\tilde u\|_{2}+ \eps\mu_b^{\frac12}|\tilde u|_{s+1,2}\\
&\le r_h\big(1+|h|_{s+2,2}\big)\mu_b^{\frac12}\|\nabla\tilde u\|_{2}+ \eps\mu_b^{\frac12}|\tilde u|_{s+1,2},
\end{aligned}
\end{equation}
since $\sum\theta_ip_i\le \delta/s$ and
\[\frac{\delta}{s}\frac{1}{1-\gamma}=\frac{\delta}{s}\frac{s}{s-n+1/2}=1.\]

Next, we consider the terms involving $|[\tilde u]|_{l,r,2}$ with $0\le l\le s-1/2$. In this case we have $\sum(i-1)p_i\le s+1/2-l$. We choose $r=4$, $r_j=4$, and $r_i=\infty$ for all $i\not=j$ and apply \eqref{inth} with $t=i$ and $\alpha=\theta_i$, \eqref{inth} with $t=j$ and $\alpha=\theta_j'$, and \eqref{intu} with $t=l$ and $\gamma=l/s$. We see exactly like in the analysis of the $\hat f_3^\perp$-terms involving $|[\tilde u]|_{l,r,2}$ with $0\le l\le s-3/2$ that $\theta_i\in (0,1)$ and $\theta_j'\le 1$. Since $\sum\theta_ip_i+\theta_j'p_j\le (s-l)/s$
and \[\frac{s-l}{s}\frac{1}{1-\beta}=1\]
we can dominate our terms by the right hand side of \eqref{f5}.

\medskip
\noindent(b.3) Now, we consider the terms $\hat f_3^\top$ and $\hat f_4$. Let $s=n+\delta$ with $\delta\in [0,1)$. We have
\begin{equation*}
\begin{aligned}
\frac{1}{\mu^{\frac12}}|\hat f_3^\top|&_{s-1,2}+\mu^{\frac12}|\hat f_4|_{s,2}\le \frac{1}{\mu^{\frac12}}|\tilde f_3^\top|_{s-1,2}+\mu^{\frac12}|\tilde f_4|_{s,2}\\
&\quad+\|\tilde e-\hat e\|_{\infty}r_h \big(\mu^{\frac12}|\tilde v|_{s+1,2}+\frac{1}{\mu^{\frac12}}|\tilde q|_{s,2}\big)\\
&\quad+ \mu^{\frac12}r_h\bigg(\sum_{(p_i,l)\in J_{s+2}^{s+1,s}}\prod_{i\in I_{s+2}}|h|_{i,p_ir_i,2}^{p_i}|[\tilde u]|_{l,r,2}+b(\delta)|h|_{1+\delta,r_{1+\delta},2}|\tilde v|_{n+1,r,2}\bigg)\\
&\quad+\frac{1}{\mu^{\frac12}}r_h\bigg(\sum_{(p_i,l)\in J_{s}^{s,s-1}}\prod_{i\in I_{s}}|h|_{i,p_ir_i,2}^{p_i}|\tilde q|_{l,r,2}+b(\delta)|h|_{1+\delta,r_{1+\delta},2}|\tilde q|_{n,r,2}\bigg)\\
&\quad+ \frac{\mu_b}{\mu^{\frac12}}\sum_{(p_i,l)\in J_{s+1}^{s,s}}r_h \prod_{i\in I_{s+1}}|h|_{i,p_ir_i,2}^{p_i}|[\tilde u]|_{l,r,2},
% &\quad+ \mu^{\frac12}\sum_{(p_i,l)\in J_{s+1}^{i_*}}r_h \prod_{i\in I_{s}}|h|_{i,p_ir_i,*}^{p_i}|\tilde v|_{l,r,*}\\
% &\quad+ \frac{1}{\mu^{\frac12}}\sum_{(p_i,l)\in J_{s}^{i_{**}}}r_h \prod_{i\in I_{s}}|h|_{i,p_ir_i,*}^{p_i}|\tilde q|_{l,r,*}\\
\end{aligned}
\end{equation*}
where $b(\delta)$ and $2\le r,r_i\le\infty$ are chosen as above.
%where $i_*=n+1, i_{**}=n$ if $\delta>0$ and $i_*=s, i_{**}=s-1$ if $\delta=0$.

We begin by assuming $\delta>0$ and analyze the term involving $|\tilde v|_{n+1,r,2}$. Choosing $r=2$ and $r_{1+\delta}=\infty$ and applying \eqref{inth} with $t=1+\delta$ and $\alpha=\theta_{1+\delta}$ and \eqref{intv} with $t=n+1$ and $\beta=(s+1/2-\delta)/(s+1/2)$ we obtain
\begin{equation*}
\begin{aligned}
\mu^{\frac12}r_h|h|_{1+\delta,\infty,2}|\tilde v|_{n+1,2}\le r_h\Big(1+\Big(\frac{\mu}{\mu_b}\Big)^{\frac{1}{2}}|h|_{s+2,2}^{1+\frac{1}{2s}}\Big)\mu_b^{\frac12}\|\nabla\tilde u\|_{2} + \eps\mu^{\frac12}|\tilde v|_{s+1,2},
\end{aligned}
\end{equation*}
since
\[\frac{\theta_{1+\delta}}{1-\beta}=1+\frac{1}{2s}.\]

Next, we consider the terms involving $\mu^{\frac12}|[\tilde u]|_{s,r,2}$. In this case we have $\sum (i-1)p_i\le 1$ and hence all $p_i$ with $i>2$ vanish. We choose $r=2$ and $r_i=\infty$ for all $i$; again, the case $s=1$ and $p_2=1$ has to be treated separately. For $s>1$ we apply \eqref{inth} with $t=i$ and $\alpha=\theta_i$ and \eqref{intu} with $t=s$ and $\gamma=(s-1/2)/s$ and obtain
\begin{equation}\label{us2}
\begin{aligned}
\mu^{1/2}\prod_{i\in I_{s+2}}|h|_{i,\infty,2}^{p_i}&|[\tilde u]|_{s,2,2}\\
&\le r_h\Big(1+\Big(\frac{\mu}{\mu_b}\Big)^{\frac12}|h|_{s+2,2}^{\sum\theta_ip_i}\Big)^{\frac{1}{1-\gamma}}\mu_b^{\frac12}\|\nabla\tilde u\|_{2}+ \eps\mu^{\frac12}_b|\tilde u|_{s+1,2}\\
&\le r_h\Big(1+\Big(\frac{\mu}{\mu_b}\Big)^{s}|h|_{s+2,2}^{2}\Big)\mu_b^{\frac12}\|\nabla\tilde u\|_{2}+ \eps\mu^{\frac12}_b|\tilde u|_{s+1,2}.
\end{aligned}
\end{equation}
Here, the last inequality follows from $\sum\theta_ip_i\le 1/s$ and 
\[\frac1s\frac{1}{1-\gamma}=2.\]
We omit the case $s=1$.

Now, we assume $\delta>0$ and consider the terms involving $\mu^{\frac12}|[\tilde u]|_{n,r}$. In this case we have $\sum(i-1)p_i\le 1+\delta$ and hence all $p_i$ with $i\ge 3$ vanish. For $s\ge 2$ or $p_{2+\delta}=0$ we can again choose $r=2$ and $r_i=\infty$ for all $i$ since then we have $\theta_i\in (0,1)$. We use \eqref{inth} with $t=i$ and $\alpha=\theta_i$ and \eqref{intu} with $t=n$ and $\gamma=(n-1/2)/s$. Since $\sum\theta_ip_i\le (1+\delta)/s$ and 
\[\frac{1+\delta}{s}\frac{1}{1-\beta}=\frac{1+\delta}{1/2+\delta}\le 2\]
we can dominate our terms by the right hand side of \eqref{us2}. For $s<2$ and $p_{2+\delta}=1$ we choose $r=\min(2/(1-\delta),3)$ and $r_{2+\delta}=2r/(r-2)$. Using \eqref{inth} with $t=2+\delta$, $p=r_{2+\delta}$, and $\alpha=(\delta+2/r)/s$ and \eqref{intu} with $t=n$, $p=r$, and $\gamma=(n+1/2-2/r)/s$, we can dominate our term by the right hand side of \eqref{us2} since 
\[\frac{\alpha}{1-\gamma}=\frac{\delta+2/r}{\delta+2/r-1/2}\le 2.\]

Next, we consider the terms involving $\mu^{\frac12}|[\tilde u]|_{l,r,2}$ with $0\le l<n$. In this case we choose $r=\infty$, $r_j=2$, and $r_i=\infty$ for all $i\not=j$. We use \eqref{inth} with $t=i$ and $\alpha=\theta_i$, \eqref{inth} with $t=j$ and $\alpha=\theta_j'':=(j-1-1/p_j)/s$, and \eqref{intu} with $t=l$ and $\gamma=(l+1/2)/s$. Note that $\gamma\in (0,1)$ since $l\le s-1$, and $\theta''_j\le 1$ since $p_{s+2},p_{n+2}\in\{0,1\}$. Since $\sum\theta_ip_i+\theta_j''p_j\le(s-l)/s$ and
\[\frac{s-l}{s}\frac{1}{1-\gamma}=\frac{s-l}{s-l-1/2}\le 2\]
we can dominate our terms by the right hand side of \eqref{us2}.

The terms involving $\tilde q$ can be handled analogously; we only present the case $l=n$. Using \eqref{inth} with $t=i$, $p=\infty$, and $\alpha=\theta_i$ and \eqref{q2} with $t=n$ and $\beta=n/s$ we obtain
\begin{equation*}
\begin{aligned}
\frac{1}{\mu^{\frac12}}r_h|h|_{1+\delta,\infty,2}|\tilde q|_{n,2}\le r_h\big(1+|h|_{s+2,2}\big)\frac{1}{\mu^{\frac12}}\|\tilde q\|_{2} + \eps\frac{1}{\mu^{\frac12}}|\tilde q|_{s,2}
\end{aligned}
\end{equation*}
since 
\[\frac{\theta_{1+\delta}}{1-\beta}=\frac{\delta}{s}\frac{s}{s-n}=1.\]

Next, we consider the terms involving $\mu_b/\mu^{\frac12}|[\tilde u]|_{l,r,2}$. We choose $r=3$, $r_j=6$, and $r_i=\infty$ for all $i\not=j$ and apply \eqref{inth} with $t=i$ and $\alpha=\theta_i$, \eqref{inth} with $t=j$ and $\alpha=\theta_j''':=(j-1-1/(3p_j))/s$, and \eqref{intu} with $t=l$ and $\gamma=(l-1/6)/s$. This gives
\begin{equation*}
\begin{aligned}
\frac{\mu_b}{\mu^{\frac12}}r_h \prod_{i\in I_{s+1}}|h|_{i,p_ir_i,2}^{p_i}|[\tilde u]|_{l,r,2}&\le r_h\Big(1+\Big(\frac{\mu_b}{\mu}\Big)^{\frac12}|h|_{s+2,2}^{\sum\theta_ip_i+\theta'''_jp_j}\Big)^{\frac{1}{1-\gamma}}\mu_b^{\frac12}\|\nabla\tilde u\|_{2}\\
&\quad+ \eps\mu^{\frac12}_b|\tilde u|_{s+1,2}.
\end{aligned}
\end{equation*}
Since $\sum\theta_ip_i+\theta'''_jp_j\le (s-l-1/3)/s$ and
\[\frac{s-l-1/3}{s}\frac{1}{1-\gamma}=\frac{s-l-1/3}{s-l+1/6}\le 1\]
we can dominate our terms by the right hand side of \eqref{us2}.

\medskip
\noindent (b.4) It remains to analyse the terms $\hat f_1$ and $\hat f_2$. Let $s=n+\delta$ with $n\in\N_{\ge1}$ and $\delta\in [0,1)$. We have
\begin{equation*}
\begin{aligned}
&\frac{1}{\mu_b^{\frac12}}|\hat f_1|_{s-1,2}+\mu^{\frac12}_b|\hat f_2|_{s,2}\le \frac{1}{\mu_b^{\frac12}}|\tilde f_1|_{s-1,2}+\mu_b^{\frac12}|\tilde f_2|_{s,2}\\
&\quad+\|\tilde e-\hat e\|_{\infty}r_h \Big(\mu_b^{\frac12}|\tilde u|_{s+1,2}+\frac{1}{\mu_b^{\frac12}}|\tilde\pi|_{s,2}\Big)\\
&\quad+ \mu_b^{\frac12}r_h\bigg(\sum_{(p_i,l)\in J_{s+2}^{s+1,s}}\prod_{i\in I_{s+2}}\kappa^{\frac{1}{r_i}}|h|_{i,p_ir_i,2}^{p_i}|\tilde u|_{l,r,2}+b(\delta)\kappa^{\frac{1}{r_{1+\delta}}}|h|_{1+\delta,r_{1+\delta},2}|\tilde u|_{n+1,r}\bigg)\\
&\quad+\frac{1}{\mu_b^{\frac12}}r_h\bigg(\sum_{(p_i,l)\in J_{s}^{s,s-1}}\prod_{i\in I_{s}}\kappa^{\frac{1}{r_i}}|h|_{i,p_ir_i,2}^{p_i}|\tilde \pi|_{l,r,2}+b(\delta)\kappa^{\frac{1}{r_{1+\delta}}}|h|_{1+\delta,r_{1+\delta},2}|\tilde \pi|_{n,r}\bigg),
%&\quad+ \mu^{\frac12}\sum_{(p_i,l)\in J_{s+1}^{i_*}}r_h \prod_{i\in I_{s}}|h|_{i,p_ir_i,*}^{p_i}|\tilde v|_{l,r,*}\\
%&\quad+ \frac{1}{\mu^{\frac12}}\sum_{(p_i,l)\in J_{s}^{i_{**}}}r_h \prod_{i\in I_{s}}|h|_{i,p_ir_i,*}^{p_i}|\tilde q|_{l,r,*}\\
\end{aligned}
\end{equation*}
where $b(\delta)$ and $2\le r,r_i\le\infty$ are chosen as above. Here, we can proceed analogously to (b.3). The main difference is that we have to apply the interpolation estimate \eqref{intubulk} for $a=\nabla\tilde u,\tilde \pi$.

To begin with, let us consider the terms involving $|\tilde u|_{l,r,2}$ with $l=s$ if $\delta=0$ and $l\in\{n,s,n+1\}$ if $\delta>0$. We choose $r=2$ and $r_i=\infty$ for all $i$ except for the case $l=n,s<2,p_{2+\delta}=1$ which has to treated separately as above. Excluding this case we apply \eqref{inth} with $t=i$ and $\alpha=\theta_i$ and \eqref{intubulk} with $t=l-1$ and $\gamma=(l-1)/s$ and obtain
\begin{equation*}
\begin{aligned}
\mu_b^{\frac12}r_h\prod_{i\in I_{s+2}}\kappa^{\frac{1}{r_i}}|h|_{i,p_ir_i,2}^{p_i}|\tilde u|_{l,r,2}&\le r_h\Big(1+\kappa^{\frac12}|h|_{s+2,2}^{\sum\theta_ip_i}\Big)^{\frac{1}{1-\gamma}}\mu_b^{\frac12}\|\nabla\tilde u\|_{2}\\
&\quad+ \eps\mu^{\frac12}_b|\tilde u|_{s+1,2}.
\end{aligned}
\end{equation*}
Since $\sum\theta_ip_i\le (s+1-l)/s$ and
\[\frac{s-l+1}{s}\frac{1}{1-\gamma}=1\]
we can dominate our terms by the right hand side of \eqref{us2}. In the special case we apply  \eqref{intubulk} with $t=n-1=0$, $p=r:=\min(2/(1-\delta),3)$, and $\gamma=(3/2-3/r)/s$  and \eqref{inth} with $t=2+\delta$, $p=2r/(r-2)$, and $\alpha=(\delta+2/r)/s$. Since
\[\frac{\alpha}{1-\gamma}=\frac{\delta+2/r}{\delta+3/r-1/2}\le 2\]
we can again dominate our term by the right hand side of \eqref{us2}. 

Next, we consider the terms involving $|\tilde u|_{l,r,2}$ with $0\le l\le s-1$. We choose $r=\infty$, $r_j=2$, and $r_i=\infty$ for all $i\not=j$ and apply \eqref{inth} with $t=i$ and $\alpha=\theta_i$, \eqref{inth} with $t=j$ and $\alpha=\theta_j'''$, and \eqref{intubulk} with $t=l-1$ and $\gamma=(l+1/2)/s$. Since $\sum\theta_ip_i+\theta_j'''p_j\le (s-l)/s$ and
\[\frac{s-l}{s}\frac{1}{1-\beta}=\frac{s-l}{s-l-1/2}\le 2\]
we can again dominate our term by the right hand side of \eqref{us2}.

Estimating the terms involving the bulk pressure is proceeds completely analogously.

We want to absorb the highest order terms on the right hand sides with a factor $\|\tilde e-\hat e\|_\infty$ on the left hand side of the estimate proven in step (a). The constant in the latter estimate is of the form $r_h$, so we need to assume $\|\tilde e-\hat e\|_\infty$ to be of the form $r_h^{-1}$. In this case the lemma is proved.

\medskip
\noindent (c) Now, let us consider the general case. For $x_0\in\Omega$ and $R>0$ we have
\begin{equation*}
\begin{aligned}
\|\tilde e-\tilde e(x_0)\|_{\infty,B_R(x_0)}\le R\|\nabla\tilde e\|_{\infty,\Omega}\le R\,r_h(1+\|\nabla^2h\|_{\infty}).
\end{aligned}
\end{equation*}
Since we want to apply step (b), we choose $R$ to be of the form $(r_h(1+\|\nabla^2h\|_{\infty}))^{-1}$ and cover $\Omega$ by balls of radius $R/2$; the number of balls we need to do so is of the order 
\[\frac{\diam(\Omega)}{R}=\diam(\Omega)\,r_h\,(1+\|\nabla^2h\|_{\infty}).\]
Let $\phi$ be a smooth function that is $1$ in $B_1(0)$ and $0$ in $B_2^c(0)$. Now, we fix one of the finitely many balls $B_R(x_0)$ and multiply our solution with $\phi_R:=\phi((\cdot-x_0)/R)$. The product is again a solution, however, with right hand sides perturbed by lower order derivatives of the original solution multiplied by suitable powers of $1/R$. These lower order terms can again be handled by interpolation and absorption. For dimensional reasons this must lead to an additional factor of the form $r_h(1+\|\nabla^2h\|_{\infty}^s)$ in front of the lower order terms in our estimate.

\medskip
\noindent (d) So far, we proved the estimate
\begin{equation*}
\begin{aligned}
&\mu_b^\frac12|\tilde u|_{s+1} + \mu^\frac12 |\tilde v|_{s+1} + \frac{1}{\mu_b^\frac12} |\tilde \pi|_{s} + \frac{1}{\mu^\frac12} |\tilde q|_{s}\\
&\quad\le r_h\Big(\frac{1}{\mu_b^\frac12}|\tilde f_1|_{s-1}+ \mu_b^\frac12|\tilde f_2|_{s}+\frac{1}{\mu^\frac12}|\tilde f_3^\top|_{s-1} + \frac{1}{\mu_b^\frac12} |\tilde f_3^\perp|_{s-\frac12} + \mu^\frac12|\tilde f_4|_{s}\Big)\\
&\quad\quad + r_h\big(1+|h|_{2,\infty}^s+c_{s+\frac32}|h|_{s+2,2}^{2+\frac{1}{s}}\big)\\
&\quad\quad\quad\times\Big(\mu_b^\frac12\|\nabla \tilde u\|_{2} + \mu^\frac12\|\nabla^g \tilde v\|_{2} + \frac{1}{\mu_b^\frac12}\|\tilde \pi\|_{2}+\frac{1}{\mu^\frac12}\|\tilde q\|_{2}\Big).
\end{aligned}
\end{equation*}
We have to transform this estimate back to the original domain. Equivalently, we can replace the Euclidean norms by norms which are defined with respect to the perturbed metric $\tilde e$; we denote such norms with an index $\tilde e$. We begin by replacing the norms of the data. From \eqref{Tk} and \eqref{Tdelta} as well as Lemma \ref{lemma:prod} and Corollary \ref{cor:analytic} we infer that
\begin{equation*}
\begin{aligned}
|\tilde f_1|_{s-1,2}&\le r_h\bigg(\sum_{(p_i,l)\in J_{s}^{s-1,s-2}}\prod_{i\in I_{s+2}}\kappa^{\frac{1}{r_i}}|h|_{i,p_ir_i,2}^{p_i}|\tilde f_1|_{l,r,2}\\
&\quad\quad+b(\delta)\kappa^{\frac{1}{r_{1+\delta}}}|h|_{1+\delta,r_{1+\delta},2}|\tilde f_1|_{n-1,r}+{_{\tilde e}|\tilde f_1|_{s-1,2}}\bigg),
\end{aligned}
\end{equation*}
where $s=n+\delta$ with $n\in\N_{\ge1}$, $\delta\in [0,1)$, and $b(\delta)$ is defined as above. Now, in the terms involving $|\tilde f_1|_{l,r}$ with $l<s-1$ we replace $\tilde f_1$ by the left hand side of \eqref{eqn:tildesystemh}$_1$ and take into account the representation of this left hand side derived in \eqref{eqn:tildesystem2}$_1$ (with $\hat e$ replaced by $e$). This way we see that
\begin{equation*}
\begin{aligned}
&|\tilde f_1|_{s-1,2}\le r_h\,{_{\tilde e}|\tilde f_1|_{s-1,2}}\\
&\quad+ \mu_br_h\bigg(\sum_{(p_i,l)\in J_{s+2}^{s+1,s}}\prod_{i\in I_{s+2}}\kappa^{\frac{1}{r_i}}|h|_{i,p_ir_i,2}^{p_i}|\tilde u|_{l,r,2}+b(\delta)\kappa^{\frac{1}{r_{1+\delta}}}|h|_{1+\delta,r_{1+\delta},2}|\tilde u|_{n+1,r}\bigg)\\
&\quad+\frac{1}{\mu_b}r_h\bigg(\sum_{(p_i,l)\in J_{s}^{s,s-1}}\prod_{i\in I_{s}}\kappa^{\frac{1}{r_i}}|h|_{i,p_ir_i,2}^{p_i}|\tilde \pi|_{l,r,2}+b(\delta)\kappa^{\frac{1}{r_{1+\delta}}}|h|_{1+\delta,r_{1+\delta},2}|\tilde \pi|_{n,r}\bigg).
\end{aligned}
\end{equation*}
These are exactly the terms we dealt with in (b.4). While in (b.4) the highest order derivatives were absorbed on the left hand side, here, these terms are implicitly contained in the expression ${_{\tilde e}|\tilde f_1|_{s-1,r}}$. The terms involving $\tilde f_2,\ldots,\tilde f_4$ can be handled completely analogously. We obtain
\begin{equation}\label{hurtz}
\begin{aligned}
&\mu_b^\frac12|\tilde u|_{s+1,2} + \mu^\frac12 |\tilde v|_{s+1,2} + \frac{1}{\mu_b^\frac12}|\tilde \pi|_{s,2} + \frac{1}{\mu^\frac12} |\tilde q|_{s,2}\\
&\quad\le r_h\Big(\frac{1}{\mu_b^\frac12}{_{\tilde e}|\tilde f_1|_{s-1,2}}+ \mu_b^\frac12{_{\tilde e}|\tilde f_2|_{s,2}}+\frac{1}{\mu^\frac12}{_{\tilde e}|\tilde f_3^\top|_{s-1,2}} + \frac{1}{\mu_b^\frac12} {_{\tilde e}|\tilde f_3^\perp|_{s-\frac12,2}}\\
&\quad\quad\quad+ \mu^\frac12{_{\tilde e}|\tilde f_4|_{s,2}}\Big) + r_h\big(1+|h|_{2,\infty}^s+c_{s+\frac32}|h|_{s+2,2}^{2+\frac{1}{s}}\big)\\
&\quad\quad\quad\times\Big(\mu_b^\frac12\|\nabla \tilde u\|_{2} + \mu^\frac12\|\nabla^g \tilde v\|_{2} + \frac{1}{\mu_b^\frac12}\|\tilde \pi\|_{2}+\frac{1}{\mu^\frac12}\|\tilde q\|_{2}\Big).
\end{aligned}
\end{equation}
Next, we take care of the terms on the left hand side. Using again \eqref{Tk} and \eqref{Tdelta} we prove that
\begin{equation*}
\begin{aligned}
{_{\tilde e}|\tilde u|_{s+1,2}}\le r_h\sum_{(p_i,l)\in J_{s+2}^{s+1,s+1}}\prod_{i\in I_{s+2}}\kappa^{\frac{1}{r_i}}|h|_{i,p_ir_i,2}^{p_i}|\tilde u|_{l,r,2}.
\end{aligned}
\end{equation*}
Again, these are exactly the terms we dealt with in (b.4); here, however the highest order derivative is estimated using \eqref{hurtz} instead of being absorbed. The term involving $\tilde \pi$ can be handled analogously using (b.4), while the terms involving $\tilde v$ and $\tilde q$ can be dealt with using (b.3). From these arguments and by a transformation of the left hand side and the data terms to the original domain $\Omega\setminus\Gamma_h$ and $\Gamma_h$, respectively, we obtain
\begin{equation*}
\begin{aligned}
&\mu_b^\frac12|u|_{s+1,2} + \mu^\frac12 |v|_{s+1,2} + \frac{1}{\mu_b^\frac12} |\pi|_{s,2} + \frac{1}{\mu^\frac12} |q|_{s,2}\\
&\quad\le r_h\Big(\frac{1}{\mu_b^\frac12}|f_1|_{s-1,2}+ \mu_b^\frac12|f_2|_{s,2}+\frac{1}{\mu^\frac12}|f_3^\top|_{s-1,2} + \frac{1}{\mu_b^\frac12}|f_3^\perp|_{s-\frac12,2} + \mu^\frac12|f_4|_{s,2}\Big)\\
&\quad\quad + r_h\big(1+|h|_{2,\infty}^s+c_{s+\frac32}|h|_{s+2,2}^{2+\frac{1}{s}}\big)\\
&\quad\quad\quad\times\Big(\mu_b^\frac12\|\nabla \tilde u\|_{2} + \mu^\frac12\|\nabla^g \tilde v\|_{2} + \frac{1}{\mu_b^\frac12}\|\tilde \pi\|_{2}+\frac{1}{\mu^\frac12}\|\tilde q\|_{2}\Big).
\end{aligned}
\end{equation*}
It remains to treat the lower order terms on the right hand side. Obviously, for the pressure terms we have 
\begin{equation*}
\begin{aligned}
\frac{1}{\mu_b^\frac12}\|\tilde \pi\|_{2}+\frac{1}{\mu^\frac12}\|\tilde q\|_{2}\le r_h\Big(\frac{1}{\mu_b^\frac12}{_{\tilde e}\|\tilde \pi\|_{2}}+\frac{1}{\mu^\frac12}{_{\tilde e}\|\tilde q\|_{2}}\Big)
\end{aligned}
\end{equation*}
Concerning the velocities, we compute using \eqref{Tk}
\begin{equation*}
\begin{aligned}
\|\nabla\tilde u\|_2&\le r_h\,{_{\tilde e}\|\nabla^{\tilde e}\tilde u\|_2}+r_h\big(1+\kappa^{\frac13}|h|_{2,3}\big)\|u\|_6\\
&\le r_h\,{_{\tilde e}\|\nabla^{\tilde e}\tilde u\|_2}+r_h(1+\kappa^{\frac13}|h|_{s+2,2}^{\frac{1}{3s}}\big)\|\nabla u\|_2,
\end{aligned}
\end{equation*}
where we used \eqref{inth} and Sobolev's embedding in $\Omega$ for the second inequality, and
\begin{equation*}
\begin{aligned}
\mu^\frac12\|\nabla^g\tilde v\|_2&\le \mu^\frac12r_h\,{_{\tilde g}\|\nabla^{\tilde g}\tilde v\|_2}+r_h\big(1+|h|_{2,4}\big)\mu^\frac12\|v\|_4\\
&\le \mu^\frac12 r_h\,{_{\tilde g}\|\nabla^{\tilde g}\tilde v\|_2}+r_h\Big(1+\frac{\mu^\frac12}{\mu_b^\frac12}|h|_{s+2,2}^{\frac{1}{2s}}\Big)\mu_b^\frac12\|\nabla u\|_2,
\end{aligned}
\end{equation*}
where we used \eqref{inth} and Lemma \ref{lemma:trace} for the second inequality. Here, we avoided absorption since this would raise the exponent of $|h|_{s+2,2}$ unnecessarily.

This, finally, finishes the proof.
\qed

\medskip
At first sight, it might seem easier to prove the preceding theorem for integer $s$ first and apply interpolation afterwards. However, first, this approach would slighty worsen the bound on the constant $c_{-s}$. And second, this would not avoid the need for estimating products of functions in $|\cdot|_{t,2}$ for non-integer $t$ since the equations \eqref{eqn:tildesystemh}$_{4,6}$ must be analyzed in $H^{s-\frac12}(\Gamma)$ and $H^{s+\frac12}(\Gamma)$, respectively.

Finally, we analyze the constants in \eqref{est:regs2}.
\begin{theorem}\label{s2h}
 Let $s\ge 1$ and $h\in H^{s+2}(\Gamma)$ such that $\|h\|_{L^\infty(\Gamma)}<\kappa/2$. For all $(f_1,\ldots,f_4)\in\F^s_\nu(\Gamma_h)$ there exists a unique solution $(u,\pi,q)\in\E^s_\nu(\Gamma_h)$ of \eqref{s2}. The estimate \eqref{est:regs2} holds with the constant $c$ being bounded by $r_h$ and the constant $c_{-s}$ being bounded by
\[r_h\big(1+c_{\frac12}|h|_{2,\infty}^s|h|_{s+2,2}^{\frac{1}{2s}}+c_{s+\frac12}|h|_{s+2,2}^{2+\frac{1}{2s}}\big);\]
 note that for $s>1$ the we can get rid of the term involving $|h|_{2,\infty}^s$ by interpolation.
 \end{theorem}
\proof We can proceed exactly like in the proof of Theorem \ref{s1h}. Here, however, we don't need step (b.1), leading to a slightly lower exponent, and step (b.2) becomes a bit simpler since $\tilde w$ is prescribed and, in particular, absorption is not needed.
\qed

\appendix
\section{Mechanics of fluid interfaces}\applabel{app:int}
In this appendix, for the sake of completeness, we repeat the basic facts concerning the mechanics of fluid interfaces; more detailed information can be found for instance in \cite{scriven,aris,guven}. 

We begin with the general kinematics of fluid interfaces. Let $\Gamma_t\subset\R^3$ be an interface moving with velocity $u=v+w\,\nu$, where $v$ is tangential and $\nu$ is a fixed unit normal field on $\Gamma_t$, and let $g=g_t$ be the Riemannian metric on $\Gamma_t$ induced by the ambient Euclidean space. Consider the familiy of diffeomorphisms $\phi_{t,s}:\Gamma_{t}\rightarrow\Gamma_s$ associated with the vector field $u$, that is, $\phi_{t,t}=\id_{\Gamma_{t}}$ and $\pa_s\phi_{t,s}=u\circ\phi_{t,s}$. We denote by $Df/Dt$ the \emph{material derivative} of a scalar field $f$ on $\Gamma_t$ with respect to the vector field $u$, that is,
\[\frac{Df}{Dt}\Big|_t:=\pa_s|_{s=t}f\circ\phi_{t,s}.\]
Throughout this appendix we shall work in \emph{convected coordinates} $(x^\alpha)$ on $\Gamma_t$, that is, $x^\alpha|_t=x^\alpha|_s\circ\phi_{t,s}$; let $(\pa_\alpha)$ be the associated coordinate vector fields. Note that \[\pa_\alpha|_s=(d\phi_{t,s}\pa_\alpha|_t)\circ\phi_{t,s}^{-1}=(\pa_\alpha\phi_{t,s})\circ\phi_{t,s}^{-1}\]
and thus, by Schwarz' theorem,
\begin{equation*}
\frac{D}{Dt}\Big|_t\,\pa_\alpha=\pa_s|_{s=t}\pa_\alpha\phi_{t,s}=\pa_\alpha u; 
\end{equation*}
here, we take the material derivative of each Cartesian component $\pa_\alpha^i$ of the vector field $\pa_\alpha$. Analogously to the bulk case, the \emph{surface rate-of-strain tensor}, $\D u$, is given by the rate of change of (infinitesimal) lengths along the flow or, in other words, by the Lie derivative of the Riemannian metric. More precisely, if the tangent vector field $X_t$ on $\Gamma_{t}$ is transported by the flow, that is, $X_s=(d\phi_{t,s}X_t)\circ\phi_{t,s}^{-1}$, then
\[\D u_t(X_t,X_t):=\frac12\frac{D}{Dt}g_t(X_t,X_t),\]
where $|\cdot|$ denotes the Euclidean norm.
%  For tangent vectors $X,Y$ on $\Gamma_t$ consider the pullback $g_{t,s}(X,Y):=d\phi_{t,s}X\cdot d\phi_{t,s}Y$ of the metric on $\Gamma_s$ (induced by the ambient Euclidean space) under $\phi_{t,s}$. 
In convected coordinates, we have 
\begin{equation*}
\begin{aligned}
(\D u)_{\alpha\beta}&=\frac12\frac{D}{Dt}g_{\alpha\beta}=\frac12\frac{D}{Dt}\ \pa_\alpha\cdot\pa_\beta\\
&=\frac12\big(\pa_\alpha(v+w\,\nu)\cdot\pa_\beta+\pa_\alpha\cdot\pa_\beta(v+w\,\nu)\big)\\
&=\frac12\big(v_{\beta;\alpha}+w\,\nu_{,\alpha}\cdot\pa_\beta+v_{\alpha;\beta}+w\,\nu_{,\beta}\cdot\pa_\alpha\big)\\
&=\frac12(v_{\beta;\alpha}+v_{\alpha;\beta}) - w\,k_{\alpha\beta},
\end{aligned}
\end{equation*}
where we used the identity $k_{\alpha\beta}=-\nu_{,\alpha}\cdot\pa_\beta$ for the third equality. Note that the material derivative of the components in convected coordinates of any tangential tensor field yields the Lie derivative of this tensor field along the flow; however, this is not true for hybrid tensor fields which are partly expressed in Cartesian components, see below. Since $g^{\alpha\beta}g_{\beta\gamma}=\delta_\gamma^\alpha$ is a constant we have
\[g_{\beta\gamma}\frac{D}{Dt}g^{\alpha\beta}=-g^{\alpha\beta}\frac{D}{Dt}g_{\beta\gamma},\quad\text{and thus}\quad\frac{D}{Dt}g^{\alpha\beta}=-g^{\alpha\gamma}g^{\beta\delta}\frac{D}{Dt}g_{\gamma\delta}.\]
Again, since $|\nu|^2=1$ and $\nu\cdot\pa_\alpha=0$ are constants we have
\[\nu\cdot\frac{D}{Dt}\nu=0\quad\text{and}\quad\pa_\alpha\cdot\frac{D}{Dt}\nu=-\nu\cdot\frac{D}{Dt}\pa_\alpha=-(\nu\cdot\pa_\alpha v+ w_{,\alpha})=-k_{\alpha\beta}v^\beta-w_{,\alpha}.\]
From the chain rule $\pa_\alpha f|_s\circ\phi_{t,s}=\pa_\alpha(f\circ\phi_{t,s})$ we see that $D_t$ and $\pa_\alpha$ commute, and thus
\[\pa_\beta\cdot\frac{D}{Dt}\pa_\alpha\nu=\pa_\beta\cdot\pa_\alpha\frac{D}{Dt}\nu=-(k_{\alpha\gamma}v^\gamma+w_{,\alpha})_{;\beta}.\]
Hence, we have
\[\frac{D}{Dt}k_{\alpha\beta}=-\pa_\beta\cdot\frac{D}{Dt}\pa_\alpha\nu-\pa_\alpha\nu\cdot\frac{D}{Dt}\pa_\beta=w_{;\alpha\beta}+(k_{\alpha\gamma}v^\gamma)_{;\beta}+k_\alpha^\gamma v_{\gamma;\beta}-wk_{\alpha\gamma}k^{\gamma}_\beta.\]
From this we see that the material derivative of the mean curvature $H=g^{\alpha\beta}k_{\alpha\beta}$ reads
\begin{equation*}
\begin{aligned}
\frac{D}{Dt}H&=k_{\alpha\beta}\frac{D}{Dt}g^{\alpha\beta} + g^{\alpha\beta}\frac{D}{Dt}k_{\alpha\beta}=\Delta_g w +w|k|^2 + k_{\alpha\beta;}^{\ \ \ \alpha}v^\beta\\
&=\Delta_g w +w (H^2-2K) + dH(v).
\end{aligned}
\end{equation*}
Here, we used the fact that the tensor $k_{\alpha\beta;\gamma}$ is symmetric with respect to all indices, which is the content of the \emph{Codazzi-Mainardi equations}. Finally, we compute the material derivative of the volume element $dA=dA_t$ on $\Gamma_t$. To this end, we note that
\[\phi_{t,s}^{-1}(dA_s)=\frac{\sqrt{\det g_s}\circ\phi_{t,s}}{\sqrt{\det g_t}}\,dA_t,\]
where the determinants are taken with respect to convected coordinates. Hence, we compute
\[\pa_{s}|_{s=t}\sqrt{\det g_{s}}\circ\phi_{t,s}=\sqrt{\det g_{t}}\,(g_t)^{\alpha\beta}\frac12\frac{D}{Dt}(g_t)_{\alpha\beta}=\sqrt{\det g_{t}}\,(g_t)^{\alpha\beta}(\D u_t)_{\alpha\beta}\]
and obtain
\[\frac{D}{Dt}\,dA=(\Div_g v-w\,H)\,dA=\DIV u\,dA.\]

Next, let us have a look at conservation of mass and linear momentum on a fluid interface. Consider a mass density $\rho$ on $\Gamma_t$ that is being conserved and transported by the vector field $u$; this means that for all subsets $M_t\subset\Gamma_t$ moving along the flow we have
\[\frac{d}{dt}\int_{M_t}\rho\,dA=0.\]
Interchanging differentiation in time and integration in space gives
\[\frac{d}{dt}\int_{M_t}\rho\,dA=\int_{M_t}\Big(\frac{D\rho}{Dt}+\rho\,\DIV u\Big)\,dA;\]
since $M_t$ is arbitrary, this gives conservation of mass in differential form
\[\frac{D\rho}{Dt}+\rho\,\DIV u=0\text{ on }\Gamma_t.\]
Again, let $M_t\subset\Gamma_t$ be (a smooth bounded domain) moving along the flow. Then, for $i=1,2,3$, balance of linear momentum in integral form reads
\begin{equation}\label{conm}
\frac{d}{dt}\int_{M_t}\rho\,u^i\,dA=\int_{\pa M_t}T^i_{\alpha}(\nu_{\pa M_t})^\alpha\,ds+\int_{M_t}f^i\,dA,
\end{equation}
where $T$ is the \emph{surface stress tensor}, $\nu_{\pa M_t}$ is the outer unit normal to $M_t$, and $f$ is a given surface force density. Gauss' theorem gives
\[\int_{\pa M_t}T^i_{\alpha}(\nu_{\pa M_t})^\alpha\,ds=\int_{M_t} g^{\alpha\beta}T^i_{\alpha;\beta}\,dA.\]
The integrand on the right hand side, however, is nothing but the $i$-th component of the surface divergence $\DIV T$. Interchanging differentiation and integration on the left hand side of \eqref{conm} gives
\[\frac{d}{dt}\int_{M_t}\rho\,u\,dA=\int_{M_t}\Big(\frac{D\rho}{Dt} u + \rho \frac{Du}{Dt} +\rho\, u\,\DIV u\Big)\,dA=\int_{M_t}\rho \frac{Du}{Dt}\,dA,\]
where we used conservation of mass for the second identity. Since $M_t$ is arbitrary, we obtain conservation of linear momentum in differential form
\[\rho\frac{Du}{Dt}=\DIV T + f\text{ on }\Gamma_t.\]
It is interesting to note that this expression has to be interpreted with respect to Cartesian coordinates in $\R^3$ since, as we already pointed out, the material derivative must be understood component-wise. In other words, while $Du/Dt$ is the Lie derivative of each Cartesian component $u^i$ along the flow, it is \emph{not} the Lie derivative of the vector field $u$ (which would vanish anyway due to skew-symmetry of the Lie bracket). In particular, balance of linear momentum is not a tensorial postulate; cf. \cite{marsdenhughes94}.

Finally, let us compute the stresses resulting from the Canham-Helfrich energy \eqref{eqn:ch}; this can be done by exploiting translational symmetry of the energy. Consider a fluid interface $\Gamma_t$ that is being translated, that is, $\phi_{t,s}(x)=x+a(t-s)$ for some $a\in\R^3$. Then, we have $v^\alpha=g^{\alpha\beta}a\cdot\pa_\beta$ and $w=a\cdot\nu$. Since the energy \eqref{eqn:ch} is invariant under translations, for all translated subsets $M_t\subset\Gamma_t$ we have
\begin{equation*}
\begin{aligned}
0&=\frac{d}{dt}\frac{\kappa}{2}\int_{M_t}(H-C_0)^2\,dA\\
&=\kappa\int_{M_t}\Big((H-C_0)\big(\Delta_g w +w (H^2-2K) + dH(v)\big)\\
&\quad\quad+\frac{(H-C_0)^2}{2}(\Div_g v-w\,H)\Big)\,dA.
\end{aligned}
\end{equation*}
Note that the integrand can be written in the form
\begin{equation*}
\begin{aligned}
&\big(\Delta_g H+ H(H^2/2-2K) + C_0(2K-HC_0/2)\big)w\\
&+ \Div_g\big((H-C_0)\grad_g w-w\grad_g (H-C_0) + v (H-C_0)^2/2\big)\\
&\quad=\sum_i a^i\big(\kappa^{-1}\grad_{L_2} F\, \nu^i +\Div_g\big((H-C_0)\grad_g \nu^i-\nu^i\grad_g (H-C_0)\\
&\quad\quad+ (H-C_0)^2/2\,g^{\alpha\beta}\pa_\beta^i\,\pa_\alpha\big).
 \end{aligned}
\end{equation*}
Since $M_t$ and $a$ are arbitrary, we conclude that the $L_2$-gradient of the Canham-Helfrich energy can be written as the surface divergence of the stress tensor
\[{^hT^i_\alpha}=\kappa\big(-(H-C_0) k_\alpha^\beta\pa_\beta^i- (H-C_0)_{,\alpha}\nu^i+(H-C_0)^2/2\,\pa_\alpha^i\big),\]
that is
\begin{equation}\label{divft}
 -\grad_{L_2} F\, \nu=\DIV {^hT}.
\end{equation}
While this formula looks very natural from a physical point of view, from a mathematical point of view it is not obvious at all that the $L_2$-gradient can be written in divergence form. In fact, when this and related identities were discovered in \cite{riviere08} by purely mathematical (and rather complicated) considerations, it was a major advance in the analysis of Willmore surfaces. Finally, let us verify the identity \eqref{divft} directly:
\begin{equation*}
\begin{aligned}
\kappa^{-1}\DIV {^hT}&=-\big((H-C_0)k^{\alpha\beta}\big)_{;\alpha}\pa_\beta-(H-C_0)k^{\alpha\beta}k_{\alpha\beta}\,\nu - \Delta_g H\,\nu\\
&\quad\,+ (H-C_0)_{,\alpha}k^{\alpha\beta}\pa_\beta + g^{\alpha\beta}\big((H-C_0)^2/2\big)_{,\beta}\pa_\alpha+ (H-C_0)^2/2\,H\nu\\
&=-\nu\,\big(\Delta_gH + H(H^2/2-2K)+C_0(2K-HC_0/2)\big)\\
&\quad\, + (H-C_0) (g^{\alpha\beta}H_{,\alpha}\pa_\beta -k^{\alpha\beta}_{\ \ \ ;\alpha}\pa_\beta)\\
&=-\nu\,\big(\Delta_gH + H(H^2/2-2K)+C_0(2K-HC_0/2)\big).
\end{aligned}
\end{equation*}
Here, we used the Codazzi-Mainardi equations for the second identity. 

\section{Covariant differentiation and curvature}\applabel{app:cov}
Let $e_{ij}$ be a Riemannian metric on a sufficiently smooth manifold $M$, and let $e^{ij}$ denote its matrix inverse. Then the corresponding Christoffel symbols are given by
\begin{equation*}%\label{eqn:christ}
{^{e}\Gamma}_{ij}^k=\frac{1}{2}e^{kl}(\pa_ie_{jl} + \pa_j e_{il}-\pa_le_{ij}),
\end{equation*}
Furthermore, for the total covariant derivatives of a vector field $Y$ and a $(1,1)$-tensor field $T$ we have
\begin{equation*}
\begin{aligned}
 (\nabla^e Y)^j_i&=\pa_iY^j+{^{e}\Gamma_{ik}^j}Y^k,\\
(\nabla^e T)_{ki}^j&=\pa_iT^j_k+{^{e}\Gamma_{il}^j}T^l_k-{^{e}\Gamma_{ik}^l}T^j_l,
\end{aligned}
\end{equation*}
while their divergences are given by contracting $i$ and $j$ (or $k$ using $e$). Furthermore, the $e$-symmetric part $D^e Y$ of $\nabla^e Y$ can be written as
\[(D^e Y)_i^j=\frac12 e^{jk}(e_{kl}(\nabla^e Y)_{i}^l+e_{il}(\nabla^e Y)_{k}^l).\]

Now, let $\tilde e$ be a second metric with corresponding Christoffel symbols $^{\tilde e}\Gamma_{ij}^k$. Recall that the difference $\Sigma_{ij}^k:={^{\tilde e}\Gamma_{ij}^k}-{^{e}\Gamma_{ij}^k}$ of two connections is a tensor, and note that\footnote{Here and in the following, we have a slight ambiguity of notation: It is not specified wether the covariant tensor $\tilde e_{ij}$ or the contravariant tensor $\tilde e^{ij}$ enters the formulae. This ambiguity, however, is irrelevant in the analysis.}
\begin{equation*}
\begin{aligned}
\Sigma &= \tilde e*\nabla^e\tilde e,\\
\nabla^e\Sigma &= \tilde e*(\nabla^e)^2\tilde e + (\nabla^e\tilde e)^2.
\end{aligned}
\end{equation*}
For a scalar function $f$ we have
\begin{equation*}
\begin{aligned}
(\grad_{\tilde e}f)^i&=\tilde e^{ij}\pa_j f=(\grad_{e}f)^i + (\tilde e^{ij}-e^{ij})\pa_j f,\\
\Div_{\tilde e}Y&=\pa_iY^i+{^{\tilde e}\Gamma_{ik}^i}Y^k=\Div_{e}Y + \Sigma_{ik}^iY^k=\Div_e Y + \tilde e*\nabla^e\tilde e*Y,\\
\Delta_{\tilde e} f &= \Div_{\tilde e}\grad_{\tilde e}f=\Delta_e f + \Div_e(\grad_{\tilde e}-\grad_{e})f + \Sigma^i_{ik}(\grad_{\tilde e}f)^k\\
&= \Delta_e f + (\tilde e-e)*r(\tilde e,e)*(\nabla^e)^2 f + r(\tilde e,e)*\nabla^e\tilde e*\nabla f.
\end{aligned}
\end{equation*}
Furthermore, we have
\[(D^{\tilde e} Y)_i^j=\frac12 \tilde e^{jk}(\tilde e_{kl}(\nabla^e Y)_{i}^l+\tilde e_{il}(\nabla^e Y)_{k}^l)+(\tilde e^2*\Sigma*Y)_i^j,\]
and thus
\begin{equation*}
D^{\tilde e}Y=D^{e}Y+(\tilde e-e)*r(\tilde e,e)*\nabla^e Y + r(\tilde e,e)*\nabla^e \tilde e *Y.
\end{equation*}
Similarly, we compute
\begin{equation*}
 \Div_{\tilde e} T=\Div_{e} T + (\tilde e-e)*r(\tilde e,e)*\nabla^e T + r(\tilde e,e)*\nabla^e\tilde e*T.
\end{equation*}
Furthermore, we can prove by induction that for any tensor field $T$ we have
\begin{equation}\label{Tk}
 (\nabla^{\tilde e})^k T=(\nabla^e)^k T + \sum_{p_i,l}\prod_{i=1}^k r(\tilde e,e)*((\nabla^e)^i \tilde e)^{p_i} (\nabla^e)^l T,
\end{equation}
where $p_i,l\in \N_{\ge 0}$ and $\sum_i i p_i+l=k$.
Using Riemannian normal coordinates with respect to $e$, we may assume in the following calculations that the Christoffel symbols ${^{e}\Gamma}_{ij}^k$ (but not their derivatives) vanish. 
Hence, we compute
\begin{equation*}
\begin{aligned}
((\nabla^{\tilde e})^2 Y)_{ij}^k&= \pa^2_{ij}Y^k + \pa_i{^{\tilde e}\Gamma_{jl}^k}\, Y^l + {^{\tilde e}\Gamma_{jl}^k}\,\pa_i Y^l + {^{\tilde e}\Gamma_{im}^k}(\pa_j Y^m+{^{\tilde e}\Gamma_{jl}^m}Y^l)\\
&\quad - {^{\tilde e}\Gamma_{ij}^l}(\pa_lY^k+{^{\tilde e}\Gamma_{lm}^k}Y^m)\\
&=((\nabla^{e})^2 Y)_{ij}^k + \pa_i({^{\tilde e}\Gamma_{jl}^k}-{^{e}\Gamma_{jl}^k})\, Y^l + {^{\tilde e}\Gamma_{jl}^k}\,\pa_i Y^l\\
&\quad + {^{\tilde e}\Gamma_{im}^k}(\pa_j Y^m+{^{\tilde e}\Gamma_{jl}^m}Y^l) - {^{\tilde e}\Gamma_{ij}^l}(\pa_lY^k+{^{\tilde e}\Gamma_{lm}^k}Y^m)\\
&=((\nabla^{e})^2 Y)_{ij}^k + (\nabla^e\Sigma)_{jli}^k Y^l + \Sigma_{jl}^k ({\nabla^e}Y)_i^l + \Sigma_{im}^k(({\nabla^e}Y)_j^m+\Sigma_{jl}^mY^l)\\
&\quad -\Sigma_{ij}^l(({^e\nabla}Y)_l^k+\Sigma_{lm}^kY^m).
\end{aligned}
\end{equation*}
In particular, we have
\begin{equation*}
\begin{aligned}
\Delta_{\tilde e} Y&=\Delta_{e} Y + (\tilde e-e)*(\nabla^e)^2 Y+ \tilde e *(\nabla^e\Sigma*Y+ \Sigma*{\nabla^e}Y + \Sigma*\Sigma*Y)\\
&=\Delta_{e} Y + (\tilde e-e)*r(\tilde e,e)*(\nabla^e)^2 Y + r(\tilde e,e)*(\nabla^e)^2\tilde e*Y\\ 
&\quad+ r(\tilde e,e)*(\nabla^e\tilde e)^2*Y
+ r(\tilde e,e)*\nabla^e\tilde e*\nabla^e Y.
\end{aligned}
\end{equation*}

Now, let $\Gamma$ be an orientable submanifold of $M$ of codimension $1$, and let $\nu_{e}$ and $\nu_{\tilde e}$ be equally oriented unit normal fields on $\Gamma$ with respect to $e$ and $\tilde e$, respectively. Employing a Gram-Schmidt orthonormalization it is not hard to see that
\begin{equation*}
 \nu_{\tilde e}=\nu_e + (\tilde e-e)*r(\tilde e,e).
\end{equation*}
Thus, we have
\begin{equation*}
\nabla^{\tilde e} \nu_{\tilde e}=\nabla^e\nu_{\tilde e} + \Sigma *\nu_{\tilde e} = \nabla^e\nu_{e} + r(\tilde e,e)*\nabla^e\tilde e.
\end{equation*}
In view of $(k_{\tilde e})_{\alpha\beta}=-\tilde e_{ij}(\nabla^{\tilde e}_\alpha\nu_{\tilde e})^i\pa_\beta^j$, where greek and latin indices refer to coordinates on $\Gamma$ and in $M$, respectively, we deduce
\begin{equation*}
\begin{aligned}
k_{\tilde e}&=k_e+(\tilde e-e)*k_e + r(\tilde e,e)*\nabla^e\tilde e,\\
H_{\tilde e}&= H_e + (\tilde e-e)*r(\tilde e,e)*k_e + r(\tilde e,e)*\nabla^e\tilde e,\\
K_{\tilde g}&=\det(\tilde g^{\alpha\delta}(k_{\tilde e})_{\delta\beta})=\det\big(g^{\alpha\delta}(k_e)_{\delta\beta} + (\tilde e-e)*r(\tilde e,e)*k_e + r(\tilde e,e)*\nabla^e\tilde e\big)\\
&=K_e+r(\tilde e,e)*\big((\tilde e-e)*k^2,k*\nabla\tilde e,(\nabla\tilde e)^2\big).
\end{aligned}
\end{equation*}

Finally, we consider the special case of $\Gamma$ being a closed surface in $\R^3$, and we let $\Gamma_h$ be the graph of a height function $h$ on $\Gamma$ as in \Subsecref{stokesvar}. Note that
\begin{equation}\label{intah}
\begin{aligned}
\phi_h^{-1}(dA_h)=|\det d\phi_h|\,dA\quad\text{and}¸\quad(\phi_h|_\gamma)^{-1}(ds_h)=|\det d(\phi_h|_\gamma)|\,ds,
%\int_{\Gamma_h}f\,dA_h&=\int_{\Gamma}f\circ \phi_h\,|\det d\phi_h|\,dA,\\
%\int_{\phi_h(\gamma)}f\,ds_h&=\int_{\gamma}f\circ \phi_h\,|\det d(\phi_h|_{\gamma})|\,ds,
\end{aligned}
\end{equation}
where $ds$ and $ds_h$ denote the line elements of arbitrary curves $\gamma\subset\Gamma$
 and $\phi_h(\gamma)\subset\Gamma_h$, respectively; the determinants must be taken with respect to orthonormal bases in the respective tangent spaces. If we define $\tilde e$ as in the proof of Theorem \ref{s1h}, then we have
\begin{equation}\label{ket}
\begin{aligned}
|k_{\tilde e}|_{\tilde e}=|k_h|\circ\phi_h,\quad H_{\tilde e}=H_h\circ\phi_h, \quad\text{etc.} 
\end{aligned}
\end{equation}
From \eqref{intah}, \eqref{ket}, and \eqref{nablare} we easily obtain for $p\in [1,\infty]$
\begin{equation}\label{krh}
\begin{aligned}
|\phi_h^{-1}(dA_h)|&\le r_h |dA|\le r_h |\phi_h^{-1}(dA_h)|,\\
|(\phi_h|_\gamma)^{-1}(ds_h)|&\le r_h |ds|\le r_h |(\phi_h|_\gamma)^{-1}(ds_h)|,\\
\|k_h\|_p&\le r_h\big(1+|h|_{2,p}\big),\\
\|H_h\|_p&\le r_h\big(1+|h|_{2,p}\big),\\
\|\nabla H_h\|_p&\le r_h\big(1+|h|_{3,p}\big),\\
\|K_h\|_p&\le r_h\big(1+|h|_{2,2p}^2\big).
\end{aligned}
\end{equation}
For the fifth estimate we used the interpolation inequality \eqref{inth}.

\section{Function spaces}\applabel{app:spaces}
In this section we give the definition and some properties of the function spaces we are using in the present paper, namely the \emph{Sobolev-Slobodetskij spaces} $W^s_p$ and the \emph{Besov spaces} $B^s_{p,q}$. Let $e$ be a Riemannian metric on a smooth $d$-dimensional manifold $M$. For measurable tensor fields $T$  and $1\le p <\infty$ let
\[\|T\|_p^p:=\int_M |T|_e^p\,dV_e\quad\text{and}\quad \|T\|_\infty:=\esssup_M |T|_e,\]
where $dV_e$ is the volume element corresponding to $e$ and $|\cdot|_e$ is the norm induced by $e$ on all tensor bundles. We write $T\in L_p(M)$, $1\le p\le\infty$, if $\|T\|_p$ is finite. For $k\in\N_0$ and $1\le p\le \infty$ we write $T\in W^k_p(M)$ if $(\nabla^e)^l T\in L^p(M)$ for all $l\in\N_0$ with $l\le k$, and we define
\[|T|_{k,p}:=\|(\nabla^e)^kT\|_p.\]
Furthermore, for $s=k+\delta$ with $k\in\N_0$, $\delta\in (0,1)$ and $1\le p<\infty$ we write $T\in W^s_p(M)$ if $T\in W^k_p(M)$ and
\begin{equation*}%\label{deltanorm}
\begin{aligned}
|T|_{s,p}^p:=\int_M\int_M \frac{|(\nabla^e)^k T(x)-(\nabla^e)^k T(y)|^p_e}{d_{e}(x,y)^{d+\delta p}}\,dV_{e}(x)\,dV_{e}(y)<\infty,
\end{aligned}
\end{equation*}
where $d_e$ denotes the distance function corresponding to $e$. We will use the notation $H^s(M):=W^s_2(M)$. For arbitrary real $s>0$ we define
\[\|T\|_{s,p}:=\sum_{0\le l\le [s]} |T|_{l,p}+|T|_{s,p},\]
where $[s]$ is the largest integer smaller than $s$; note that this summation doesn't make sense from the point of view of physical dimensions. 
%Now, let $M=\R^d$, and let $\phi_j$, $j\in\N_0$, be a sequence of Schwartz functions such that $\phi_0$ is supported in the closed ball $\{x\in\R^d\,|\,|x|\le 2\}$, $\phi_j$ is supported in the the closed annulus $\{x\in\R^d\,|\,2^{j-1}\le |x|\le 2^{j+1}\}$, $\sum_{j}\phi_j=1$ in $\R^d$, and $2^{jk}|\nabla^k\phi_j|\le c_j$ in $\R^d$ for positive constants $c_j$. 
% \[\|T\|_{s,p,q}^q:=\|2^{sj}\|\mathcal{F}^{-1}\phi_j\mathcal{F}T\|_p\|_{l_q}<\infty,\]
% where $\mathcal{F}$ is the Fourier transform and $\|\cdot\|_{l_q}$ denotes the usual $l_q$-norm of sequences. It is well-known, see \cite{Triebel:Function-Spaces-1}, that for this norm is equivalent to
Now, let T be a measurable tensor field defined in $M=\R^d$, and consider the operator $\delta_hT:=T(\cdot+h)-T$. For $k\in\N$ and $1\le p,q\le\infty$ we write $T\in B^k_{p,q}(\R^d)$ if $T\in W^{k-1}_p(\R^d)$ and 
\[|T|_{k,p,q}:=\big\||h|^{-1}\|\delta_h^2\nabla^{k-1}T\|_p\big\|_{q,\frac{dh}{|h|^d}}<\infty,\]
where the $L_q$-norm is taken with respect to the measure $dh/|h|^d$. For $s=k+\delta$ with $k\in\N_0$, $\delta\in (0,1)$ and $1\le p,q\le\infty$ we write $T\in B^s_{p,q}(\R^d)$ if $T\in W^k_p(\R^d)$ and
\[|T|_{s,p,q}:=\big\||h|^{-\delta}\|\delta_h\nabla^kT\|_p\big\|_{q,\frac{dh}{|h|^d}}<\infty.\] 
For arbitrary real $s>0$ we define
\[\|T\|_{s,p,q}:=\sum_{0\le l\le [s]} |T|_{l,p}+|T|_{s,p,q}.\]
By the theorem of Section 2.5.12 in \cite{Triebel:Function-Spaces-1} this definition of the spaces $B^s_{p,q}(\R^d)$ is equivalent to the usual definition based on the Fourier transform and Littlewood-Paley theory; the latter can be extended to arbitrary $s\in\R$. For non-integer $s$ and $1\le p<\infty$ we have $B^s_{pp}(\R^d)=W^s_{p}(\R^d)$ algebraically and topologically, since the semi-norms $|\cdot|_{s,p}$ and $|\cdot|_{s,p,p}$ coincide in this case. By Proposition 2 (iii) of Section 2.3.2 and the theorem of Section 2.5.6 in \cite{Triebel:Function-Spaces-1} we have $B^k_{2,2}(\R^d)=H^k(\R^d)$ even for integer $k$.  Now, let $M$ be a smooth, closed manifold with some finite atlas $\psi_j$ and a subordinate partition of unity $\phi_j$. For real $s>0$, $1\le p,q\le\infty$, and measurable tensor fields $T$ defined on $M$ we write $T\in B^s_{p,q}(M)$ if $(\phi_j\,T)\circ\psi_j^{-1}\in B^s_{p,q}(\R^d)$ for all $j$, and we define
\[\|T\|_{s,p,q}:=\sum_j\|(\phi_j\,T)\circ\psi_j^{-1}\|_{s,p,q},\quad|T|_{s,p,q}:=\sum_j|(\phi_j\,T)\circ\psi_j^{-1}|_{s,p,q}.\]
It is not hard to see that the definition of the spaces $W^s_p(M)$ by this localization procedure is equivalent to the one given above; hence, we have $B^s_{pp}(M)=W^s_{p}(M)$ for non-integer $s$ and $1\le p<\infty$ and $B^k_{2,2}(M)=H^k(M)$ for integer $k$. Next, let $M$ be a smooth bounded domain in $\R^d$. For real $s>0$, $1\le p,q\le\infty$, and measurable tensor fields $T$ defined in $M$ we write $T\in B^s_{p,q}(M)$ if $T$ is the restriction of a tensor field $\tilde T\in B^s_{p,q}(\R^d)$. We define
\[\|T\|_{s,p,q}:=\inf_{\tilde T}\|\tilde T\|_{s,p,q},\]
where the infimum is taken over all extensions $\tilde T$ of $T$. By the theorem of Section 3.3.4 in \cite{Triebel:Function-Spaces-1} there exists a common linear and continuous extension operator $E:B^s_{p,q}(M)\rightarrow B^s_{p,q}(\R^d)$ for all $1\le p,q\le\infty$ and uniformly bounded real $s>0$. It follows from the theorem of Section 5.2.2 in \cite{Triebel:Function-Spaces-2} that the definition of the spaces $W^s_p(M)$ by this restriction approach is equivalent to the one given above; hence, again, we have $B^s_{pp}(M)=W^s_{p}(M)$ for non-integer $s$ and $1\le p<\infty$ and $B^k_{2,2}(M)=H^k(M)$ for integer $k$. We define
\[|T|_{s,p,q}:=|ET|_{s,p,q}\]
for some fixed instance of $E$. Furthermore, by Proposition 2 of Section 2.3.3 in \cite{Triebel:Function-Spaces-1} and by the construction of the spaces given above we have
\begin{equation}\label{qeinbettung}
B^s_{p,q_1}(M)\hookrightarrow B^s_{p,q_2}(M)
\end{equation}
if $q_1\le q_2$ for $M=\R^d$ and for $M$ a smooth, closed manifold or a smooth bounded domain.

\begin{lemma}
 Let $\Gamma$ be a closed surface contained in a smooth bounded domain $\Omega\subset\R^3$. For $s> 0$, $h\in H^{s+2}(\Gamma)$, $t>0$, and $2\le p\le\infty$ such that $\alpha:=(t-2/p-1)/s \in (0,1)$ we have
\begin{equation}\label{inth}
|h|_{t,p,2}\le c\big(|h|_{s+2,2}^{\alpha}\|h\|_{2,2}^{1-\alpha}+\|h\|_{2,2}\big).
\end{equation}
For $s> 0$, $q\in H^s(\Gamma)$, $0<t<s$, and $\alpha=t/s$ we have
\begin{equation}\label{q2}
 |q|_{t,2}\le c\big(|q|_{s,2}^{\alpha}\|q\|_{2}^{1-\alpha}+\|q\|_{2}\big).
\end{equation}
For $s\ge 0$, $u\in H^1_0(\Omega)$, $v=P_\Gamma[u]_\Gamma\in H^{s+1}(\Gamma)$, $t> 0$, and $2\le p\le\infty$ such that $\alpha:=(t-2/p+1/2)/(s+1/2)\in (0,1)$ we have
\begin{equation}\label{intv}
 |v|_{t,p,2}\le c\big(|v|_{s+1,2}^{\alpha}\|\nabla u\|_{2}^{1-\alpha}+\|\nabla u\|_{2}\big).
\end{equation}
For $s\ge 0$, $a\in H^{s}(\Omega)$, $t> 0$, and $2\le p\le\infty$ such that $\alpha:=(t-3/p+3/2)/s\in (0,1)$ we have
\begin{equation}\label{intubulk}
|a|_{t,p,2}\le c\big(|a|_{s,2}^{\alpha}\|a\|_{2}^{1-\alpha}+\|a\|_{2}\big).
\end{equation}
For $s\ge 0$, $u\in H^1_0(\Omega)\cap H^{s+1}(\Omega)$, $t> 0$, and $2\le p\le\infty$ such that $\alpha:=(t-2/p+1/2)/s\in (0,1)$ we have
\begin{equation}\label{intu}
 |[u]_\Gamma|_{t,p,2}\le c\big(|u|_{s+1,2}^{\alpha}\|\nabla u\|_{2}^{1-\alpha}+\|\nabla u\|_{2}\big).
\end{equation}
\end{lemma}
% \begin{lemma}
%  Let $\Gamma$ be a closed surface contained in a smooth bounded domain $\Omega\subset\R^3$. For $s>0$, $h\in H^{s+2}(\Gamma)$, $1<t < s+2$, $1\le p\le\infty$, and $\alpha:=(t-2/p-1)/s\in (0,1)$ we have
% \begin{equation}\label{inth}
% |h|_{t,p,2}\le c\big(|h|_{s+2,2}^{\alpha}\|\nabla h\|_{\infty}^{1-\alpha}+\|\nabla h\|_{\infty}\big).
% \end{equation}
% For $s>0$, $q\in H^s(\Gamma)$, $0<t<s$, and $\alpha=t/s$ we have
% \begin{equation}\label{q2}
%  |q|_{t,2}\le c\big(|q|_{s,2}^{\alpha}\|q\|_{2}^{1-\alpha}+\|q\|_{2}\big).
% \end{equation}
% For $s\ge 0$, $u\in H^1_0(\Omega)$, $v=P_\Gamma[u]_\Gamma\in H^{s+1}(\Gamma)$, $1/2<t<s+1$, $1\le p\le\infty$, and $\alpha:=(t-2/p+1/2)/(s+1/2)\in (0,1)$ we have
% \begin{equation}\label{intv}
%  |v|_{t,p,2}\le c\big(|v|_{s+1,2}^{\alpha}\|\nabla u\|_{2}^{1-\alpha}+\|\nabla u\|_{2}\big).
% \end{equation}
% For $s>0$, $a\in H^{s}(\Omega)$, $0<t<s$, $1\le p\le\infty$, and $\alpha:=(t-3/p+3/2)/s\in (0,1)$ we have
% \begin{equation}\label{intubulk}
% |a|_{t,p,2}\le c\big(|a|_{s,2}^{\alpha}\|a\|_{2}^{1-\alpha}+\|a\|_{2}\big).
% \end{equation}
% For $s>0$, $u\in H^1_0(\Omega)\cap H^{s+1}(\Omega)$, $1/2<t<s+1$, $1\le p\le\infty$, and $\alpha:=(t-2/p+1/2)/s\in (0,1)$ we have
% \begin{equation}\label{intu}
%  |[u]_\Gamma|_{t,p,2}\le c\big(|u|_{s+1,2}^{\alpha}\|\nabla u\|_{2}^{1-\alpha}+\|\nabla u\|_{2}\big).
% \end{equation}
% \end{lemma}
\proof
Let us begin with the proof of \eqref{inth}. From the theorem of Section 2.4.2 in \cite{Triebel:Function-Spaces-1} we have
\[(H^{s+2}(\R^2),H^2(\R^2))_{1-\alpha,2}=H^{s_\alpha}(\R^2),\]
where $s_\alpha=\alpha s +2$. Hence, using the embedding $H^{s_\alpha}(\R^2)\hookrightarrow B^t_{p,2}(\R^2)$ from the theorem in Section 2.7.1 in \cite{Triebel:Function-Spaces-1} with $t,p$ as in the assertion, we obtain
\[\|h\|_{t,p,2}\le c\,\|h\|_{s+2,2}^{\alpha}\|h\|_{2,2}^{1-\alpha}.\] 
% Fix some arbitrary numbers $\theta,\tilde\theta\in (0,1)$. By the theorem in section 2.7.1 in \cite{Triebel:Function-Spaces-1} we have $H^{s+1}(\R^2)=B^{s+1}_{2,2}(\R^2)\hookrightarrow B^{s}_{\infty,2}(\R^2)$, and Theorem 6.2.4 (9) in \cite{bergh76} gives us $L_\infty(\R^2)\hookrightarrow B^{0}_{\infty,\infty}(\R^2)$; for the definition of $B^0_{\infty,\infty}(\R^2)$ see \cite{Triebel:Function-Spaces-1}. Thus, by the theorem in section 2.4.2 in \cite{Triebel:Function-Spaces-1} we obtain
% \[(H^{s+1}(\R^2),L_{\infty}(\R^2))_{1-\theta,2} \hookrightarrow B^{\theta s}_{\infty,2}(\R^2).\]
% Interpolating this with the trivial embedding $H^{s+1}(\R^2)\hookrightarrow H^{s+1}(\R^2)$ yields
% \[\big(H^{s+1}(\R^2),(H^{s+1}(\R^2),L_{\infty}(\R^2))_{1-\theta,2}\big)_{1-\tilde\theta,2} \hookrightarrow (B^{s+1}_{2,2}(\R^2),B^{\theta s}_{\infty,2}(\R^2))_{1-\tilde\theta,2}.\]
% By Theorem 3.5.3 in \cite{bergh76} the left hand side equals \[(H^{s+1}(\R^2),L_{\infty}(\R^2))_{(1-\theta)(1-\tilde\theta),2},\] while Theorem 6.4.5 in \cite{bergh76} shows that the right hand side is nothing but $B^{t-1}_{p,2}(\R^2)$ with $t-1 = \tilde\theta(s+1) + (1-\tilde\theta)\theta s$ and $1/p=\tilde\theta/2$. Hence, choosing $\alpha=(1-\theta)(1-\tilde\theta)$ we proved that
% \[\|w\|_{t-1,p,2}\le c\,\|w\|_{s+1,2}^{\alpha}\|w\|_{\infty}^{1-\alpha}\]
% for all $w\in H^{s+1}(\R^2)$. 
In view of the definition of the Besov norms on $\Gamma$ the same estimate is true for all $h\in H^{s+2}(\Gamma)$; in fact, we even have \eqref{inth}.
% \begin{equation}\label{eqn:gnb}
% \|h\|_{t,p,2}\le c\big(|h|_{s+2,2}^{\alpha}\|h\|_{2,2}^{1-\alpha}+\|h\|_{2,2}\big). 
% \end{equation}
Indeed, if this estimate was false, there would exist a sequence $(h_n)\subset H^{s+2}(\Gamma)$ such that
\[|h_n|_{s+2,2}^{\alpha}\|h_n\|_{2,2}^{1-\alpha}+\|h_n\|_{2,2}<\frac1n\, \|h_n\|_{t,p,2}\le \frac cn\, \|h_n\|_{s+2,2}^{\alpha}\|h_n\|_{2,2}^{1-\alpha}.\]
Dividing by $\|h_n\|_{2,2}^{1-\alpha}$ and taking the $\alpha$-th root we obtain
\[|h_n|_{s+2,2}+\|h_n\|_{2,2}\le \Big(\frac cn\Big)^{\frac{1}{\alpha}}\, \|h_n\|_{s+2,2}.\]
We may assume without restriction that $\|h_n\|_{s+1,2}=1$. Hence, we have $h_n \rightarrow 0$ strongly in $H^2(\Gamma)$, weakly in $H^{s+2}(\Gamma)$, and thus strongly in $H^{[s+2]}(\Gamma)$ by the compact embedding $H^{s+2}(\Gamma)\hookrightarrow H^{[s+2]}(\Gamma)$ which follows from Remark 1 of Section 4.3.2 in \cite{Triebel:Function-Spaces-1}; here, $[s+2]$ is the largest integer smaller than $s+2$. Since also $|h_n|_{s+2,2}\rightarrow 0$ this is a contradiction to $\|h_n\|_{s+2,2}=1$. This proves \eqref{inth}. The proof of \eqref{q2} proceeds essentially along the same lines.
% Finally, the question arises if the set of numbers $(t,p)$ that can be reached by our construction involving $\theta,\tilde\theta$ coincides with the set of
% numbers $(t,p)$ from the assertion. This, however, is easily seen from Figure\marginpar{Figure!}. 
% From the theorem of section 2.4.2 in \cite{Triebel:Function-Spaces-1} we obtain
% \[(H^{s}(\R^2),L_2(\R^2))_{1-\alpha,2}=H^{s_\alpha}(\R^2),\]
% where $s_\alpha=\alpha s$; hence, 
% \[\|q\|_{t,2}\le c\,\|q\|_{s,2}^{\alpha}\|q\|_{2}^{1-\alpha}.\]
% From this, we can proceed as above to prove \eqref{q2}. 

The proof of \eqref{intv} is very similar using
\[(H^{s+1}(\R^2),H^{\frac12}(\R^2))_{1-\alpha,2}=H^{s_\alpha}(\R^2)\hookrightarrow B^t_{p,2}(\R^2)\]
with $s_\alpha=\alpha(s+1)+(1-\alpha)1/2$ and $t,p$ as in the assertion, the continuity of the trace operator $u\mapsto [u]_\Gamma$, $H^1(\Omega)\rightarrow H^{\frac12}(\Gamma)$ which follows from the theorem of Section 3.3.3 in \cite{Triebel:Function-Spaces-1}, and Poincar\'e's inequality. 

For the proof of \eqref{intubulk} we note that
\[(H^{s}(\R^3),L_2(\R^3))_{1-\alpha,2}=H^{s_\alpha}(\R^3)\hookrightarrow B^t_{p,2}(\R^3),\]
with $s_\alpha=\alpha s$ and $t,p$ as in the assertion, and we make use of the common extension operator $E:H^s(\Omega)\rightarrow H^s(\R^3)$, $E:L_2(\Omega)\rightarrow L_2(\R^3)$ introduced above. 

Finally, for the proof of \eqref{intu} we proceed very similarly by making use of
\[(H^{s+1}(\R^3),H^1(\R^3))_{1-\alpha,2}=H^{s_\alpha}(\R^3)\hookrightarrow B^{t+\frac1p}_{p,2}(\R^3),\]
where $s_\alpha=\alpha s+1$ and $t,p$ as in the assertion, the trace operator $u\mapsto [u]_\Gamma$, $B^{t+1/p}_{p,2}(\Omega)\rightarrow B^{t}_{p,2}(\Gamma)$, see again the theorem of Section 3.3.3 in \cite{Triebel:Function-Spaces-1}, and Poincar\'e's inequality.
\qed

\medskip
In the proof of Theorem \ref{s1h} we are considering perturbations of a some fixed metric $e$. Whenever we want to emphasize the dependence of the semi-norms on the metric we write ${_e|T|_{s,p}}$. Now, it is not hard to see that for $s>0$ and $1\le p <\infty$ we have
\begin{equation}\label{Tdelta}
  {_e|T|_{\delta,p}}\le r(\|\tilde e\|_\infty,\|\tilde e^{-1}\|_\infty,\|e\|_\infty,\|e^{-1}\|_\infty)\, {_{\tilde e}|T|_{\delta,p}},
\end{equation}
where the $L^\infty$-norms are taken with respect to a fixed background metric.

\begin{lemma}\label{lemma:prod}
 Let $M$ be a smoothly bounded domain or a smooth, closed manifold. Then, for sufficiently regular tensor fields $T_1,\ldots,T_k$ defined in $M$, and $\delta\in (0,1)$ we have
 \begin{equation*}
\begin{aligned}
\Big|\prod_{i=1}^kT_i\Big|_{\delta,2}\le c\sum_{j=1}^k\prod_{i=1\atop i\not=j}^k\|T_i\|_{p_{ij}}|T_j|_{\delta,p_{jj},2},
\end{aligned}
\end{equation*}
where $p_{ij}\in [2,\infty]$ such that $\sum_{i=1}^k1/p_{ij}=1/2$ for all $j$.
\end{lemma}
\proof It is sufficient to prove the estimate in the case $M=\R^d$. Indeed, then we deduce that for a closed manifold $M$
\begin{equation*}
\begin{aligned}
\Big|\prod_{i=1}^kT_i\Big|_{\delta,2}&\le \sum_{l}\Big|\Big(\phi_l\prod_{i=1}^k T_i\Big)\circ\psi_l^{-1}\Big|_{\delta,2}\\
&\le\sum_l\sum_{j=1}^k\prod_{i=1\atop i\not=j}^k\|T_i\circ\psi_l^{-1}\|_{p_{ij}}|(\phi_lT_j)\circ\psi_l^{-1}|_{\delta,p_{jj},2}\\
&\le\sum_{j=1}^k\prod_{i=1\atop i\not=j}^k\|T_i\|_{p_{ij}}|T_j|_{\delta,p_{jj},2},
\end{aligned}
\end{equation*}
while for a bounded domain $M$ we have
\begin{equation*}
\begin{aligned}
\Big|\prod_{i=1}^kT_i\Big|_{\delta,2}&\le \Big|\prod_{i=1}^k ET_i\Big|_{\delta,2}\le\sum_{j=1}^k\prod_{i=1\atop i\not=j}^k\|ET_i\|_{p_{ij}}|ET_j|_{\delta,p_{jj},2}\\
&\le\sum_{j=1}^k\prod_{i=1\atop i\not=j}^k\|T_i\|_{p_{ij}}|T_j|_{\delta,p_{jj},2},
\end{aligned}
\end{equation*}
where $E$ is the extension operator introduced above, and the first estimate is a trivial consequence of the definition of the $|\cdot|_{\delta,2}$ semi-norm. Now, let $M=\R^d$. We assume for simplicity that $k=2$; the general case follows by induction. We have
\begin{equation*}%\label{deltanorm}
\begin{aligned}
|T_1T_2|_{\delta,2}&=\big\||h|^{-\delta}\|\delta_h(T_1T_2)\|_2\big\|_{2,\frac{dh}{|h|^d}}\\
&\le \big\||h|^{-\delta}\|T_1\delta_h T_2\|_2\big\|_{2,\frac{dh}{|h|^d}}+\big\||h|^{-\delta}\|\delta_hT_1\,T_2(\cdot+h)\|_2\big\|_{2,\frac{dh}{|h|^d}}\\
&\le \|T_1\|_{p_{12}}\big\||h|^{-\delta}\|\delta_h T_2\|_{p_{22}}\big\|_{2,\frac{dh}{|h|^d}}+\|T_2\|_{p_{21}}\big\||h|^{-\delta}\|\delta_hT_1\|_{p_{11}}\big\|_{2,\frac{dh}{|h|^d}}\\
&= \|T_1\|_{q_{12}}|T_2|_{\delta,q_{22},2} + \|T_2\|_{q_{21}}|T_1|_{\delta,q_{11},2}.
\end{aligned}
\end{equation*}
Here, we used the triangle inequality for the first estimate and H\"older's inequality for the second estimate.
% \begin{equation*}%\label{deltanorm}
% \begin{aligned}
% |f_1f_2|_{\delta,2}^2&=\int_{\R^d}\int_{\R^d} |(f_1f_2)(x)-(f_1f_2)(x-h)|^2\,dx\frac{1}{|h|^{d+2\delta}}dh\\
% &\le 2\int_{\R^d}\int_{\R^d} |f_1(x)|^2|f_2(x)-f_2(x-h)|^2\,dx\frac{1}{|h|^{d+2\delta}}dh\\
% &\quad+2\int_{\R^d}\int_{\R^d} |f_1(x)-f_1(x-h)|^2|f_2(x-h)|^2\,dx\frac{1}{|h|^{d+2\delta}}dh\\
% &\le 2\|f_1\|_{2\tilde q_{12}}^2\int_{\R^d}\|f_2(x)-f_2(x-h)\|_{2\tilde q_{22}}^2\frac{1}{|h|^{d+2\delta}}dh\\
% &\quad +2\|f_2\|_{2\tilde q_{21}}^2\int_{\R^d}\|f_1(x)-f_1(x-h)\|_{2\tilde q_{11}}^2\frac{1}{|h|^{d+2\delta}}dh\\
% &= 2\|f_1\|_{2\tilde q_{12}}^2|f_2|_{\delta,2\tilde q_{22},2}^2 + 2\|f_2\|_{2\tilde q_{21}}^2|f_1|_{\delta,2\tilde q_{11},2}^2,
% \end{aligned}
% \end{equation*}
%where we used H\"older's inequality for the second estimate with $\sum_i 1/\tilde q_{ij}=1$.
\qed

\medskip
\begin{corollary}\label{cor:analytic}
 Let $M$ be a smoothly bounded domain or a smooth, closed manifold. Then, for a sufficiently regular tensor field $T$ of class $(k,l)$ defined in $M$, an analytic bundle homomorphism $f$ mapping the tensor bundle of class $(k,l)$ to some other tensor bundle, and $\delta\in (0,1)$ we have
\begin{equation*}
\begin{aligned}
|f(T)|_{\delta,2}\le \tilde f(\|T\|_{\infty})|T|_{\delta,2},
\end{aligned}
\end{equation*}
where $\tilde f$ is an increasing analytic function depending on $f$.
\end{corollary}

\bibliographystyle{plain}
\bibliography{references}
\end{document}